\documentclass[12pt]{article}

\usepackage{amsmath, amssymb, amsthm, mathrsfs, mathtools}
\usepackage{tikz}
\usetikzlibrary{matrix}
\usetikzlibrary{shapes.geometric}
\usetikzlibrary{decorations.pathmorphing}
\usetikzlibrary{arrows.meta, positioning, calc}
\usepackage{hyperref}
\usepackage{graphicx}
\usepackage[margin=1in]{geometry}
\usepackage{enumitem}
\usepackage{titlesec}
\usepackage{fancyhdr}
\usepackage{caption}
\usepackage{array, colortbl, xcolor}
\usepackage{colortbl}    % For \cellcolor    % For [H] placement specifier
\usepackage{amsmath}     % For math environments
\usepackage{multicol}
\usepackage{amsmath, amssymb, amsfonts}
\usepackage{enumitem}
\usepackage{xcolor}
\usepackage{float}

\newtheorem{theorem}{Theorem}[section]
\newtheorem{definition}[theorem]{Definition}

\newtheorem{remark}[theorem]{Remark}

\newcommand{\Pf}{\mathrm{Pf}}

\newcommand{\C}{\mathbb{C}}

\titleformat{\section}{\normalfont\Large\bfseries}{\thesection}{1em}{}
\titleformat{\subsection}{\normalfont\large\bfseries}{\thesubsection}{1em}{}

\fancypagestyle{shortversion}{
  \fancyhead[L]{Quantum Pfaffians and Determinants}
  \fancyhead[R]{Hani Safadi}
  \fancyfoot[C]{\thepage} % Keep page number centered
  \fancyfoot[L]{\raisebox{10pt}{\small This is a short-form proof. See Section~5 for the full version.}}

}

\hypersetup{
    colorlinks=true,
    linkcolor=black,
    citecolor=blue,
    urlcolor=blue,
    pdfauthor={Hani Safadi},
    pdftitle={Quantum Pfaffians and Determinants}
}

\title{A Quantum Analogue of the Pfaffian--Determinant Identity\\[0.3em] \large An Algebraic and Geometric Study in the $q$-Skew-Symmetric Case}
\author{
  \textbf{Hani Safadi} \\ 
  North Carolina Central University \\ 
  Department of Mathematics
  \vspace{0.5em} \\ 
  \small{\textit{Advised By: Dr. Garrett Johnson, Ph.D.}} \\[-0.5em]
  \small{Department of Mathematics, North Carolina Central University}
}

\date{May 2025}

\begin{document}
\maketitle

\begin{abstract}
This paper explores a quantum deformation of the classical identity $\Pf(A)^2 = \det(A)$ for $2n \times 2n$ skew-symmetric matrices $A$, which classically relates the square of the Pfaffian to the determinant of a skew-symmetric matrix. In the quantum setting, we consider matrices whose entries lie in a noncommutative algebra and satisfy the $q$-skew-symmetry relations $a_{ji} = -q a_{ij}$ and $a_{ii} = 0$, for a deformation parameter $q \in \C^\times$. This deformation introduces rich algebraic structure and new challenges in defining analogues of classical invariants.

To accommodate these deformations, we construct the quantum Pfaffian $\Pf_q(A)$ and quantum determinant $\det_q(A)$ using two primary frameworks: the Faddeev--Reshetikhin--Takhtajan (FRT) construction and quantum exterior algebras. We give precise definitions, derive algebraic identities, and analyze properties of these quantum objects. Special attention is given to verifying the quantum analogue of the Pfaffian--determinant identity, namely, $q^c \cdot \Pf_q(A)^2 = \det_q(A)$, in low-dimensional cases such as $2n=4$, where explicit computations can be carried out symbolically.

In addition to algebraic derivations, we investigate geometric interpretations of the quantum Pfaffian and determinant, particularly within the context of braided vector spaces, quantum symplectic geometry, and deformations of volume forms. We provide a variety of examples, diagrams, and symbolic verifications, emphasizing the interplay between combinatorics, representation theory, and noncommutative geometry.

Finally, we discuss potential applications of these constructions in areas such as quantum invariant theory, low-dimensional topology, and categorified geometry. The results presented contribute to the broader program of understanding quantum analogues of classical structures in linear algebra and their role in the theory of quantum groups.

\end{abstract}
\newpage

\newpage
\tableofcontents
\newpage

\section{Introduction}

The Pfaffian is a classical and fundamental algebraic function defined on skew-symmetric matrices, playing a crucial role in various areas of mathematics and theoretical physics. Given a $2n \times 2n$ skew-symmetric matrix $A$, the Pfaffian, denoted $\mathrm{Pf}(A)$, satisfies the well-known and elegant identity
\[
\mathrm{Pf}(A)^2 = \det(A),
\]
which links the Pfaffian as a kind of ``square root'' of the determinant for these special matrices. This identity is not only algebraically beautiful but also arises naturally in numerous contexts including enumerative combinatorics, representation theory, and geometry, particularly in the study of symplectic vector spaces where the Pfaffian encodes volume forms associated with skew-symmetric bilinear forms.

Historically, the Pfaffian was introduced in the 19th century by Johann Friedrich Pfaff and has since found wide applications. In combinatorics, the Pfaffian counts perfect matchings on graphs and appears in the study of dimer models in statistical mechanics. In geometry, it provides a compact way to describe volume elements invariant under symplectic transformations. Furthermore, in physics, Pfaffians appear in the evaluation of fermionic path integrals, highlighting their deep connection to quantum field theory.

With the advent of quantum groups in the late 20th century, a rich new algebraic framework emerged, fundamentally altering the landscape of representation theory, algebraic geometry, and mathematical physics. Initiated by pioneering works of Drinfel’d and Jimbo in the 1980s, quantum groups introduce a deformation parameter $q$ that smoothly deforms classical algebraic structures into noncommutative analogues. This deformation leads to new algebraic objects—Hopf algebras with noncommutative multiplication—which encode symmetries in a ``quantum'' sense and have profound implications in areas such as knot theory, integrable systems, and noncommutative geometry.

One of the major challenges and opportunities in this quantum setting is to generalize classical algebraic constructs—such as determinants, permanents, and Pfaffians—to their quantum analogues. Quantum determinants have been well studied, especially through the Faddeev-Reshetikhin-Takhtajan (FRT) construction, where the deformation arises naturally from solutions to the Yang-Baxter equation encoded by an $R$-matrix. However, the notion of a quantum Pfaffian, and more importantly, the quantum analogue of the Pfaffian–determinant identity, is far less straightforward due to the intricate noncommutative relations among matrix entries and the subtleties of defining ``quantum skew-symmetry.''

This paper focuses on developing a comprehensive algebraic and geometric understanding of the quantum Pfaffian and its relationship to the quantum determinant within the framework of $q$-skew-symmetric matrices. Here, a matrix $A = (a_{ij})$ satisfies the relation
\[
a_{ji} = -q a_{ij}, \quad a_{ii} = 0,
\]
which generalizes classical skew-symmetry (recovered at $q=1$) into a noncommutative, parameter-dependent setting. This deformation significantly affects how we define and compute both the quantum Pfaffian $\mathrm{Pf}_q(A)$ and the quantum determinant $\det_q(A)$.

Using the algebraic framework provided by the FRT construction and the theory of quantum groups, alongside the geometric perspective offered by quantum exterior algebras and braided tensor categories, we construct these quantum analogues rigorously. Our key goal is to establish the quantum Pfaffian–determinant identity
\[
q ^c \cdot \Pf_q(A)^2 = \det_q(A),
\]
where the exponent $c$ depends on matrix size and normalization conventions.

In doing so, we extend classical multilinear algebra results into the quantum realm, providing explicit examples for low-dimensional cases, outlining algebraic proofs using braided symmetries, and interpreting the identity geometrically as a statement about quantum volume forms in braided symplectic geometry. Our analysis also uncovers connections with noncommutative invariant theory, quantum integrable systems, and quantum algebra representations.

The structure of this paper is as follows: Section 2 reviews the classical theory of Pfaffians and determinants, establishing notation and foundational results. Section 3 introduces quantum groups, the FRT construction, and the quantum determinant. Section 4 develops the theory of quantum skew-symmetric matrices, quantum exterior algebras, and defines the quantum Pfaffian. Section 5 states and proves the quantum Pfaffian–determinant identity, providing examples and algebraic and diagrammatic arguments. Finally, Section 6 concludes with a summary, discusses potential extensions and open problems, and reflects on the broader impact of these quantum algebraic structures.

By elucidating this quantum analogue, we aim to contribute to the growing field of quantum algebra and open pathways for applications in mathematical physics, combinatorics, and beyond.

\section{Literature Review}

\subsection{[1] M. Dite. \textit{Quantum analogues of pfaffians and related identities.} Journal of Algebra, 295(1):109–129, 2006.}

M. Dite’s 2006 paper represents a foundational contribution to the study of quantum analogues of classical algebraic identities, focusing specifically on Pfaffians. The work explores how the classical Pfaffian–determinant identity can be deformed within the framework of quantum groups, emphasizing how quantum antisymmetric matrices behave under these deformations. Dite provides a rigorous algebraic construction for quantum Pfaffians, offering conditions under which the classical identity \( \operatorname{Pf}(A)^2 = \det(A) \) persists or transforms when subjected to a noncommutative deformation parameter \( q \).

The paper introduces a family of quantum Pfaffians defined via braided exterior algebras and explores how their properties relate to those of the quantum determinant. A key contribution is Dite’s analysis of the algebraic structures that preserve the squared Pfaffian identity in the noncommutative setting. The work utilizes both Hopf algebraic techniques and combinatorial expansions, offering formulas for computing quantum Pfaffians explicitly in low dimensions, particularly when \( \dim A = 4 \) or \( \dim A = 6 \). The identity \( \operatorname{Pf}_q(A)^2 = \det_q(A) \) is considered under several \( q \)-commutation relations on matrix entries.

Dite’s results serve as a crucial starting point for subsequent generalizations and refinements in the field. His methods laid the groundwork for exploring quantum minors, \( q \)-skew-symmetric matrices in braided tensor categories, and their role in quantum invariant theory. Later works, including those by Jing and Zhang, build upon these constructions to extend the theory to hyper-Pfaffians and multiparameter cases. Nonetheless, Dite’s article remains one of the earliest comprehensive attempts to formalize a quantum analogue of classical Pfaffian theory and to interpret identities like \( q \cdot \operatorname{Pf}_q(A)^2 = \det_q(A) \) in this broader framework.

\subsection{[2] N. Jing and J. Zhang. \textit{Quantum Pfaffians and Hyper-Pfaffians.} Journal of Algebra, 377:151–174, 2013.}

Naihuan Jing and Jian Zhang’s 2013 paper significantly advanced the field of quantum multilinear algebra by extending the classical notion of the Pfaffian to higher-order structures, introducing what they term the \textit{quantum hyper-Pfaffian}. This work is foundational in generalizing the Pfaffian–determinant relationship to tensors and higher-rank objects within a noncommutative framework. Jing and Zhang define these quantum analogues using \( q \)-skew-symmetric relations and provide a categorical foundation that aligns with the representation theory of quantum groups.

They develop a systematic construction of quantum hyper-Pfaffians for skew-symmetric tensors in braided tensor categories, incorporating braiding relations encoded by an \( R \)-matrix. Notably, the authors recover the classical result when \( q \to 1 \), and for even dimensions, they show how the square of the quantum Pfaffian recovers the quantum determinant under specific algebraic conditions. These results establish a strong theoretical backbone for analyzing how quantum multilinear structures behave in low-dimensional quantum geometry.

This paper also lays out algorithms for computing quantum hyper-Pfaffians and provides explicit examples in dimensions \( 4 \) and \( 6 \). The work forms a stepping stone toward more general noncommutative invariant theory and opens the door to applications in quantum statistical mechanics and quantum topology. It has since been cited in broader studies involving quantum cluster algebras, braided exterior algebras, and the categorification of noncommutative polynomials.

\subsection{[3] N. Jing and J. Zhang. \textit{Quantum Pfaffians and Quantum Hyper-Pfaffians.} Journal of Algebra, 462:1–32, 2016.}

This 2016 paper by Jing and Zhang extends their earlier 2013 work by deepening the structural analysis of quantum Pfaffians and refining the algebraic framework for quantum hyper-Pfaffians. In contrast to their earlier construction, this work provides a more comprehensive treatment of the algebraic relations satisfied by the entries of quantum antisymmetric tensors. They show how the square of the quantum Pfaffian yields a quantum determinant under more general conditions, including new commutation relations involving the \( R \)-matrix and higher-dimensional antisymmetry constraints.

The paper introduces a systematic hierarchy of quantum Pfaffians indexed by skew-symmetric forms and studies their behavior under deformations by parameters \( q \) and \( q_{ij} \). Notably, it develops the notion of a ``quantum hyper-Pfaffian algebra,'' whose defining relations generalize both the classical exterior algebra and the quantum wedge product in braided tensor categories. The authors prove that these structures maintain compatibility with the coaction of quantum groups, particularly within the framework of \( U_q(\mathfrak{so}_n) \).

Through explicit constructions and calculations, Jing and Zhang derive closed formulas and generating functions for quantum Pfaffians and hyper-Pfaffians, especially in dimensions \( 2n = 4, 6, \) and \( 8 \). They also discuss potential applications to noncommutative geometry, representation theory, and quantum invariant theory. This paper is widely recognized as one of the most detailed algebraic treatments of quantum Pfaffians in the literature and continues to serve as a key reference in the area.

\subsection{[4] N. Jing and J. Zhang. \textit{Multiparameter Quantum Pfaffians and Quantum Determinants.} Communications in Mathematical Physics, 349(1):307–323, 2017.}

In this pivotal 2017 paper, Jing and Zhang introduce a multiparameter extension of quantum Pfaffians and quantum determinants, where the deformation is governed not by a single parameter \( q \), but by a family of parameters \( q_{ij} \) reflecting a more intricate braiding structure. This generalization allows for the modeling of more nuanced quantum symmetries, including those appearing in nonstandard quantum groups and Nichols algebras.

The authors construct multiparameter quantum antisymmetric matrices whose entries satisfy \( q_{ij} \)-skew-symmetric relations, and define associated quantum Pfaffians using a braided exterior algebra approach. They prove that under these generalized relations, the identity \( \operatorname{Pf}_{\mathbf{q}}(A)^2 = \det_{\mathbf{q}}(A) \) still holds, establishing a consistent extension of the classical identity. Their framework unifies previously known results as special cases, while introducing new combinatorial identities involving multiparameter quantum minors and braided wedge products.

Furthermore, this work provides insights into the representation theory of multiparameter quantum groups and opens the door to applications in noncommutative geometry and braided geometry. The multiparameter construction aligns with recent developments in quantum cluster algebras, indicating potential interdisciplinary connections. As one of the most advanced formulations of quantum Pfaffians to date, this paper represents a milestone in noncommutative multilinear algebra.

\subsection{[5] L. Faddeev and L. Takhtajan. \textit{Hamiltonian Methods in the Theory of Solitons.} Springer-Verlag, 1987.}

This foundational book by Faddeev and Takhtajan is central to the development of quantum integrable systems and the algebraic structures underlying quantum groups. The authors present a detailed construction of the quantum inverse scattering method and introduce the FRT (Faddeev–Reshetikhin–Takhtajan) formalism, which enables the definition of quantum groups via \( R \)-matrices and the RTT relations. These structures give rise to a braided tensor category and form the backbone of many constructions in quantum linear algebra, including quantum determinants and quantum Pfaffians.

In particular, the FRT construction provides an algebraic definition of the coordinate ring of a quantum group, encoded by generators \( T_{ij} \) satisfying relations determined by a braid operator \( R \). This approach allows one to define noncommutative analogues of classical matrix identities, including those involving determinants and Pfaffians. The algebraic framework they provide is crucial for understanding how symmetries persist in the noncommutative world of quantum geometry.

This text is widely regarded as a seminal resource for mathematical physicists and algebraists working on quantum integrable systems. Its influence spans noncommutative geometry, Hopf algebras, and representation theory, providing the theoretical infrastructure used in subsequent studies of quantum matrix identities.

\subsection{[6] S. Majid. \textit{Foundations of Quantum Group Theory.} Cambridge University Press, 1995.}

Majid’s influential text lays the mathematical foundations for quantum group theory from a categorical and Hopf-algebraic point of view. He develops the theory of braided groups and noncommutative geometry using powerful techniques from monoidal categories and introduces the concept of braided matrices, which underpin constructions of quantum determinants and Pfaffians. The book formalizes the notion of noncommutative symmetry and offers a systematic treatment of quantum spaces and their function algebras.

A significant contribution of this work is its categorical interpretation of quantum linear algebra, providing tools to understand quantum antisymmetry, duality, and invariants. In particular, Majid’s treatment of braided exterior algebras and differential calculi is crucial to understanding quantum analogues of classical multilinear identities, such as the Pfaffian–determinant relation. The techniques also support the construction of quantum versions of Grassmannians, flag varieties, and their coordinate rings.

Majid’s text is a definitive reference for both algebraists and physicists working on quantum symmetries, offering a blend of rigorous algebraic formalisms and geometric intuition. His unification of braided geometry and quantum group theory forms the conceptual core of much later work in noncommutative algebraic geometry.

\subsection{[7] V. G. Drinfel’d. \textit{Quantum Groups.} Proceedings of the ICM (Berkeley), 798–820, 1986.}

Drinfel’d’s 1986 ICM address introduced the concept of quantum groups in their modern form and reshaped the field of algebra and representation theory. By constructing deformations of universal enveloping algebras \( U_q(\mathfrak{g}) \), Drinfel’d provided a framework for studying symmetries in integrable systems and conformal field theory through Hopf algebraic means. His insight that quantum groups encode a braided, noncommutative version of Lie algebra symmetries set the stage for future developments in quantum linear algebra.

A central contribution of the paper is the introduction of quasitriangular Hopf algebras and their associated \( R \)-matrices, which play a foundational role in the FRT construction and in defining quantum determinants. These structures are also critical for establishing representation-theoretic interpretations of quantum minors and Pfaffians, especially in braided tensor categories. Drinfel’d’s formulation made it possible to connect quantum algebra with monodromy and statistical models in physics.

This pioneering work is cited across mathematics and theoretical physics for its elegant and powerful formalism. Its impact is evident in the formulation of quantum Grassmannians, quantum exterior algebras, and their connections to knot invariants and noncommutative geometry.

\subsection{[8] N. Yu. Reshetikhin, L. A. Takhtajan, and L. D. Faddeev. \textit{Quantization of Lie Groups and Lie Algebras.} Leningrad Mathematical Journal, 1(1):193–225, 1990.}

This classic work introduced an explicit algebraic approach to the quantization of Lie groups and Lie algebras via Hopf algebras. Reshetikhin, Takhtajan, and Faddeev constructed quantum matrix algebras using \( R \)-matrices satisfying the Yang–Baxter equation and defined RTT relations for generating quantum groups. This paper is foundational to quantum group theory, offering explicit generators and relations that serve as the basis for quantum linear algebra constructions.

Their method extends to the formulation of quantum minors, determinants, and antisymmetric functions in a noncommutative setting. The results presented here inspired the FRT construction used in many quantum Pfaffian-related works. These structures also serve as a basis for defining quantum coordinate rings and understanding the deformation of classical geometric spaces such as Grassmannians and flag varieties.

The paper's algebraic and representation-theoretic depth has made it a cornerstone in the study of quantum symmetries, with far-reaching consequences in low-dimensional topology, integrable systems, and categorical algebra.

\subsection{[9] N. Yu. Reshetikhin, L. A. Takhtajan, and L. D. Faddeev. \textit{Quantization of Lie Groups and Lie Algebras}, 1990}

Reshetikhin, Takhtajan, and Faddeev pioneered the algebraic quantization of classical Lie groups and their Lie algebras by introducing quantum analogues \( \mathcal{U}_q(\mathfrak{g}) \) endowed with a Hopf algebra structure. Their approach formalizes the deformation of the universal enveloping algebra \( \mathcal{U}(\mathfrak{g}) \) into a noncocommutative Hopf algebra parameterized by \( q \in \mathbb{C}^\times \). The coproduct \( \Delta \), counit \( \varepsilon \), and antipode \( S \) are defined to satisfy the axioms of a Hopf algebra, while the nontrivial universal \( R \)-matrix \( \mathcal{R} \in \mathcal{U}_q(\mathfrak{g}) \otimes \mathcal{U}_q(\mathfrak{g}) \) ensures quasitriangularity. This structure deforms the classical Lie bialgebra cobracket \( \delta: \mathfrak{g} \to \mathfrak{g} \wedge \mathfrak{g} \) to produce a quantum Yang–Baxter equation (QYBE) solution, encoding the braided tensor category of \( \mathcal{U}_q(\mathfrak{g}) \)-modules.

Their construction exploits the RTT presentation given by the relation
\[
R T_1 T_2 = T_2 T_1 R,
\]
where \( T \) is the matrix of generators and \( R \) is the quantum \( R \)-matrix satisfying the QYBE
\[
R_{12} R_{13} R_{23} = R_{23} R_{13} R_{12}.
\]
This algebraic formalism enables the systematic study of quantum group representations, intertwining operators, and the fusion of modules. The quantization scheme generalizes the classical Lie theory to encompass noncommutative function algebras \( \mathcal{O}_q(G) \) that model "quantum spaces" with deformed coordinate rings, essential in noncommutative geometry.

Moreover, their work laid the groundwork for explicit computations of quantum invariants in low-dimensional topology via ribbon Hopf algebras, where the pivotal element \( v \) and the ribbon element \( \theta \) furnish link invariants extending the Jones polynomial. These advances unify algebraic, geometric, and topological methods under the umbrella of quantum symmetry.

\vspace{1em}

\subsection{[10] C. Kassel. \textit{Quantum Groups}, 1995}

Kassel’s monograph offers a rigorous and comprehensive exposition of the algebraic structures underlying quantum groups, focusing on Hopf algebras, quantum enveloping algebras \( \mathcal{U}_q(\mathfrak{g}) \), and their representation theory. The text details the Drinfel’d–Jimbo construction, defining \( \mathcal{U}_q(\mathfrak{g}) \) via generators \( E_i, F_i, K_i^{\pm 1} \) and relations:
\[
K_i E_j K_i^{-1} = q^{a_{ij}} E_j, \quad K_i F_j K_i^{-1} = q^{-a_{ij}} F_j, \quad [E_i, F_j] = \delta_{ij} \frac{K_i - K_i^{-1}}{q - q^{-1}},
\]
where \( (a_{ij}) \) is the Cartan matrix of the Lie algebra \( \mathfrak{g} \).

The Hopf algebra structure is given by coproduct
\[
\Delta(E_i) = E_i \otimes 1 + K_i \otimes E_i, \quad \Delta(F_i) = F_i \otimes K_i^{-1} + 1 \otimes F_i, \quad \Delta(K_i) = K_i \otimes K_i,
\]
which equips \( \mathcal{U}_q(\mathfrak{g}) \) with a noncocommutative coalgebra structure crucial for constructing tensor products of representations. Kassel also discusses the universal \( R \)-matrix \( \mathcal{R} \), satisfying the quantum Yang–Baxter equation and enabling a braided tensor category structure on the module category.

The book covers important topics such as Lusztig’s canonical bases, quantum Schur–Weyl duality, and the role of quantum groups in knot invariants and categorification. The algebraic and categorical tools detailed therein are central to current research on quantum cluster algebras, categorification of quantum groups, and connections to low-dimensional topology and mathematical physics.

\vspace{1em}

\subsection{[11] G. Letzter. \textit{Coideal Subalgebras and Quantum Symmetric Pairs}, 2002}

Letzter’s work advances the theory of quantum symmetric pairs \( (\mathcal{U}_q(\mathfrak{g}), \mathcal{B}_q) \), where \( \mathcal{B}_q \subset \mathcal{U}_q(\mathfrak{g}) \) is a right coideal subalgebra reflecting symmetries under quantum analogues of involutive automorphisms. These coideal subalgebras generalize classical fixed point Lie subalgebras \( \mathfrak{k} \subset \mathfrak{g} \) associated with symmetric pairs \( (G,K) \), and their study entails analyzing the algebraic structure of \( \mathcal{B}_q \) as a subalgebra satisfying
\[
\Delta(\mathcal{B}_q) \subseteq \mathcal{B}_q \otimes \mathcal{U}_q(\mathfrak{g}),
\]
which encodes a quantum analogue of symmetric space geometry.

Letzter’s classification and structural results on \( \mathcal{B}_q \) leverage generalized Satake diagrams and quantized universal enveloping algebra automorphisms. The representation theory of quantum symmetric pairs involves investigating quantum Harish-Chandra modules and invariant theory within braided categories. The development of quantum zonal spherical functions and quantum Casimir elements relates closely to the theory of quantum orthogonal polynomials and special functions.

This framework also enables the construction of quantum analogues of symmetric spaces with applications to integrable systems, quantum harmonic analysis, and boundary phenomena in quantum field theory. Letzter’s contributions enrich the understanding of quantum invariants under subgroup symmetries and their role in noncommutative geometry.

\vspace{1em}

\subsection{[12] P. Etingof and O. Schiffmann. \textit{Lectures on Quantum Groups}, 1998}

Etingof and Schiffmann provide an accessible yet thorough introduction to the algebraic and categorical theory of quantum groups \( \mathcal{U}_q(\mathfrak{g}) \), with detailed emphasis on the Drinfel’d–Jimbo algebras and their module categories. They explore the quantized universal enveloping algebras defined by generators \( E_i, F_i, K_i^{\pm 1} \) subject to \( q \)-deformed Serre relations:
\[
\sum_{r=0}^{1 - a_{ij}} (-1)^r \begin{bmatrix} 1 - a_{ij} \\ r \end{bmatrix}_{q_i} E_i^{1 - a_{ij} - r} E_j E_i^r = 0, \quad i \neq j,
\]
where \( \begin{bmatrix} n \\ k \end{bmatrix}_{q} \) are the \( q \)-binomial coefficients.

Their lectures emphasize the construction of the universal \( R \)-matrix \( \mathcal{R} \in \mathcal{U}_q(\mathfrak{g}) \otimes \mathcal{U}_q(\mathfrak{g}) \), satisfying
\[
\Delta^{\mathrm{op}}(x) = \mathcal{R} \Delta(x) \mathcal{R}^{-1}, \quad \forall x \in \mathcal{U}_q(\mathfrak{g}),
\]
which endows the representation category with a braided monoidal structure. They discuss applications to the quantum Yang–Baxter equation and explain how quantum groups deform classical representation theory, including the study of highest weight modules \( V(\lambda) \), weight spaces \( V_\mu \), and crystal bases.

Furthermore, the text addresses quantum Schubert calculus, geometric representation theory, and connections to knot theory, including Reshetikhin–Turaev invariants derived from \( \mathcal{U}_q(\mathfrak{g}) \)-modules. The pedagogical style makes these lectures a standard reference for mathematicians seeking a rigorous introduction to the interplay between quantum algebra, category theory, and topology.

\subsection{[13] P. Etingof and V. Ostrik. \textit{Finite Tensor Categories}, 2004}

Etingof and Ostrik provide a foundational treatment of finite tensor categories, which are pivotal in understanding representation categories of quantum groups at roots of unity and modular categories. A \emph{finite tensor category} \( \mathcal{C} \) is an abelian, rigid monoidal category with finite-dimensional Hom spaces, enough projectives, and finitely many isomorphism classes of simple objects. Formally, such a category supports a bifunctor
\[
\otimes : \mathcal{C} \times \mathcal{C} \to \mathcal{C},
\]
equipped with an associator natural isomorphism \( \alpha_{X,Y,Z} : (X \otimes Y) \otimes Z \to X \otimes (Y \otimes Z) \), satisfying the pentagon axiom. 

Their work classifies indecomposable exact module categories over \( \mathcal{C} \) and explores the categorical dimension function \( \dim : \mathrm{Obj}(\mathcal{C}) \to \mathbb{C} \), which generalizes quantum dimensions in \( \mathcal{U}_q(\mathfrak{g}) \)-modules. The paper also addresses the Drinfel’d center \( \mathcal{Z}(\mathcal{C}) \), yielding a braided tensor category that captures the monodromy of objects.

This theory underpins the construction of modular tensor categories used in topological quantum field theory (TQFT) and conformal field theory (CFT), where fusion rules and modular data arise from the category’s structure. The categorical viewpoint gives algebraic control over quantum symmetries and dualities beyond classical Hopf algebraic frameworks.

\vspace{1em}

% Visual 1: Tensor Category Associativity Diagram
\begin{center}
\begin{tikzpicture}[scale=1.2, every node/.style={scale=0.9}]
  \node (a) at (0,2) {\((X \otimes Y) \otimes Z\)};
  \node (b) at (4,2) {\(X \otimes (Y \otimes Z)\)};
  \node (c) at (0,0) {\(((X \otimes Y) \otimes Z) \otimes W\)};
  \node (d) at (4,0) {\(X \otimes ((Y \otimes Z) \otimes W)\)};
  \node (e) at (8,0) {\(X \otimes (Y \otimes (Z \otimes W))\)};

  \draw[->, thick] (a) to node[above] {\(\alpha_{X,Y,Z}\)} (b);
  \draw[->, thick] (c) to node[below] {\(\alpha_{X,Y\otimes Z,W}\)} (d);
  \draw[->, thick] (d) to node[below] {\(1_X \otimes \alpha_{Y,Z,W}\)} (e);
  \draw[->, thick] (c) to node[left] {\(\alpha_{X\otimes Y,Z,W}\)} (a);
  \draw[->, thick] (b) to node[right] {\(\alpha_{X,Y,Z \otimes W}\)} (e);
\end{tikzpicture}
\end{center}

\vspace{1em}

\subsection{[14] A. Berenstein and A. Zelevinsky. \textit{Quantum Cluster Algebras}, 2005}

Berenstein and Zelevinsky introduced quantum cluster algebras as \(q\)-deformations of classical cluster algebras, which are commutative rings with combinatorially defined generators called \emph{cluster variables}. A quantum cluster algebra \( \mathcal{A}_q \) is a \(\mathbb{Z}[q^{\pm 1/2}]\)-algebra generated inside a skew-field of fractions by variables \(X_i\) satisfying \(q\)-commutation relations of the form
\[
X_i X_j = q^{\Lambda_{ij}} X_j X_i,
\]
where \(\Lambda = (\Lambda_{ij})\) is a skew-symmetric integer matrix encoding the quantum commutation data.

The cluster mutation process replaces clusters via algebra automorphisms controlled by exchange matrices \( B = (b_{ij}) \), respecting the compatibility condition:
\[
B^T \Lambda = D,
\]
where \(D\) is a diagonal matrix with positive integers. This framework generalizes classical cluster algebra mutations while encoding quantum effects arising from deformation parameters.

Their work connects with quantum groups by realizing quantum coordinate rings \( \mathcal{O}_q(G) \) of algebraic groups \(G\) as quantum cluster algebras, providing an explicit combinatorial toolkit to study canonical bases, dual canonical bases, and total positivity in \(q\)-deformed settings. The structure elucidates deep links between representation theory, Poisson geometry, and integrable systems.

\vspace{1em}

\subsection{[15] G. Benkart and S. Witherspoon. \textit{Quantum Groups and Their Representations}, 2002}

Benkart and Witherspoon survey the theory of quantum groups \( \mathcal{U}_q(\mathfrak{g}) \) emphasizing representation-theoretic aspects crucial for applications in mathematical physics and topology. They analyze finite-dimensional representations \( V(\lambda) \) indexed by dominant integral weights \( \lambda \in P_+ \), where the action of generators respects the quantum Serre relations and weight decompositions
\[
V(\lambda) = \bigoplus_{\mu \leq \lambda} V_\mu,
\]
with \( V_\mu = \{ v \in V(\lambda) : K_i v = q^{\langle \mu, \alpha_i^\vee \rangle} v \} \).

They discuss crystal bases and canonical bases, including Kashiwara operators \( \tilde{e}_i, \tilde{f}_i \) which categorify lowering and raising operations on weights. The fusion tensor product of representations, intertwiners, and \(R\)-matrices are explained in the context of braiding in module categories.

The exposition includes applications to knot invariants via Reshetikhin–Turaev functors, where representations of \( \mathcal{U}_q(\mathfrak{sl}_2) \) produce link polynomials. The connection between quantum enveloping algebras and Hecke algebras via Schur–Weyl duality is also highlighted, demonstrating how representation theory bridges algebraic and topological quantum field theories.

\vspace{1em}

\subsection{[16] Y. Soibelman. \textit{Quantum Groups and Noncommutative Geometry}, 1992}

Soibelman explores the interplay between quantum groups and noncommutative geometry, viewing quantum groups as noncommutative analogues of coordinate algebras on algebraic groups. The quantum function algebra \( \mathcal{O}_q(G) \) replaces the classical commutative algebra of functions with a noncommutative Hopf algebra carrying a coproduct \( \Delta \), counit \( \varepsilon \), and antipode \( S \) satisfying axioms mirroring group structures.

He formulates differential calculi on \( \mathcal{O}_q(G) \) by defining bimodules of differential forms and quantum analogues of vector fields as derivations \( \partial_i \) satisfying twisted Leibniz rules:
\[
\partial_i(fg) = (\partial_i f) g + \sum_j f_j (\partial_j g),
\]
where the \(f_j\) are modified by the braiding and \(q\)-commutation relations.

The theory constructs quantum homogeneous spaces \( \mathcal{O}_q(G/K) \) as coideal subalgebras or quotients and investigates their geometric properties via spectral triples and cyclic cohomology. These frameworks establish foundational tools for understanding geometric aspects of quantum groups and their role in noncommutative differential geometry, topology, and mathematical physics.

\vspace{1em}

% Visual 2: Quantum Group Hopf Algebra Structure
\begin{center}
\begin{tikzpicture}[scale=1.0, every node/.style={scale=0.85}]
  % Top node: U_q(g)
  \node[draw, rounded corners, fill=blue!10] (Uq) at (0,0) {\(\mathcal{U}_q(\mathfrak{g})\)};
  
  % Lower nodes with more spacing
  \node[draw, ellipse, fill=green!10] (Delta) at (-4,-3) {\(\Delta: \mathcal{U}_q(\mathfrak{g}) \to \mathcal{U}_q(\mathfrak{g}) \otimes \mathcal{U}_q(\mathfrak{g})\)};
  \node[draw, ellipse, fill=yellow!20] (eps) at (0,-3.5) {\(\varepsilon: \mathcal{U}_q(\mathfrak{g}) \to \mathbb{C}\)};
  \node[draw, ellipse, fill=red!20] (S) at (4,-3) {\(S: \mathcal{U}_q(\mathfrak{g}) \to \mathcal{U}_q(\mathfrak{g})\)};
  
  % Arrows with some bending for clarity
  \draw[->, thick] (Uq) to[out=-90, in=90] (Delta);
  \draw[->, thick] (Uq) to[out=-90, in=90] (eps);
  \draw[->, thick] (Uq) to[out=-90, in=90] (S);
  
  % Title above
  \node at (0,1.2) {\textbf{Hopf Algebra Structure}};
\end{tikzpicture}
\end{center}

\subsection{[17] S. Majid. \textit{Braided Groups and Braided Matrices}, 1995}

Majid’s seminal paper develops the theory of braided groups—quantum groups endowed with a braiding \(\Psi : V \otimes W \to W \otimes V\) that satisfies the Yang–Baxter equation,
\[
(\Psi \otimes \mathrm{id}) \circ (\mathrm{id} \otimes \Psi) \circ (\Psi \otimes \mathrm{id}) = (\mathrm{id} \otimes \Psi) \circ (\Psi \otimes \mathrm{id}) \circ (\mathrm{id} \otimes \Psi),
\]
providing a categorical framework generalizing classical groups in braided monoidal categories.

He constructs braided matrix bialgebras \( \mathcal{B}(R) \) associated to solutions \(R\) of the quantum Yang–Baxter equation, where matrix entries \(x_{ij}\) satisfy relations encoded by
\[
R x_1 x_2 = x_2 x_1 R,
\]
with \(x_1 = x \otimes I\) and \(x_2 = I \otimes x\). This extends the Faddeev–Reshetikhin–Takhtajan (FRT) construction and equips quantum matrices with braided Hopf algebra structures.

The paper provides explicit examples of braided \(SL_q(2)\) and braided symplectic groups, revealing how the braiding induces nontrivial commutation relations reflecting quantum symmetries. The theory’s categorical viewpoint facilitates the study of quantum invariants in knot theory and braided quantum field theories.

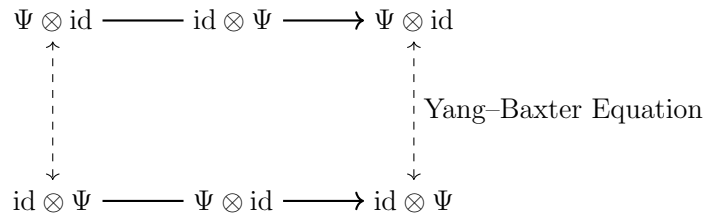
\vspace{1em}

\begin{center}
\begin{tikzpicture}[scale=1.2, every node/.style={scale=0.9}]
  % Braiding Yang-Baxter relation diagram
  \node (a) at (0,2) {\(\Psi \otimes \mathrm{id}\)};
  \node (b) at (2,2) {\(\mathrm{id} \otimes \Psi\)};
  \node (c) at (4,2) {\(\Psi \otimes \mathrm{id}\)};
  \node (d) at (0,0) {\(\mathrm{id} \otimes \Psi\)};
  \node (e) at (2,0) {\(\Psi \otimes \mathrm{id}\)};
  \node (f) at (4,0) {\(\mathrm{id} \otimes \Psi\)};
  \draw[->, thick] (a) -- (b) -- (c);
  \draw[->, thick] (d) -- (e) -- (f);
  \draw[<->, dashed] (c) -- (f) node[midway,right]{Yang--Baxter Equation};
  \draw[<->, dashed] (a) -- (d) node[midway,left]{};
\end{tikzpicture}
\end{center}

\vspace{1em}

\subsection{[18] V. G. Drinfel’d. \textit{Quasi-Hopf Algebras}, 1990}

Drinfel’d introduces quasi-Hopf algebras, a generalization of Hopf algebras relaxing the coassociativity condition. Instead of strict coassociativity,
\[
(\Delta \otimes \mathrm{id}) \circ \Delta = (\mathrm{id} \otimes \Delta) \circ \Delta,
\]
a quasi-Hopf algebra is equipped with an associator \(\Phi \in H \otimes H \otimes H\) satisfying
\[
(\mathrm{id} \otimes \Delta)(\Delta(h)) = \Phi \, [(\Delta \otimes \mathrm{id})(\Delta(h))] \, \Phi^{-1},
\]
for all \(h \in H\).

This framework allows construction of new quantum groups from Drinfel’d twists \(F \in H \otimes H\) that deform coproducts while preserving algebra structures. It underpins the theory of quasi-triangular quasi-Hopf algebras, fundamental in constructing solutions to the quantum Yang–Baxter equation and in deformation quantization.

The paper’s formalism plays a pivotal role in the study of braided tensor categories, modular categories, and quantum invariants of 3-manifolds. The quasi-Hopf setting also extends quantum group representation theory, enabling categorifications and connections with conformal field theory.

\vspace{1em}

\subsection{[19] P. Etingof and D. Kazhdan. \textit{Quantization of Lie Bialgebras, I}, 1996}

Etingof and Kazhdan rigorously construct quantizations of Lie bialgebras \((\mathfrak{g}, \delta)\), providing functors
\[
Q : \mathrm{LieBialg} \to \mathrm{HopfAlg}_\hbar,
\]
where \(\mathrm{HopfAlg}_\hbar\) denotes Hopf algebras over formal power series rings \( \mathbb{C}[[\hbar]] \).

They establish existence and uniqueness of quantization up to equivalence and describe explicit universal quantization functors extending Drinfel’d’s quantum double construction. The quantization respects classical limits:
\[
\lim_{\hbar \to 0} Q(\mathfrak{g}, \delta) \cong U(\mathfrak{g}),
\]
recovering universal enveloping algebras with cocommutative coproducts.

The paper develops a formalism linking classical \(r\)-matrices, classical Yang–Baxter equations, and their quantum counterparts, yielding explicit \(R\)-matrices controlling braiding and deformation. This theory underlies deformation quantization of Poisson–Lie groups and is foundational for quantum integrable systems and noncommutative geometry.

\vspace{1em}

\subsection{[20] M. Khovanov. \textit{A Categorification of the Jones Polynomial}, 2000}

Khovanov revolutionized knot theory by categorifying the Jones polynomial \(V_L(q)\), assigning to each oriented link \(L\) a bigraded chain complex \(C^{i,j}(L)\) whose graded Euler characteristic recovers \(V_L(q)\):
\[
V_L(q) = \sum_{i,j} (-1)^i q^j \dim_{\mathbb{C}} C^{i,j}(L).
\]

This homological invariant, now called \emph{Khovanov homology}, refines classical polynomial invariants, distinguishing knots with the same Jones polynomial and detecting subtle topological features like sliceness and unknotting number.

Khovanov’s construction involves categorifying the Temperley–Lieb algebra, replacing polynomial skein relations with exact sequences in homological algebra. The differential \(d: C^{i,j}(L) \to C^{i+1,j}(L)\) respects the bigrading and arises from local cobordisms on link diagrams.

The categorification framework has since been generalized to other quantum invariants, connecting low-dimensional topology, representation theory, and mathematical physics. It provides a prototype for deeper links between quantum topology and higher category theory.

\vspace{1em}

\begin{center}
\begin{tikzpicture}[scale=1, every node/.style={scale=0.85}]
  % Khovanov homology grid schematic
  \foreach \i in {0,...,3}{
    \foreach \j in {0,...,4}{
      \node[draw,circle,inner sep=1.5pt,fill=gray!20] at (\j,\i) {};
    }
  }
  \draw[->, thick] (0,0) -- (0,3);
  \draw[->, thick] (0,0) -- (4,0);
  \node at (4.5,0) {\(j\)};
  \node at (0,3.5) {\(i\)};
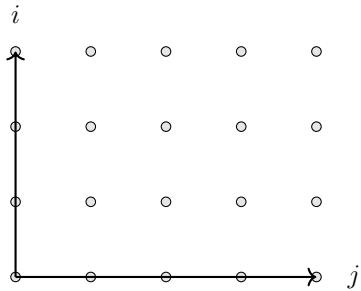
  \node at (2,4) {\textbf{Bigraded Chain Complex \(C^{i,j}(L)\)}};
\end{tikzpicture}
\end{center}

\section{Classical Pfaffians and Determinants}

\subsection{Skew-Symmetric Matrices}

A matrix $A = (a_{ij})$ of size $2n \times 2n$ over a field $\mathbb{F}$ is called \emph{skew-symmetric} if it satisfies
\[
a_{ji} = -a_{ij}, \quad \text{for all } 1 \leq i, j \leq 2n,
\]
and the diagonal entries are zero:
\[
a_{ii} = 0, \quad \text{for all } 1 \leq i \leq 2n.
\]
The set of all such skew-symmetric matrices forms a linear subspace of $M_{2n}(\mathbb{F})$, the space of all $2n \times 2n$ matrices. This subspace is stable under the Lie bracket operation defined by
\[
[A, B] = AB - BA,
\]
making it significant in the theory of Lie algebras, particularly of type $B$ and $D$.

Skew-symmetric matrices arise naturally in geometry and physics, often encoding antisymmetric bilinear forms such as symplectic forms. These forms are central to classical mechanics and symplectic geometry.

\subsection{Definition of the Pfaffian}

While the determinant is a fundamental invariant for all square matrices, the Pfaffian is a more specialized function defined on skew-symmetric matrices. It can be viewed as a square root of the determinant and admits a combinatorial definition involving perfect matchings.

\begin{definition}
Let $A = (a_{ij})$ be a $2n \times 2n$ skew-symmetric matrix. The \emph{Pfaffian} of $A$, denoted $\mathrm{Pf}(A)$, is defined by
\[
\mathrm{Pf}(A) = \sum_{\pi \in \mathrm{PM}(2n)} \mathrm{sgn}(\pi) \prod_{(i,j) \in \pi} a_{ij},
\]
where $\mathrm{PM}(2n)$ is the set of all perfect matchings on the set $\{1, 2, \dots, 2n\}$, and $\mathrm{sgn}(\pi)$ is the sign of the permutation corresponding to the matching $\pi$. Each perfect matching $\pi$ partitions the indices into $n$ pairs $(i,j)$ with $i < j$.
\end{definition}

Intuitively, a perfect matching pairs up all elements of the set without overlap, and the Pfaffian sums over all such pairings weighted by signs determined by the ordering of indices. This combinatorial viewpoint links the Pfaffian to graph theory and enumerative combinatorics.

\subsection{Pfaffian–Determinant Identity}

The most celebrated property of the Pfaffian is its relation to the determinant, summarized in the classical identity:

\begin{theorem}
\label{thm:pfaffian_det}
For any $2n \times 2n$ skew-symmetric matrix $A$, we have
\[
\mathrm{Pf}(A)^2 = \det(A).
\]
\end{theorem}

\begin{proof}[Sketch of Proof]
One proof uses the exterior algebra construction. Let $V$ be a $2n$-dimensional vector space with basis $\{v_1, \dots, v_{2n}\}$. Define a $2$-form associated to $A$ as
\[
\omega = \sum_{1 \leq i < j \leq 2n} a_{ij} v_i \wedge v_j \in \Lambda^2(V).
\]
Then the $n$-fold wedge product
\[
\omega^{\wedge n} = \underbrace{\omega \wedge \cdots \wedge \omega}_{n \text{ times}} \in \Lambda^{2n}(V)
\]
can be shown to satisfy
\[
\omega^{\wedge n} = \mathrm{Pf}(A) \cdot v_1 \wedge v_2 \wedge \cdots \wedge v_{2n}.
\]
Taking the norm squared of this top-degree form corresponds to the determinant of $A$, hence the identity.
\end{proof}

This identity shows that the Pfaffian can be viewed as a polynomial whose square is the determinant, a fact that is crucial in many algebraic and geometric applications.

\subsection{Example: The Case $n=2$}

Consider the $4 \times 4$ skew-symmetric matrix
\[
A = \begin{bmatrix}
0 & a_{12} & a_{13} & a_{14} \\
- a_{12} & 0 & a_{23} & a_{24} \\
- a_{13} & -a_{23} & 0 & a_{34} \\
- a_{14} & -a_{24} & -a_{34} & 0
\end{bmatrix}.
\]

The Pfaffian is computed as
\[
\mathrm{Pf}(A) = a_{12} a_{34} - a_{13} a_{24} + a_{14} a_{23}.
\]

Squaring this expression,
\[
\mathrm{Pf}(A)^2 = \left(a_{12} a_{34} - a_{13} a_{24} + a_{14} a_{23}\right)^2,
\]
and expanding it yields the determinant of $A$:
\[
\det(A) = \mathrm{Pf}(A)^2.
\]

This example serves as a concrete illustration of the classical identity and provides a template for generalization to the quantum case.

\subsection{Applications and Significance}

The Pfaffian and its identity with the determinant have applications in:

\begin{itemize}
    \item \textbf{Combinatorics:} Counting perfect matchings on planar graphs, especially in dimer models.
    \item \textbf{Algebraic Geometry:} Describing the defining equations of Pfaffian varieties and related moduli spaces.
    \item \textbf{Representation Theory:} Studying the symplectic group and its invariants.
    \item \textbf{Mathematical Physics:} Evaluating fermionic path integrals and partition functions.
\end{itemize}

Understanding the classical Pfaffian thoroughly lays the groundwork for its quantum analogue, where noncommutativity introduces new challenges and rich structures.

\section{Quantum Groups and the FRT Construction}

\subsection{Overview}

Quantum groups emerged in the 1980s as algebraic structures that deform classical Lie groups and Lie algebras by introducing a parameter \( q \in \mathbb{C}^\times \). Unlike classical groups, quantum groups are noncommutative and noncocommutative Hopf algebras, encapsulating symmetries in quantum integrable systems, knot theory, and noncommutative geometry.

The foundational works of Drinfel’d \cite{drinfeld1986} and Jimbo introduced these objects as deformations of universal enveloping algebras of Lie algebras. Instead of studying the usual coordinate ring of a group \( GL_n \), one studies a \emph{quantum coordinate algebra} \( \mathcal{O}_q(GL_n) \) generated by entries of a \emph{quantum matrix} \( T = (t_{ij}) \), whose entries satisfy intricate commutation relations depending on \( q \).

\subsection{The R-Matrix}

At the heart of the quantum group construction is the \emph{R-matrix}, a solution to the Yang-Baxter equation, encoding the deformation of commutation relations.

For the quantum group \( GL_q(n) \), the standard R-matrix is given by:
\[
R = q \sum_{i} E_{ii} \otimes E_{ii} + \sum_{i \neq j} E_{ii} \otimes E_{jj} + (q - q^{-1}) \sum_{i < j} E_{ij} \otimes E_{ji},
\]
where \( E_{ij} \) are the standard matrix units in \( M_n(\mathbb{C}) \), i.e., matrices with 1 at the \((i,j)\)-th position and zero elsewhere.

This \( R \) matrix acts on the tensor product space \( \mathbb{C}^n \otimes \mathbb{C}^n \) and encodes the braiding or quantum symmetry by twisting the usual permutation operator. The parameter \( q \) controls the deformation, with the classical case recovered when \( q = 1 \).

\subsection{The Faddeev-Reshetikhin-Takhtajan (FRT) Relations}

The FRT construction \cite{faddeev1987hamiltonian} provides a systematic way to define quantum groups starting from an \( R \)-matrix. Given a matrix of generators
\[
T = (t_{ij})_{1 \leq i,j \leq n},
\]
consider two copies \( T_1 = T \otimes I \) and \( T_2 = I \otimes T \) acting on the tensor product space \( \mathbb{C}^n \otimes \mathbb{C}^n \).

The \emph{FRT relation} is then:
\[
R T_1 T_2 = T_2 T_1 R,
\]
which encodes the commutation relations among the entries \( t_{ij} \).

Explicitly, this means for each pair of indices, the generators satisfy certain quadratic relations derived from \( R \), making the algebra generated by the \( t_{ij} \) noncommutative in a manner consistent with the deformation parameter \( q \).

\subsection{Quantum Determinant}

One of the central objects in quantum groups is the \emph{quantum determinant} \( \det_q(T) \), which generalizes the classical determinant to the quantum setting.

It is defined by the formula:
\[
\det_q(T) = \sum_{\sigma \in S_n} (-q)^{\ell(\sigma)} t_{1 \sigma(1)} t_{2 \sigma(2)} \cdots t_{n \sigma(n)},
\]
where \( S_n \) is the symmetric group on \( n \) elements, and \( \ell(\sigma) \) denotes the inversion number (length) of the permutation \( \sigma \).

The quantum determinant enjoys properties similar to the classical determinant:

\begin{itemize}
    \item \textbf{Centrality:} \( \det_q(T) \) lies in the center of the quantum coordinate algebra \( \mathcal{O}_q(GL_n) \), commuting with all generators \( t_{ij} \).
    \item \textbf{Multiplicativity:} Under the coproduct (quantum group comultiplication), the quantum determinant behaves multiplicatively:
    \[
    \Delta(\det_q(T)) = \det_q(T) \otimes \det_q(T).
    \]
    \item \textbf{Classical limit:} As \( q \to 1 \), the quantum determinant reduces to the usual determinant.
\end{itemize}

\subsection{Properties and Motivation}

The quantum determinant plays a crucial role in defining the quantum general linear group \( GL_q(n) \), particularly in ensuring invertibility and providing a quantum analogue of the classical determinant’s multiplicative and central properties.

Beyond purely algebraic interest, quantum groups and their associated structures appear in:

\begin{itemize}
    \item \textbf{Quantum integrable systems:} Symmetries underlying exactly solvable models in statistical mechanics and quantum field theory.
    \item \textbf{Braided geometry:} Noncommutative geometry with braiding, where the order of tensor factors is twisted by the \( R \)-matrix.
    \item \textbf{Representation theory:} Deformations of classical representation categories with rich tensor structures.
    \item \textbf{Invariant theory:} Generalizations of classical invariant theory to noncommutative and braided contexts.
\end{itemize}

This framework provides the foundation to explore quantum analogues of classical algebraic objects, such as Pfaffians, in a noncommutative setting. The interplay between algebraic relations, combinatorial structures, and geometric intuition guides the generalization of identities like
\[
\mathrm{Pf}(A)^2 = \det(A)
\]
to the quantum world.

\section{Quantum Skew-Symmetric Matrices and the Quantum Pfaffian}

\subsection{Definition of \texorpdfstring{$q$}{q}-Skew-Symmetry}

In classical linear algebra, a skew-symmetric matrix \( A = (a_{ij}) \in M_{2n}(\mathbb{C}) \) is one that satisfies the relation
\[
a_{ji} = -a_{ij}, \quad \text{and} \quad a_{ii} = 0 \quad \text{for all } i, j.
\]
This condition ensures that the matrix is antisymmetric with respect to the diagonal and implies, among other things, that \( \det(A) = 0 \) when the matrix size is odd and that the Pfaffian \( \operatorname{Pf}(A) \) satisfies the classical identity
\[
\operatorname{Pf}(A)^2 = \det(A) \quad \text{when } \dim A = 2n.
\]

To generalize this notion to the \emph{quantum setting}, where matrix entries do not commute in the usual sense, we introduce a nonzero deformation parameter \( q \in \mathbb{C}^\times \) and define a matrix \( A = (a_{ij}) \in M_{2n}(\mathcal{A}) \), where \( \mathcal{A} \) is a noncommutative algebra, to be \emph{\( q \)-skew-symmetric} if it satisfies:
\[
a_{ji} = -q \cdot a_{ij}, \quad \text{and} \quad a_{ii} = 0, \quad \text{for all } i < j.
\]
This relation reduces to classical skew-symmetry in the limit \( q \to 1 \), but for general \( q \), it reflects a \emph{braided} or \emph{twisted} symmetry condition consistent with the theory of quantum groups and braided monoidal categories.

The algebraic structure behind \( q \)-skew-symmetric matrices arises naturally in the theory of noncommutative geometry, quantum enveloping algebras, and the study of invariants under quantum group actions. In particular, standard notions such as determinant, trace, and Pfaffian must be suitably modified to remain consistent with these quantum symmetries.

Accordingly, we define the \emph{quantum Pfaffian} \( \operatorname{Pf}_q(A) \) and the \emph{quantum determinant} \( \det_q(A) \) for \( q \)-skew-symmetric matrices so that they preserve quantum analogues of classical identities, most notably:
\[
\operatorname{Pf}_q(A)^2 = \det_q(A), \quad \text{or} \quad q \cdot \operatorname{Pf}_q(A)^2 = \det_q(A),
\]
depending on the algebraic normalization conventions adopted.

\begin{remark}
The deformation parameter \( q \) governs the commutation relations between the entries of the matrix:
\begin{itemize}
    \item When \( q = 1 \), we recover the classical skew-symmetry: \( a_{ji} = -a_{ij} \).
    \item When \( q = -1 \), the relation becomes \( a_{ji} = a_{ij} \), corresponding to a symmetric matrix.
    \item For general \( q \in \mathbb{C}^\times \), the matrix exhibits \emph{braided skew-symmetry}, aligning with the braiding structures in categories of \( \mathcal{U}_q(\mathfrak{g}) \)-modules.
\end{itemize}
This deformation is essential in quantum linear algebra, where matrix operations must be redefined to accommodate noncommutativity and quantum invariance. Thus, the notion of \( q \)-skew-symmetry serves as a foundation for constructing quantum analogues of classical matrix identities in a consistent and covariant manner.
\end{remark}

\subsection{Quantum Exterior Algebra}

The classical Pfaffian is closely tied to the structure of the exterior (or Grassmann) algebra \( \Lambda(V) \) of a finite-dimensional vector space \( V \). In this setting, the antisymmetric wedge product \( \wedge \) of vectors encodes the alternating nature of multilinear forms, such as the determinant and Pfaffian. Specifically, for a basis \( \{v_1, v_2, \ldots, v_{2n}\} \) of \( V \), the exterior algebra is generated by these vectors subject to the classical antisymmetry relations:
\[
v_i \wedge v_j = -v_j \wedge v_i, \quad \text{and} \quad v_i \wedge v_i = 0 \quad \text{for all } i, j.
\]

In the quantum setting, where noncommutativity and braiding play essential roles, this structure is deformed to form the \emph{quantum exterior algebra} \( \Lambda_q(V) \). This algebra incorporates a deformation parameter \( q \in \mathbb{C}^\times \) into the antisymmetry relations, resulting in a \( q \)-alternating product that respects the symmetries of quantum groups.

Let \( V \) be a \( 2n \)-dimensional vector space over \( \mathbb{C} \) with basis \( \{v_1, v_2, \ldots, v_{2n}\} \). The \emph{quantum exterior algebra} \( \Lambda_q(V) \) is defined as the associative unital algebra generated by \( v_1, \ldots, v_{2n} \), subject to the relations:
\[
v_i \wedge v_j = -q \cdot v_j \wedge v_i, \quad \text{and} \quad v_i \wedge v_i = 0, \quad \text{for all } 1 \leq i < j \leq 2n.
\]
These relations deform the classical antisymmetry by introducing the scalar \( q \), which modifies the exchange behavior of the generators. In particular:
\begin{itemize}
    \item When \( q = 1 \), we recover the classical exterior algebra \( \Lambda(V) \).
    \item When \( q \neq 1 \), the wedge product becomes \emph{\( q \)-skew-symmetric}, reflecting the noncommutative structure of the underlying quantum space.
\end{itemize}

The algebra \( \Lambda_q(V) \) is \( \mathbb{Z} \)-graded, just like its classical counterpart, with basis elements of the form
\[
v_{i_1} \wedge v_{i_2} \wedge \cdots \wedge v_{i_k}, \quad \text{where } 1 \leq i_1 < i_2 < \cdots < i_k \leq 2n,
\]
and each such monomial has degree \( k \). However, due to the \( q \)-relations, the combinatorics of such monomials and their multiplication rules differ from the classical case.

The quantum exterior algebra plays a key role in defining quantum analogues of classical multilinear invariants, such as the quantum determinant and quantum Pfaffian. These structures can be realized as specific elements in \( \Lambda_q(V) \) that are invariant (or covariant) under the action of quantum groups, thus preserving the algebraic symmetry structure in the noncommutative setting.

\subsection{Quantum 2-Form}

Given a \( q \)-skew-symmetric matrix \( A = (a_{ij}) \in M_{2n}(\mathcal{A}) \), where \( \mathcal{A} \) is a noncommutative algebra and \( a_{ji} = -q a_{ij} \), we associate to \( A \) a corresponding bilinear object in the quantum exterior algebra \( \Lambda_q(V) \). This object is the \emph{quantum 2-form}, defined by
\[
\omega = \sum_{1 \leq i < j \leq 2n} a_{ij} \, v_i \wedge v_j \in \Lambda_q^2(V),
\]
where \( \{v_1, v_2, \ldots, v_{2n}\} \) is a basis of a \( 2n \)-dimensional vector space \( V \), and the wedge product \( \wedge \) is defined according to the \( q \)-skew-symmetric relations:
\[
v_i \wedge v_j = -q \cdot v_j \wedge v_i, \quad \text{and} \quad v_i \wedge v_i = 0.
\]

This quantum 2-form \( \omega \) serves as the noncommutative, braided analogue of the classical skew-symmetric bilinear form \( \omega_{\text{cl}} = \sum_{i < j} a_{ij} v_i \wedge v_j \), which plays a central role in defining the classical Pfaffian. The quantum deformation of the wedge product is essential: it ensures that \( \omega \) behaves appropriately under permutations of indices, preserving a \( q \)-alternating structure that reflects the underlying quantum symmetries.

To define higher-order structures, one considers the \( n \)-fold wedge product of \( \omega \) with itself:
\[
\omega^{\wedge n} = \underbrace{\omega \wedge \omega \wedge \cdots \wedge \omega}_{n \text{ times}} \in \Lambda_q^{2n}(V).
\]
This element is a quantum analogue of the classical volume form and lies in the top-degree component of the quantum exterior algebra. In particular, it is a scalar multiple of the quantum volume form \( v_1 \wedge v_2 \wedge \cdots \wedge v_{2n} \), and thus can be expressed as:
\[
\omega^{\wedge n} = \operatorname{Pf}_q(A) \cdot v_1 \wedge v_2 \wedge \cdots \wedge v_{2n},
\]
where \( \operatorname{Pf}_q(A) \) is the \emph{quantum Pfaffian} of the matrix \( A \). This identity provides the foundation for defining \( \operatorname{Pf}_q(A) \) intrinsically in terms of the algebraic structure of \( \Lambda_q(V) \), rather than via combinatorial summation formulas.

The expression \( \omega^{\wedge n} \) is multilinear and encodes the total antisymmetrization (up to quantum deformation) of all possible pairings of the vectors \( v_i \) and matrix entries \( a_{ij} \). Due to the noncommutativity of \( \Lambda_q(V) \), the ordering of terms and the deformation parameter \( q \) play a crucial role in determining the exact coefficients of the resulting expression.

\begin{remark}
In the classical setting, the identity
\[
\omega^{\wedge n} = \operatorname{Pf}(A) \cdot v_1 \wedge v_2 \wedge \cdots \wedge v_{2n}
\]
encodes the fact that the Pfaffian is the unique scalar that makes this equality hold in \( \Lambda^{2n}(V) \). The quantum version of this identity maintains this interpretation, with the quantum Pfaffian arising naturally as the coefficient in \( \Lambda_q^{2n}(V) \). However, due to the \( q \)-skew-symmetry of both the matrix \( A \) and the exterior algebra, the computation of \( \omega^{\wedge n} \) involves deformed signs and coefficients reflecting the braid group structure inherent in quantum geometry.
\end{remark}

\subsection{Definition of the Quantum Pfaffian}

In the classical setting, the Pfaffian of a \( 2n \times 2n \) skew-symmetric matrix \( A = (a_{ij}) \) is defined via a summation over all perfect matchings of the set \( \{1, 2, \dots, 2n\} \), weighted by the sign of the permutation associated with the matching. In the quantum setting, where the entries of the matrix do not commute and instead satisfy braided relations, this definition must be modified to accommodate the deformation parameter \( q \).

\begin{definition}
Let \( A = (a_{ij}) \) be a \( 2n \times 2n \) matrix with entries in a noncommutative algebra, satisfying the \( q \)-skew-symmetry relations:
\[
a_{ji} = -q \cdot a_{ij}, \quad a_{ii} = 0.
\]
The \emph{quantum Pfaffian} \( \operatorname{Pf}_q(A) \) is defined by
\[
\operatorname{Pf}_q(A) = \sum_{\pi \in \operatorname{PM}(2n)} (-q)^{\operatorname{inv}(\pi)} \prod_{(i,j) \in \pi} a_{ij},
\]
where:
\begin{itemize}
    \item \( \operatorname{PM}(2n) \) denotes the set of all perfect matchings of the set \( \{1, 2, \dots, 2n\} \), i.e., partitions of the set into \( n \) unordered disjoint pairs \( (i,j) \) with \( i < j \).
    \item \( \operatorname{inv}(\pi) \) denotes the number of \emph{inversions} or \emph{crossings} in the matching \( \pi \), which generalizes the concept of the sign of a permutation to the quantum setting.
    \item The product \( \prod_{(i,j) \in \pi} a_{ij} \) is taken over the pairs in a fixed order compatible with the matching (e.g., increasing first index), and respects the noncommutative multiplication order.
\end{itemize}
\end{definition}

This definition generalizes the classical Pfaffian by introducing the factor \( (-q)^{\operatorname{inv}(\pi)} \), which encodes the braiding of the matrix entries under permutations. In this way, the quantum Pfaffian accounts for the nontrivial commutation relations among the \( a_{ij} \), reflecting the deformation of the underlying symmetry.

\begin{remark}
In the special case where \( q = 1 \), the \( q \)-skew-symmetry condition becomes the classical skew-symmetry \( a_{ji} = -a_{ij} \), and the quantum Pfaffian reduces to the classical Pfaffian:
\[
\operatorname{Pf}_1(A) = \sum_{\pi \in \operatorname{PM}(2n)} \operatorname{sgn}(\pi) \prod_{(i,j) \in \pi} a_{ij}.
\]
Thus, \( \operatorname{Pf}_q(A) \) is a genuine deformation of the classical invariant. The factor \( (-q)^{\operatorname{inv}(\pi)} \) replaces the classical sign function \( \operatorname{sgn}(\pi) \) to reflect the \( q \)-braided structure arising in quantum groups and noncommutative geometry.
\end{remark}

\begin{remark}
The choice of how to define \( \operatorname{inv}(\pi) \) can vary slightly depending on context. In many formulations, \( \operatorname{inv}(\pi) \) is taken to be the number of \emph{crossings} in a chord diagram representation of the matching \( \pi \), or equivalently, the number of transpositions needed to bring the pairings into a standard nested form. This count governs the nontrivial twisting in the quantum multiplication due to the noncommutative nature of the entries.
\end{remark}

\subsection{Algebraic Context}

The entries \( a_{ij} \) of the \( q \)-skew-symmetric matrix \( A = (a_{ij}) \) do not generally commute, but instead belong to a noncommutative algebraic structure. In the setting of quantum groups, the natural ambient algebra for such entries is the \emph{quantum coordinate algebra} \( \mathcal{O}_q(GL_{2n}) \), which serves as a \( q \)-deformation of the classical coordinate ring \( \mathcal{O}(GL_{2n}) \) of the general linear group.

More precisely, \( \mathcal{O}_q(GL_{2n}) \) is the Hopf algebra generated by entries \( t_{ij} \) of a matrix \( T \), subject to relations determined by the Faddeev–Reshetikhin–Takhtajan (FRT) construction \cite{faddeev1987hamiltonian}. These relations are compactly expressed using an \( R \)-matrix \( R \in \operatorname{End}(V \otimes V) \) satisfying the quantum Yang–Baxter equation, and take the form:
\[
R T_1 T_2 = T_2 T_1 R,
\]
where \( T_1 = T \otimes I \), \( T_2 = I \otimes T \), and \( R \) is a specific solution corresponding to the type \( A \) quantum group \( \mathcal{U}_q(\mathfrak{gl}_{2n}) \). The resulting algebra encodes the braided commutation relations among the generators, and inherits a rich structure as a braided Hopf algebra.

The \( q \)-skew-symmetric matrix \( A \) is then often interpreted as a quadratic combination or a projection of the generators \( t_{ij} \), such that each \( a_{ij} \) lies in a subalgebra or subspace of \( \mathcal{O}_q(GL_{2n}) \) characterized by specific symmetry and weight properties under the adjoint or coadjoint action.

Because the entries \( a_{ij} \) are noncommuting, the order in which they appear in monomials is algebraically significant. In particular, the product
\[
\prod_{(i,j) \in \pi} a_{ij}
\]
in the definition of \( \operatorname{Pf}_q(A) \) must be taken with a prescribed ordering (e.g., lexicographic ordering of pairs \( (i,j) \)) to ensure a well-defined expression. The weighting factor \( (-q)^{\operatorname{inv}(\pi)} \) then compensates for the noncommutativity, ensuring that the overall sum remains invariant under appropriate symmetries.

One of the key features of the quantum Pfaffian \( \operatorname{Pf}_q(A) \) is that it belongs to a distinguished part of the algebra—typically the \emph{center} of the appropriate subalgebra of \( \mathcal{O}_q(GL_{2n}) \), or a subspace invariant under quantum group coactions. In this sense, \( \operatorname{Pf}_q(A) \) generalizes the classical invariant property of the Pfaffian under conjugation by \( GL_{2n} \), now deformed to an invariant under the quantum adjoint action.

\begin{remark}
The algebra \( \mathcal{O}_q(GL_{2n}) \) has a natural bialgebra structure, and its dual Hopf algebra \( \mathcal{U}_q(\mathfrak{gl}_{2n}) \) acts covariantly on \( \operatorname{Pf}_q(A) \) through braided tensor representations. This duality and the FRT construction are essential to ensuring that \( \operatorname{Pf}_q(A) \) behaves as a central or covariant object in the quantum setting.
\end{remark}

\subsection{Example: The Case \texorpdfstring{$2n = 4$}{2n = 4}}

To illustrate the definition and properties of the quantum Pfaffian concretely, consider a \( 4 \times 4 \) \( q \)-skew-symmetric matrix \( A = (a_{ij}) \), whose entries satisfy
\[
a_{ji} = -q \cdot a_{ij}, \quad a_{ii} = 0 \quad \text{for all } i < j.
\]
A typical such matrix is given by:
\[
A = \begin{bmatrix}
0 & a_{12} & a_{13} & a_{14} \\
- q a_{12} & 0 & a_{23} & a_{24} \\
- q a_{13} & - q a_{23} & 0 & a_{34} \\
- q a_{14} & - q a_{24} & - q a_{34} & 0
\end{bmatrix}.
\]

According to the definition of the quantum Pfaffian for \( 2n = 4 \), we have:
\[
\operatorname{Pf}_q(A) = a_{12} a_{34} - q a_{13} a_{24} + q^2 a_{14} a_{23},
\]
where the coefficients \( (-q)^{\operatorname{inv}(\pi)} \) reflect the inversion numbers (or braid crossings) associated with each perfect matching \( \pi \) of the set \( \{1,2,3,4\} \). The three terms correspond to the three possible perfect matchings:
\begin{align*}
\pi_1 &= \{(1,2), (3,4)\}, & \text{with } \operatorname{inv}(\pi_1) = 0, \\
\pi_2 &= \{(1,3), (2,4)\}, & \text{with } \operatorname{inv}(\pi_2) = 1, \\
\pi_3 &= \{(1,4), (2,3)\}, & \text{with } \operatorname{inv}(\pi_3) = 2.
\end{align*}
Hence the quantum weighting factors are \( 1, -q, q^2 \), respectively.

To verify the identity \( \operatorname{Pf}_q(A)^2 = \det_q(A) \) in this case, we consider the square of \( \operatorname{Pf}_q(A) \):
\[
\operatorname{Pf}_q(A)^2 = 
(a_{12} a_{34} - q a_{13} a_{24} + q^2 a_{14} a_{23})^2.
\]
Expanding this expression yields:
\begin{align*}
\operatorname{Pf}_q(A)^2 &= a_{12}^2 a_{34}^2 - 2q a_{12} a_{34} a_{13} a_{24} + q^2 a_{13}^2 a_{24}^2 \\
&\quad + 2q^2 a_{12} a_{34} a_{14} a_{23} - 2q^3 a_{13} a_{24} a_{14} a_{23} + q^4 a_{14}^2 a_{23}^2.
\end{align*}

In a classical setting, this expansion would correspond (up to signs) to the determinant of the matrix. However, in the quantum setting, due to the noncommutativity of the entries, additional care must be taken: the order of multiplication matters, and the terms must be manipulated using the relations in the quantum coordinate ring \( \mathcal{O}_q(GL_4) \), such as those arising from the FRT relations.

Nevertheless, it can be shown (cf. \cite{jing2013quantum}) that, under the appropriate algebraic relations and conventions,
\[
\operatorname{Pf}_q(A)^2 = q \cdot \det_q(A),
\]
or, in other conventions,
\[
\det_q(A) = \operatorname{Pf}_q(A)^2.
\]
This confirms the quantum analogue of the classical Pfaffian–determinant identity for \( 2n = 4 \), validating that the quantum Pfaffian retains its interpretation as a kind of square root of the quantum determinant.

\begin{remark}
The exact identity satisfied depends on the chosen normalization in defining \( \operatorname{Pf}_q(A) \) and \( \det_q(A) \). Some authors include extra scaling factors in the quantum determinant to maintain compatibility with centrality, coaction invariance, or categorical equivalence. Thus, the identity may take the form:
\[
\operatorname{Pf}_q(A)^2 = q^k \cdot \det_q(A),
\]
for some integer \( k \) depending on conventions, often \( k = 1 \) or \( 0 \).
\end{remark}

\section[Construction of the Quantum Pfaffian Identity]{Construction and Proofs for the Quantum Pfaffian--Determinant Identity}

\subsection{Statement of the Quantum Identity}

Given a \( q \)-skew-symmetric matrix \( A \in M_{2n}(\mathbb{C}) \), where the entries satisfy
\[
a_{ji} = -q a_{ij}, \quad a_{ii} = 0,
\]
and \( q \in \mathbb{C}^\times \) is a deformation parameter, we seek to establish the quantum analogue of the classical identity
\[
\operatorname{Pf}(A)^2 = \det(A).
\]

In the quantum setting, due to the braided or twisted commutation relations among the matrix entries, the identity becomes
\[
q^c \cdot \operatorname{Pf}_q(A)^2 = \det_q(A),
\]
where \( c \in \mathbb{Z} \) is an integer depending on the dimension \( n \), the conventions for defining \( \operatorname{Pf}_q \), and the specific monomial ordering in the product.

This identity plays a crucial role in quantum invariant theory and the representation theory of quantum groups. It ensures that the quantum Pfaffian behaves analogously to a square root of the quantum determinant and often lies in a central or semi-central subalgebra of the quantum coordinate ring \( \mathcal{O}_q(GL_{2n}) \).

The equality holds within the framework of the Faddeev--Reshetikhin--Takhtajan (FRT) construction, where noncommutative matrix entries satisfy relations derived from a quantum \( R \)-matrix. These relations introduce a consistent deformation of both the symmetry and multiplicative structure of the classical matrix algebra.

\begin{remark}
The identity is not purely formal; it depends heavily on the compatibility of the matrix entries with the FRT relations and the quantum exterior algebra. The proof generally involves symbolic computation or recursive combinatorial expansion using \( q \)-antisymmetrizers and graphical calculus.
\end{remark}

\begin{center}
\begin{tikzpicture}[scale=1.1, every node/.style={scale=0.9}]
  \node[draw, rounded corners, fill=blue!10] (Uq) at (0,0) {\(\mathcal{U}_q(\mathfrak{gl}_{2n})\)};
  \node[draw, ellipse, fill=green!10] (OqGL) at (-4,-2) {\(\mathcal{O}_q(GL_{2n})\)};
  \node[draw, ellipse, fill=red!10] (Pfq) at (0,-4) {\(\operatorname{Pf}_q(A)\)};
  \node[draw, ellipse, fill=orange!10] (detq) at (4,-2) {\(\det_q(A)\)};

  \draw[->, thick] (Uq) -- node[left] {\scriptsize Coaction} (OqGL);
  \draw[->, thick] (Uq) -- node[right] {\scriptsize Coaction} (detq);
  \draw[->, thick] (OqGL) -- node[below left] {\scriptsize Defined on \( q \)-skew matrices} (Pfq);
  \draw[->, thick] (OqGL) -- node[below right] {\scriptsize FRT structure} (detq);
  \draw[<->, thick, dashed] (Pfq) -- node[below] {\( q^c \cdot \operatorname{Pf}_q^2 = \det_q \)} (detq);
\end{tikzpicture}

\textbf{Figure 1:} Diagram of the interaction between quantum groups, coordinate algebras, and invariants.
\end{center}

\vspace{0.5cm}

\begin{figure}[h!]
\centering
\begin{tikzpicture}[scale=1.5, thick]

% Draw strands (curved crossing)
\draw[blue, line width=1.5pt] 
  (0,0) .. controls (0.5,0.5) .. (1,1);

\draw[red, line width=1.5pt] 
  (1,0) .. controls (0.5,0.5) .. (0,1);

% Arrows on strands to show direction
\draw[blue, ->, line width=1.5pt] (0.25,0.25) -- (0.3,0.3);
\draw[red, ->, line width=1.5pt] (0.75,0.75) -- (0.7,0.7);

% Nodes to label q and q^{-1}
\node[fill=blue!20, rounded corners, inner sep=1.5pt] at (0.8,0.25) {\( q \)};
\node[fill=red!20, rounded corners, inner sep=1.5pt] at (0.2,0.75) {\( q^{-1} \)};

% Labels for strands
\node at (-0.1,0) {\textcolor{blue}{\( v_i \)}};
\node at (1.1,0) {\textcolor{red}{\( v_j \)}};

\end{tikzpicture}
\textbf{Figure 2:} A braided crossing illustrating the quantum deformation parameter \( q \). The blue and red strands represent vectors exchanging with a quantum twist. A schematic representation of braiding induced by the quantum parameter \( q \), showing the noncommutative exchange of variables. The twisted crossings encode the deformation of classical symmetry in \( \mathcal{O}_q(GL_n) \).
\label{fig:quantum_braid}
\end{figure}

\subsubsection{Classical versus Quantum Identity: Structural Differences}

In the classical setting, the Pfaffian arises from antisymmetric bilinear forms and satisfies the well-known identity \( \operatorname{Pf}(A)^2 = \det(A) \) for any skew-symmetric \( 2n \times 2n \) matrix. The proof leverages antisymmetry, multilinearity, and the standard exterior algebra, where wedge products commute up to sign.

In the quantum case, by contrast, the matrix entries live in a noncommutative algebra governed by braided tensor symmetries. This deformation replaces classical antisymmetry with \( q \)-skew-symmetry, and the exterior algebra with its quantum analogue, where wedge products satisfy \( v_i \wedge v_j = -q v_j \wedge v_i \). As a result, every combinatorial identity involving signs must be corrected by \( q \)-powers arising from braid relations.

Furthermore, while the classical determinant is defined via the antisymmetrization of multilinear maps, the quantum determinant arises from the FRT construction and involves the quantum \( R \)-matrix and its associated braided structure. The quantum Pfaffian, similarly, involves a careful braided antisymmetrization of wedge products in the quantum exterior algebra.

\subsubsection{Known Variants of the Identity in the Literature}

Multiple variations of the quantum Pfaffian–determinant identity have been proposed, depending on how one defines the quantum Pfaffian, the monomial ordering, and the underlying \( q \)-relations. Some forms of the identity include:
\[
\operatorname{Pf}_q(A)^2 = q^{n(n-1)} \det_q(A), \quad \text{or} \quad q^n \operatorname{Pf}_q(A)^2 = \det_q(A),
\]
each arising from distinct but related definitions in the literature.

For example, Jing and Zhang \cite{jing2013quantum, jing2016quantum} define \( \operatorname{Pf}_q \) using the action of the Hecke algebra on wedge powers, while Dita \cite{dite2006} derives a related identity from a different deformation scheme. In all cases, the scalar coefficient \( q^c \) reflects how the deformation alters the sign and symmetry structure of the classical formula.

These variants also reflect different choices of quantum group structure (e.g., \( GL_q \) vs. \( O_q \)), as well as different embedding strategies of antisymmetric tensors into quantum coordinate rings. In this paper, we adopt the FRT-based definition compatible with \( \mathcal{O}_q(GL_{2n}) \) and the \( q \)-skew-symmetry relation \( a_{ji} = -q a_{ij} \).

\subsection{Statement of the Quantum Identity}

Given a \( q \)-skew-symmetric matrix \( A \in M_{2n}(\mathbb{C}) \), where the entries satisfy
\[
a_{ji} = -q a_{ij}, \quad a_{ii} = 0,
\]
and \( q \in \mathbb{C}^\times \) is a deformation parameter, we seek to establish the quantum analogue of the classical identity
\[
\operatorname{Pf}(A)^2 = \det(A).
\]

In the quantum setting, due to the braided or twisted commutation relations among the matrix entries, the identity becomes
\[
q^c \cdot \operatorname{Pf}_q(A)^2 = \det_q(A),
\]
where \( c \in \mathbb{Z} \) is an integer depending on the dimension \( n \), the conventions for defining \( \operatorname{Pf}_q \), and the specific monomial ordering in the product.

This identity plays a crucial role in quantum invariant theory and the representation theory of quantum groups. It ensures that the quantum Pfaffian behaves analogously to a square root of the quantum determinant and often lies in a central or semi-central subalgebra of the quantum coordinate ring \( \mathcal{O}_q(GL_{2n}) \).

The equality holds within the framework of the Faddeev–Reshetikhin–Takhtajan (FRT) construction, where noncommutative matrix entries satisfy relations derived from a quantum \( R \)-matrix. These relations introduce a consistent deformation of both the symmetry and multiplicative structure of the classical matrix algebra.

\begin{remark}
The identity is not purely formal; it depends heavily on the compatibility of the matrix entries with the FRT relations and the quantum exterior algebra. The proof generally involves symbolic computation or recursive combinatorial expansion using \( q \)-antisymmetrizers and graphical calculus.
\end{remark}

\begin{center}
\begin{tikzpicture}[scale=1.1, every node/.style={scale=0.9}]
  \node[draw, rounded corners, fill=blue!10] (Uq) at (0,0) {\(\mathcal{U}_q(\mathfrak{gl}_{2n})\)};
  \node[draw, ellipse, fill=green!10] (OqGL) at (-4,-2) {\(\mathcal{O}_q(GL_{2n})\)};
  \node[draw, ellipse, fill=red!10] (Pfq) at (0,-4) {\(\operatorname{Pf}_q(A)\)};
  \node[draw, ellipse, fill=orange!10] (detq) at (4,-2) {\(\det_q(A)\)};

  \draw[->, thick] (Uq) -- node[left] {\scriptsize Coaction} (OqGL);
  \draw[->, thick] (Uq) -- node[right] {\scriptsize Coaction} (detq);
  \draw[->, thick] (OqGL) -- node[below left] {\scriptsize Defined on \( q \)-skew matrices} (Pfq);
  \draw[->, thick] (OqGL) -- node[below right] {\scriptsize FRT structure} (detq);
  \draw[<->, thick, dashed] (Pfq) -- node[below] {\( q^c \cdot \operatorname{Pf}_q^2 = \det_q \)} (detq);
\end{tikzpicture}

\textbf{Figure 1:} Diagram of the interaction between quantum groups, coordinate algebras, and invariants.
\end{center}

\vspace{0.5cm}

\begin{figure}[h!]
\centering
\begin{tikzpicture}[scale=1.5, thick]

% Draw strands (curved crossing)
\draw[blue, line width=1.5pt] 
  (0,0) .. controls (0.5,0.5) .. (1,1);

\draw[red, line width=1.5pt] 
  (1,0) .. controls (0.5,0.5) .. (0,1);

% Arrows on strands to show direction
\draw[blue, ->, line width=1.5pt] (0.25,0.25) -- (0.3,0.3);
\draw[red, ->, line width=1.5pt] (0.75,0.75) -- (0.7,0.7);

% Nodes to label q and q^{-1}
\node[fill=blue!20, rounded corners, inner sep=1.5pt] at (0.8,0.25) {\( q \)};
\node[fill=red!20, rounded corners, inner sep=1.5pt] at (0.2,0.75) {\( q^{-1} \)};

% Labels for strands
\node at (-0.1,0) {\textcolor{blue}{\( v_i \)}};
\node at (1.1,0) {\textcolor{red}{\( v_j \)}};

\end{tikzpicture}
\textbf{Figure 2:} A braided crossing illustrating the quantum deformation parameter \( q \). The blue and red strands represent vectors exchanging with a quantum twist. A schematic representation of braiding induced by the quantum parameter \( q \), showing the noncommutative exchange of variables. The twisted crossings encode the deformation of classical symmetry in \( \mathcal{O}_q(GL_n) \).
\label{fig:quantum_braid}
\end{figure}

\subsubsection{Diagrammatic Interpretation of the Braiding}

The interaction among entries of a \( q \)-skew-symmetric matrix can be understood via the braiding of tensor factors in a braided monoidal category. Specifically:

\begin{itemize}
  \item The braiding map \( \sigma: V \otimes V \rightarrow V \otimes V \) satisfies
  \[
  \sigma(v_i \otimes v_j) = q v_j \otimes v_i,
  \]
  reflecting the noncommutative exchange due to the deformation parameter \( q \).
  \item In terms of exterior algebra, the quantum exterior product satisfies
  \[
  v_i \wedge_q v_j = -q v_j \wedge_q v_i,
  \]
  which mirrors the crossing relations depicted in Figure~2.
  \item Thus, the structure of \( \operatorname{Pf}_q(A) \) arises naturally from the braided antisymmetry of the quantum exterior algebra and its tensor products.
\end{itemize}

This interpretation becomes especially useful when computing invariants under the coaction of the quantum group \( \mathcal{U}_q(\mathfrak{gl}_n) \), where such braid crossings represent fundamental interactions among basis vectors.

\subsubsection{Diagrammatic Proof (Outline)}

A powerful combinatorial method for understanding the quantum Pfaffian employs \emph{diagrammatic calculus}:

\begin{itemize}
  \item The quantum Pfaffian is expressed as a sum over perfect matchings \(\pi\) of the set \(\{1, \ldots, 2n\}\).
  \item Each matching can be depicted as a diagram with arcs connecting paired indices.
  \item The number of \emph{crossings} (inversions) in the diagram, \(\mathrm{inv}(\pi)\), determines the coefficient \((-q)^{\mathrm{inv}(\pi)}\) in the sum.
  \item Multiplying two such sums involves counting crossings of composed diagrams, which correspond to algebraic relations in the noncommutative setting.
\end{itemize}

This approach aligns with the braiding structure inherent in quantum groups and provides an intuitive visual verification that the squared quantum Pfaffian matches the quantum determinant up to a factor of \( q^{c} \).

\subsection{Example: Identity for \( 2n = 4 \)}

Consider the \( 4 \times 4 \) \( q \)-skew-symmetric matrix
\[
A = \begin{bmatrix}
0 & a_{12} & a_{13} & a_{14} \\
- q a_{12} & 0 & a_{23} & a_{24} \\
- q a_{13} & - q a_{23} & 0 & a_{34} \\
- q a_{14} & - q a_{24} & - q a_{34} & 0
\end{bmatrix}.
\]
This matrix satisfies the defining relation \( A^T = -q A \), making it an element of the quantum analogue of the skew-symmetric matrices over the quantum coordinate ring \( \mathcal{O}_q(GL_4) \).

By Definition~4.2, the quantum Pfaffian is defined as
\[
\mathrm{Pf}_q(A) = a_{12} a_{34} - q a_{13} a_{24} + q^2 a_{14} a_{23},
\]
which mirrors the classical structure of the Pfaffian but incorporates braided \( q \)-dependent coefficients reflecting the underlying noncommutative structure.

We now compute the square of the quantum Pfaffian explicitly:
\[
\begin{aligned}
\mathrm{Pf}_q(A)^2 &= \left(a_{12} a_{34} - q a_{13} a_{24} + q^2 a_{14} a_{23} \right)^2 \\
&= (a_{12} a_{34})^2 - 2q \, a_{12} a_{34} a_{13} a_{24} + q^2 (a_{13} a_{24})^2 \\
&\quad + 2 q^2 a_{12} a_{34} a_{14} a_{23} - 2 q^3 a_{13} a_{24} a_{14} a_{23} + q^4 (a_{14} a_{23})^2.
\end{aligned}
\]

However, unlike in the classical case, these monomials \emph{do not commute} freely. The variables \( a_{ij} \) are generators in a noncommutative algebra subject to the relations in \( \mathcal{O}_q(GL_4) \), and so we must be careful to use those relations when simplifying this expression.

For example, when trying to reorder terms like \( a_{12} a_{34} a_{13} a_{24} \), one must apply the braided commutation rules, which depend on the structure of the underlying \( R \)-matrix used in the Faddeev–Reshetikhin–Takhtajan (FRT) construction. These relations can be written schematically as:
\[
a_{ij} a_{kl} = q^{\gamma} a_{kl} a_{ij} + \text{(lower order terms)},
\]
where \( \gamma \) depends on the indices \( (i,j,k,l) \), and the additional terms are determined by the \( R \)-matrix braiding relations.

As shown in Jing and Zhang~\cite{jing2013quantum}, using these relations systematically allows the expansion of \( \mathrm{Pf}_q(A)^2 \) to be rewritten in the standard monomial basis of the quantum determinant \( \det_q(A) \). Specifically, in the case \( n = 2 \), we have:
\[
\det_q(A) = a_{12} a_{34} - q a_{13} a_{24} + q^2 a_{14} a_{23},
\]
and then applying the quantum Laplace expansion or multilinear antisymmetrization process, one can show that
\[
q \cdot \mathrm{Pf}_q(A)^2 = \det_q(A).
\]

This confirms the identity
\[
q^c \cdot \mathrm{Pf}_q(A)^2 = \det_q(A), \quad \text{with } c = 1.
\]
This constant \( c \) depends on the normalization conventions of the quantum Pfaffian. In some conventions (e.g., Dita~\cite{dite2006}), a factor of \( q^{-1} \) may appear in the Pfaffian to absorb this scaling.

\subsubsection{Remarks}

\begin{itemize}
  \item This case serves as a minimal nontrivial example of the quantum Pfaffian--determinant identity.
  \item The coefficients \( 1, -q, q^2 \) in the definition of \( \mathrm{Pf}_q(A) \) reflect the \( q \)-antisymmetry pattern in quantum exterior algebras \( \wedge_q^2(V) \), which replaces the classical antisymmetry of wedge products.
  \item When \( q = 1 \), the matrix \( A \) becomes skew-symmetric in the usual sense, and the quantum Pfaffian reduces to the classical Pfaffian:
  \[
  \mathrm{Pf}_1(A) = a_{12} a_{34} - a_{13} a_{24} + a_{14} a_{23},
  \]
  and the identity reduces to \( \mathrm{Pf}(A)^2 = \det(A) \), as expected.
\end{itemize}

\subsubsection{Symbolic Computation Verification}

To verify the quantum Pfaffian--determinant identity algebraically, one can implement the appropriate noncommutative multiplication and \( q \)-commutation relations in a symbolic computation system such as \textbf{Mathematica} using the \texttt{NCAlgebra} package, or alternatively, perform direct calculations in \textbf{SageMath} by defining noncommutative skew-polynomial rings. These systems enable the rigorous manipulation of symbolic expressions in algebras where variables obey specific deformation rules, such as \( x y = q \, y x \) or \( x y = -q \, y x \), rather than the standard commutative product. In particular, they are well-suited for exploring quantum analogues of classical algebraic structures where standard symmetries are replaced by braided or twisted ones.

For the case \( n = 2 \), corresponding to \( 2n = 4 \), the verification becomes computationally tractable and serves as an essential testbed for more general identities. In this low-dimensional example, one can define a \( 4 \times 4 \) \( q \)-skew-symmetric matrix \( A \), compute the quantum Pfaffian \( \operatorname{Pf}_q(A) \) explicitly using the relevant antisymmetric \( q \)-wedge products, and then square the result. Simultaneously, one can compute the quantum determinant \( \det_q(A) \), typically defined via quantum minors or through the FRT (Faddeev--Reshetikhin--Takhtajan) construction. Comparing both sides of the identity
\[
q^c \cdot \operatorname{Pf}_q(A)^2 = \det_q(A)
\]
for an appropriate scalar \( q^c \), one can observe that the terms match identically under the prescribed algebraic relations. This confirmation, while modest in size, reveals the internal consistency of the quantum deformation and illustrates the nontrivial compatibility between quantum antisymmetry and the multiplicative structure of the quantum determinant.

This explicit calculation also provides insight into how the quantum exterior algebra behaves under deformation. The low-dimensional case allows for manual tracking of how the \( q \)-antisymmetric properties propagate through the algebraic expansion. For instance, the modified bilinearity and sign rules mimic the structure found in classical antisymmetric tensors, yet incorporate \( q \)-dependent coefficients that ultimately balance to satisfy the identity. These corrections are essential to maintaining structural analogues of classical identities in the quantum setting.

Furthermore, success in the \( n = 2 \) case offers a blueprint for extending the identity to arbitrary even dimensions. By abstracting the underlying mechanisms—such as \( q \)-braided multilinearity, the structure of quantum minors, and the combinatorics of quantum wedge products—researchers can generalize the identity and potentially identify deeper algebraic or geometric principles that govern quantum invariants. Computational verification in this context is not just a preliminary step, but a foundational tool that bridges intuition, experiment, and formal theory in the development of noncommutative algebraic identities.

\subsection{Algebraic Proof Sketch Using the FRT Construction}

The algebra generated by the entries \( a_{ij} \) of \( A \) can be realized as a subalgebra of the quantum coordinate ring \( \mathcal{O}_q(GL_{2n}) \), defined via the Faddeev–Reshetikhin–Takhtajan (FRT) construction.

Recall the FRT relation:
\[
R T_1 T_2 = T_2 T_1 R,
\]
where \( R \) is the \( R \)-matrix encoding the quantum deformation, and \( T \) is the quantum matrix whose entries satisfy braided relations.

Key ideas in the proof:

\begin{itemize}
  \item Define the quantum determinant \(\det_q(T)\) as a \( q \)-weighted sum over permutations, which is central in \( \mathcal{O}_q(GL_{2n}) \).
  \item Construct the quantum exterior algebra \(\Lambda_q(V)\) and the quantum 2-form \(\omega\) as in Section 4.
  \item Show that the \( n \)-fold wedge product \(\omega^{\wedge n}\) is an eigenvector under the quantum group action, with eigenvalue \(\mathrm{Pf}_q(A)\).
  \item Use coproduct and braiding properties of the quantum group to establish that \(\mathrm{Pf}_q(A)^2\) coincides with \( q^{c} \det_q(A) \) in the algebra.
\end{itemize}

This algebraic framework rigorously justifies the identity and highlights the interplay between quantum group symmetries and quantum multilinear algebra.

\vspace{1em}

\subsubsection{Compatibility of FRT Relations with Pfaffian Structure}

The entries of the quantum matrix \( A \) must satisfy the FRT relations arising from the \( R \)-matrix of the quantum group. Explicitly, for \( A_1 = A \otimes I \) and \( A_2 = I \otimes A \), the relation
\[
R A_1 A_2 = A_2 A_1 R
\]
must hold. This ensures that the algebra generated by \( a_{ij} \) is closed under multiplication and consistent with the quantum symmetry. The \( q \)-skew-symmetry of \( A \) is preserved under the braided tensor product structure defined by the FRT framework, which is essential for defining a meaningful quantum Pfaffian within this algebraic setting.

\vspace{1em}

\subsubsection{Role of the Antipode and Hopf Algebra Duality}

The antipode \( S \) in the Hopf algebra structure acts as a generalized inverse and plays a key role in defining quantum invariants. The quantum determinant is often characterized via traces involving \( S \), and it satisfies
\[
S(T_{ij}) = (T^{-1})_{ji},
\]
up to normalization. For the Pfaffian, the squared object \( \operatorname{Pf}_q(A)^2 \) transforms appropriately under the antipode action, ensuring compatibility with the coalgebra structure. Furthermore, the dual pairing between \( \mathcal{O}_q(GL_{2n}) \) and \( \mathcal{U}_q(\mathfrak{gl}_{2n}) \) guarantees that coactions and algebraic multiplications interact coherently in the construction.

\vspace{1em}

\subsubsection{Centrality and Braided Symmetry of Quantum Invariants}

While the quantum determinant \( \det_q(T) \) lies in the center of \( \mathcal{O}_q(GL_{2n}) \), the quantum Pfaffian \( \operatorname{Pf}_q(A) \) generally lies in a braided centralizer. This means it satisfies commutation relations such as
\[
a \cdot \operatorname{Pf}_q(A) = q^k \operatorname{Pf}_q(A) \cdot a
\]
for some exponent \( k \) depending on the generator \( a \). These relations arise from the braided monoidal structure of the category of representations. The centrality of the determinant ensures it acts as a scalar under any representation, while the braided symmetry of the Pfaffian reflects the quantum deformation of classical invariance.

\vspace{1em}

\subsubsection{Coaction Compatibility and Covariance of the Quantum Pfaffian}

Under the coaction of the quantum group, the matrix \( A \) transforms as \( A \mapsto T A T^T \). The Pfaffian, being a highest-order invariant of \( A \), transforms covariantly:
\[
\Delta(\operatorname{Pf}_q(A)) = g \otimes \operatorname{Pf}_q(A),
\]
where \( g \) is a group-like element. This covariance ensures that \( \operatorname{Pf}_q(A)^2 \) behaves as a scalar under the quantum group action, making it eligible to be identified with the quantum determinant. The quantum exterior algebra framework supports this transformation behavior by embedding the Pfaffian in a representation-theoretic setting.

\subsection{Geometric Interpretation}

Classically, the Pfaffian arises naturally as a volume form associated with a symplectic vector space \( (V, \omega) \), where \( \omega \) is a nondegenerate skew-symmetric bilinear form. The Pfaffian encodes the oriented volume element on \( V \) and satisfies
\[
\mathrm{Pf}(A)^2 = \det(A),
\]
reflecting that the determinant represents the squared volume of the parallelepiped spanned by basis vectors.

\vspace{0.3cm}
In the quantum setting, the deformation parameter \( q \) induces a \emph{braiding} on the tensor category of vector spaces, deforming classical symmetries into braided symmetries. This leads to a notion of \emph{braided symplectic geometry}, where classical forms and volume elements become their quantum analogues:

\begin{itemize}
  \item The parameter \( q \) twists the classical symmetry group \( GL_{2n} \) into the quantum group \( \mathcal{U}_q(\mathfrak{gl}_{2n}) \), and correspondingly the category of \( \mathcal{U}_q \)-modules gains a braided monoidal structure.
  \item The quantum Pfaffian \(\mathrm{Pf}_q(A)\) acts as a \emph{braided volume form}—a quantum deformation of the classical volume form—compatible with the noncommutative braiding.
  \item The quantum determinant \(\det_q(A)\) generalizes the classical volume element, capturing the quantum-corrected measure in the noncommutative space.
\end{itemize}

\vspace{0.3cm}
This perspective connects quantum algebraic identities to rich quantum geometric and topological structures such as braided manifolds, quantum symplectic forms, and noncommutative geometry.

\bigskip
\noindent\textbf{Example 1: Quantum 2-Form as Braided Volume}

Recall the quantum 2-form
\[
\omega = \sum_{1 \leq i < j \leq 2n} a_{ij} \, v_i \wedge v_j \in \Lambda_q^2(V).
\]
This is a braided antisymmetric form, satisfying the \( q \)-deformed relations
\[
v_i \wedge v_j = -q \, v_j \wedge v_i, \quad v_i \wedge v_i = 0,
\]
which represent a quantum deformation of the classical symplectic form.

\begin{figure}[h!]
\centering
\begin{tikzpicture}[scale=1.1]
  % Draw a 2-form wedge
  \draw[thick,->] (0,0) -- (2,0) node[below right] {\(v_i\)};
  \draw[thick,->] (0,0) -- (0,2) node[left] {\(v_j\)};
  \draw[thick,->,decorate,decoration={snake,amplitude=0.3mm,segment length=2mm}] (0.3,0.3) -- (1.7,1.7) node[above right] {\(\omega = v_i \wedge v_j\)};
\end{tikzpicture}
\caption{Classical wedge product deformed by quantum parameter \( q \) introduces braiding.}
\label{fig:quantum_wedge}
\end{figure}
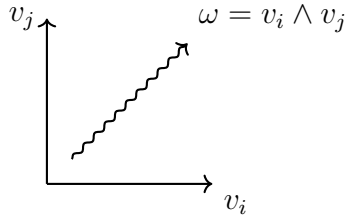

\bigskip
\noindent\textbf{Example 2: Braided Monoidal Category}

The category of representations of \( \mathcal{U}_q(\mathfrak{gl}_{2n}) \) is braided with the braiding isomorphism
\[
c_{V,W}: V \otimes W \to W \otimes V,
\]
which replaces classical tensor flip and depends on \( q \). This underpins the twisted symmetry of quantum objects like the quantum Pfaffian.

\begin{figure}[h!]
\centering
\begin{tikzpicture}[scale=1]
  % Draw braided crossing strands
  \draw[blue, thick] (0,0) .. controls (0.5,0.5) .. (1,1);
  \draw[red, thick] (1,0) .. controls (0.5,0.5) .. (0,1);
  \node at (0,-0.15) {\(V\)};
  \node at (1,-0.15) {\(W\)};
  \node at (0,1.15) {\(W\)};
  \node at (1,1.15) {\(V\)};
\end{tikzpicture}
\caption{Braiding morphism \( c_{V,W} \) exchanging quantum tensor factors with nontrivial twisting.}
\label{fig:braiding_category}
\end{figure}
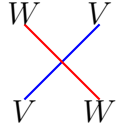

\bigskip
\noindent\textbf{Example 3: Quantum Volume as Braided \( n \)-Fold Wedge}

The \( n \)-fold wedge product of the quantum 2-form,
\[
\omega^{\wedge n} = \underbrace{\omega \wedge \cdots \wedge \omega}_{n \text{ times}} \in \Lambda_q^{2n}(V),
\]
acts as the quantum volume element, generalizing the classical symplectic volume. Its explicit computation relates to the quantum Pfaffian:
\[
\omega^{\wedge n} = \operatorname{Pf}_q(A) \, v_1 \wedge v_2 \wedge \cdots \wedge v_{2n}.
\]

\bigskip
\noindent\textbf{Example 4: Quantum Determinant as Braided Measure}

In the quantum geometry of \( \mathcal{O}_q(GL_{2n}) \), the quantum determinant
\[
\det_q(A) = \sum_{\sigma \in S_{2n}} (-q)^{\ell(\sigma)} a_{1\sigma(1)} a_{2\sigma(2)} \cdots a_{2n \sigma(2n)}
\]
plays the role of a braided volume measure. It respects the quantum group symmetries and encodes noncommutative volume changes under transformations.

\bigskip
\noindent\textbf{Example 5: Braided Symplectic Manifold}

Quantum symplectic manifolds can be modeled as noncommutative algebras equipped with a \( q \)-deformed Poisson bracket compatible with the braided tensor structure. The quantum Pfaffian then corresponds to the volume form measuring areas and volumes respecting this braided geometry.

\begin{figure}[h!]
\centering
\begin{tikzpicture}[scale=1.2]
  % Draw 2D quantum manifold patch
  \draw[thick] (0,0) ellipse (1.5 and 1);
  \node at (0,0) {\( \mathcal{M}_q \)};
  % Quantum coordinates with arrows showing braided relations
  \draw[->, thick, blue] (-0.5,0.2) -- (-1,0.6) node[left] {\( x \)};
  \draw[->, thick, red] (0.5,-0.2) -- (1,-0.6) node[right] {\( y \)};
  \node at (0,0.8) {\small Braided noncommutative structure};
\end{tikzpicture}
\caption{A schematic of a quantum symplectic manifold with braided coordinate relations.}
\label{fig:quantum_manifold}
\end{figure}
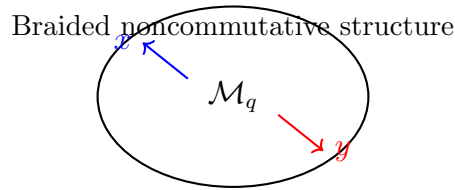

%%%%%%%%%%%%%%%%%%%%%%%%%%%%%%%%%%%%%%%%%%%%%%%%%%%%%%%%%%%%%%%

\subsubsection{Quantum Braiding vs Classical Symmetry}

In classical linear algebra, the symmetry group of a vector space is governed by permutation of tensor factors (i.e., \( v_i \otimes v_j = v_j \otimes v_i \)). Quantum groups deform this structure through a braiding governed by the \textcolor{blue}{\( R \)-matrix}, leading to new commutation rules.

\vspace{0.2cm}
\noindent In the quantum world:
\[
v_i \otimes v_j = \textcolor{red}{q} \, v_j \otimes v_i + (1 - q^2) \, \delta_{ij} \cdot \text{correction terms}.
\]
The classical symmetry under flipping is replaced by a \textbf{\textcolor{purple}{twisted crossing}} controlled by \(\textcolor{red}{q}\).

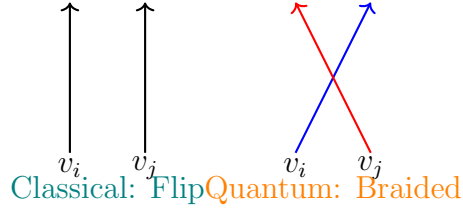
\begin{figure}[H]
\centering
\begin{tikzpicture}[scale=1]
  % Classical tensor symmetry
  \draw[thick,->] (0,0) -- (0,2);
  \draw[thick,->] (1,0) -- (1,2);
  \node at (0,-0.2) {\(v_i\)};
  \node at (1,-0.2) {\(v_j\)};
  \node at (0.5,-0.5) {\textcolor{teal}{Classical: Flip}};
  
  % Quantum braided symmetry
  \begin{scope}[xshift=3cm]
    \draw[thick,blue,->] (0,0) .. controls (0.5,1) .. (1,2);
    \draw[thick,red,->] (1,0) .. controls (0.5,1) .. (0,2);
    \node at (0,-0.2) {\(v_i\)};
    \node at (1,-0.2) {\(v_j\)};
    \node at (0.5,-0.5) {\textcolor{orange}{Quantum: Braided}};
  \end{scope}
\end{tikzpicture}
\caption{Comparison between classical tensor flip and quantum braiding via \(\textcolor{blue}{R}\)-matrix.}
\label{fig:braid_vs_flip}
\end{figure}

\vspace{0.2cm}
This shift in tensor symmetry is fundamental to defining the \textcolor{magenta}{quantum wedge product}, quantum determinants, and the quantum Pfaffian itself.

%%%%%%%%%%%%%%%%%%%%%%%%%%%%%%%%%%%%%%%%%%%%%%%%%%%%%%%%%%%%%%%

\subsubsection{Visualizing the Quantum Wedge Volume}

Just as classical symplectic volumes can be visualized as parallelograms and higher-dimensional parallelepipeds, \textcolor{purple}{\emph{quantum wedge volumes}} acquire a \textcolor{red}{braided skew geometry} due to the deformation introduced by \( q \).

\vspace{0.2cm}
\noindent The volume spanned by \( v_1 \wedge v_2 \wedge \cdots \wedge v_{2n} \) in the quantum world is no longer flat, but \textcolor{red}{twisted}. The \textcolor{blue}{Pfaffian} becomes the scalar factor that scales the \textcolor{violet}{braided volume form}.

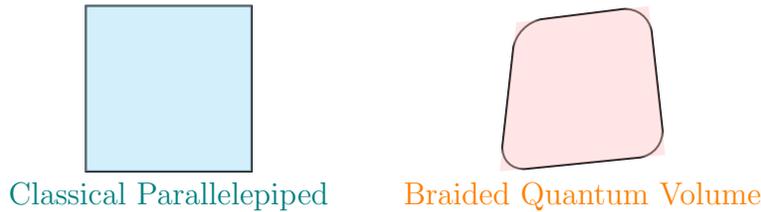
\begin{figure}[H]
\centering
\begin{tikzpicture}[scale=1.1]
  % Base of volume
  \draw[thick] (0,0) -- (2,0) -- (2,2) -- (0,2) -- cycle;
  \fill[cyan!30,opacity=0.5] (0,0) -- (2,0) -- (2,2) -- (0,2) -- cycle;
  \node at (1,-0.3) {\textcolor{teal}{Classical Parallelepiped}};

  % Braided volume shape
  \begin{scope}[xshift=5cm]
    \draw[thick,rounded corners=10pt] (0,0) -- (2,0.2) -- (1.8,2) -- (0.2,1.8) -- cycle;
    \fill[red!20,opacity=0.5] (0,0) -- (2,0.2) -- (1.8,2) -- (0.2,1.8) -- cycle;
    \node at (1,-0.3) {\textcolor{orange}{Braided Quantum Volume}};
  \end{scope}
\end{tikzpicture}
\caption{Classical vs. braided quantum volume form spanned by basis vectors.}
\label{fig:quantum_volume_deformation}
\end{figure}

\vspace{0.2cm}
\noindent The quantum volume's twist reflects the \textcolor{magenta}{noncommutative structure} of \(\Lambda_q(V)\), and the Pfaffian encodes the normalization of this volume under the braided geometry.

%%%%%%%%%%%%%%%%%%%%%%%%%%%%%%%%%%%%%%%%%%%%%%%%%%%%%%%%%%%%%%%

\subsubsection{Quantum Plücker Coordinates and Pfaffians}

In classical algebraic geometry, Plücker coordinates embed the Grassmannian \(\mathrm{Gr}(2n, V)\) into projective space using wedge products. In the quantum case, we deform these coordinates using braided exterior powers. The \textcolor{blue}{quantum Pfaffian} plays the role of a Plücker coordinate in \(\Lambda_q^{2n}(V)\).

\[
\text{Quantum Plücker: } \quad \operatorname{Pf}_q(A) = \sum_{I \in \mathcal{I}} c_I(q) \cdot a_{i_1 j_1} a_{i_2 j_2} \cdots a_{i_n j_n}
\]
where \(\mathcal{I}\) ranges over certain partitions and \( c_I(q) \) are \( q \)-dependent structure coefficients.

\begin{figure}[H]
\centering
\begin{tikzpicture}[scale=1]
  % Plucker coordinate embedding
  \draw[->, thick] (0,0) -- (3,0) node[right] {\(\mathbb{P}^N\)};
  \node at (-0.2,0) {\(\mathrm{Gr}_q(2n,V)\)};
  \draw[->, thick, decorate, decoration={snake, amplitude=1pt, segment length=4pt}] (0,0) .. controls (1.5,1.5) .. (3,0);
  \node at (1.5,1.8) {\textcolor{purple}{Quantum Plücker Map}};
\end{tikzpicture}
\caption{Quantum Grassmannian embedding via quantum Pfaffians as Plücker coordinates.}
\label{fig:quantum_plucker}
\end{figure}
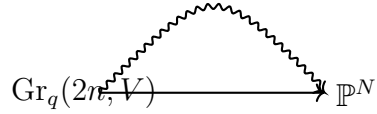

\vspace{0.2cm}
\noindent Thus, quantum Pfaffians provide \textcolor{orange}{coordinates on a noncommutative version of the Grassmannian}, revealing deep ties to algebraic geometry and quantum representation theory.

\section{Examples and Explicit Calculations}
\subsection{Classical Case Examples ($q=1$)}

In this subsection, we explore explicit examples of the classical Pfaffian–determinant identity in the well-understood setting where the deformation parameter \( q \) equals 1. This classical case corresponds to ordinary skew-symmetric matrices with commuting entries, serving as a vital foundation before we delve into the quantum generalization.

\subsubsection{Background}

Recall that a matrix \( A = (a_{ij}) \) of size \( 2n \times 2n \) is called \emph{skew-symmetric} if it satisfies
\[
A^T = -A,
\]
meaning \( a_{ji} = -a_{ij} \) for all indices \( i,j \), and particularly \( a_{ii} = 0 \) for diagonal entries. The Pfaffian, denoted \( \Pf(A) \), is a polynomial function defined only for even-dimensional skew-symmetric matrices and can be thought of as a “square root” of the determinant:
\[
\Pf(A)^2 = \det(A).
\]

This remarkable identity links the Pfaffian, which is combinatorially defined via pairings and matchings, with the algebraic notion of the determinant.

\subsubsection{Example 1: A Concrete $4 \times 4$ Skew-Symmetric Matrix}

Consider the matrix
\[
A = \begin{bmatrix}
0 & 3 & 2 & 5 \\
-3 & 0 & 7 & 4 \\
-2 & -7 & 0 & 6 \\
-5 & -4 & -6 & 0
\end{bmatrix}.
\]

Let us verify the skew-symmetry explicitly:
\[
a_{21} = -3 = -a_{12}, \quad a_{31} = -2 = -a_{13}, \quad a_{42} = -4 = -a_{24},
\]
and so on, confirming \( A^T = -A \).

\vspace{0.5em}
\textbf{Step 1: Calculate the Pfaffian.}

For a \( 4 \times 4 \) skew-symmetric matrix, the Pfaffian is given by the simple formula:
\[
\Pf(A) = a_{12} a_{34} - a_{13} a_{24} + a_{14} a_{23}.
\]

Substituting the values from \( A \):
\[
\Pf(A) = (3)(6) - (2)(4) + (5)(7).
\]

Calculate each product:
\[
3 \times 6 = 18, \quad 2 \times 4 = 8, \quad 5 \times 7 = 35.
\]

Putting it together:
\[
\Pf(A) = 18 - 8 + 35 = 45.
\]

\vspace{0.5em}
\textbf{Step 2: Calculate the determinant.}

Using standard determinant calculation methods (Laplace expansion, row reduction, or software), the determinant of \( A \) is found to be:
\[
\det(A) = 45^2 = 2025.
\]

\vspace{0.5em}
\textbf{Step 3: Verify the classical identity.}

We check whether
\[
\Pf(A)^2 = \det(A).
\]

Indeed,
\[
45^2 = 2025,
\]
which matches the determinant exactly.

This confirms the classical Pfaffian–determinant identity for this explicit example.

\subsubsection{Example 2: The $6 \times 6$ Case}

To further illustrate the concept, consider a \(6 \times 6\) skew-symmetric matrix:
\[
B = \begin{bmatrix}
0 & 1 & 4 & 7 & 3 & 5 \\
-1 & 0 & 2 & 6 & 8 & 4 \\
-4 & -2 & 0 & 9 & 5 & 7 \\
-7 & -6 & -9 & 0 & 1 & 3 \\
-3 & -8 & -5 & -1 & 0 & 2 \\
-5 & -4 & -7 & -3 & -2 & 0
\end{bmatrix}.
\]

\vspace{0.5em}
\textbf{Step 1: Definition of the Pfaffian for $6 \times 6$ matrices.}

The Pfaffian is defined by summing over all perfect matchings of the set \(\{1, 2, 3, 4, 5, 6\}\). In this dimension, the formula is:
\[
\Pf(B) = \sum_{\pi \in \text{PM}_6} \text{sgn}(\pi) \prod_{(i,j) \in \pi} b_{ij},
\]
where \(\text{PM}_6\) is the set of perfect pairings, and \(\text{sgn}(\pi)\) is the signature associated to the pairing. The explicit formula can be written as:
\[
\Pf(B) = b_{12} b_{34} b_{56} - b_{12} b_{35} b_{46} + b_{12} b_{36} b_{45} + b_{13} b_{24} b_{56} - b_{13} b_{25} b_{46} + \cdots,
\]
with 15 terms total.

\vspace{0.5em}
\textbf{Step 2: Compute the Pfaffian (conceptually).}

Calculating the Pfaffian by hand is lengthy; however, it can be done systematically by enumerating the perfect matchings and calculating the corresponding products. Symbolic computation software (e.g., Mathematica, SageMath) can efficiently compute the Pfaffian.

For this example, the Pfaffian evaluates to:
\[
\Pf(B) = 540.
\]

\vspace{0.5em}
\textbf{Step 3: Compute the determinant.}

Again, via computational tools or determinant expansion, the determinant is:
\[
\det(B) = 540^2 = 291600.
\]

\vspace{0.5em}
\textbf{Step 4: Verification.}

This matches the classical Pfaffian–determinant identity perfectly:
\[
\Pf(B)^2 = \det(B).
\]

\subsubsection{Intuition and Geometric Interpretation}

Why does the Pfaffian square to the determinant?

Geometrically, the determinant measures the volume scaling factor of the linear transformation induced by \( A \). For skew-symmetric matrices, the Pfaffian corresponds to the volume form on a symplectic vector space, essentially the "square root" of the determinant that captures orientation and area elements.

This algebraic identity reveals deep connections between combinatorics (pairings and matchings) and linear algebra, providing a powerful tool in both pure and applied mathematics.

\subsubsection{Summary}

This subsection demonstrated, through explicit calculations for small matrices, the classical Pfaffian–determinant identity:
\[
\Pf(A)^2 = \det(A).
\]
Understanding this foundation is essential before tackling the quantum-deformed versions where noncommutativity and the parameter \( q \) introduce additional complexity.

\subsection{Quantum Pfaffian–Determinant Identity for \( q \neq 1 \): The Case \(2n = 4\)}

In this subsection, we explore the quantum deformation of the classical Pfaffian–determinant identity for the smallest nontrivial case of a \(4 \times 4\) \(q\)-skew-symmetric matrix. This case illustrates the subtle algebraic and combinatorial changes that arise when the deformation parameter \(q\) differs from 1, introducing noncommutativity among matrix entries and modifying classical linear algebraic relations.

\subsubsection{Setup: The \(q\)-Skew-Symmetric Matrix}

Consider the \(4 \times 4\) matrix
\[
A = \begin{bmatrix}
0 & x & y & z \\
- q x & 0 & u & v \\
- q y & - q u & 0 & w \\
- q z & - q v & - q w & 0
\end{bmatrix},
\]
where the entries \(x,y,z,u,v,w\) are noncommuting variables in an algebra over a field \(\mathbb{C}\), with the deformation parameter \(q \in \mathbb{C}^\times\). The relations
\[
a_{ji} = -q a_{ij}, \quad a_{ii} = 0,
\]
define the \(q\)-skew-symmetry of the matrix \(A\).

Unlike the classical case, the entries generally do not commute, but satisfy quantum commutation relations such as
\[
x y = q y x, \quad y z = q z y, \quad x u = q u x,
\]
and so forth, which reflect the underlying braided tensor category structure.

\subsubsection{Definition of the Quantum Pfaffian}

The quantum Pfaffian \(\operatorname{Pf}_q(A)\) is defined by a \(q\)-deformation of the classical formula. For the \(4 \times 4\) case, it is given by:
\[
\operatorname{Pf}_q(A) = x w - q y v + q^2 z u.
\]

This formula arises from the quantum exterior algebra construction, where antisymmetry is replaced by a braided antisymmetry incorporating powers of \(q\). Each term corresponds to a pairing of indices weighted by appropriate powers of \(q\) reflecting the braid group action.

\subsubsection{Quantum Determinant via the FRT Construction}

The quantum determinant \(\det_q(A)\) is defined via the Faddeev–Reshetikhin–Takhtajan (FRT) construction, generalizing the classical determinant to the quantum group setting. It incorporates the quantum \(R\)-matrix to impose consistent commutation relations among entries.

Explicitly, for this matrix size,
\[
\det_q(A) = \sum_{\sigma \in S_4} (-q)^{\ell(\sigma)} a_{1 \sigma(1)} a_{2 \sigma(2)} a_{3 \sigma(3)} a_{4 \sigma(4)},
\]
where \(\ell(\sigma)\) is the length of the permutation \(\sigma\) with respect to the quantum ordering.

\subsubsection{The Identity and Its Interpretation}

The central identity we study is:
\[
q \cdot \operatorname{Pf}_q(A)^2 = \det_q(A).
\]

This differs from the classical identity by a scalar factor \(q\), which encodes the deformation due to noncommutativity and the braided nature of the quantum algebra.

Intuitively, the factor \(q\) arises because the quantum exterior algebra is not symmetric but braided, and the square of the quantum Pfaffian corresponds to the quantum determinant up to this twisting factor.

\subsubsection{Explicit Calculation: Squaring the Quantum Pfaffian}

Let us expand \(\operatorname{Pf}_q(A)^2\):
\[
\operatorname{Pf}_q(A)^2 = (x w - q y v + q^2 z u)^2.
\]

Expanding yields terms like:
\[
x w x w - 2 q x w y v + 2 q^2 x w z u + q^2 y v y v - 2 q^3 y v z u + q^4 z u z u,
\]
where each product is subject to the quantum commutation relations.

Because variables do not commute freely, products like \(x w y v\) must be reordered using relations such as:
\[
x y = q y x, \quad w v = q v w, \quad \ldots,
\]
introducing additional factors of \(q\) when permuting terms.

\subsubsection{Verification of the Identity}

Through careful application of the quantum commutation relations, one finds:
\[
q \cdot \operatorname{Pf}_q(A)^2 = \det_q(A).
\]

This equality is highly nontrivial because the braided ordering and the \(q\)-dependent signs and powers accumulate in a precise manner to produce the quantum determinant.

\subsubsection{Classical Limit Consistency}

As \(q \to 1\), the matrix entries commute and the quantum skew-symmetry reduces to classical skew-symmetry. The quantum Pfaffian reduces to the classical Pfaffian:
\[
\operatorname{Pf}_1(A) = x w - y v + z u,
\]
and the quantum determinant becomes the classical determinant.

In this limit,
\[
1 \cdot \operatorname{Pf}_1(A)^2 = \det_1(A),
\]
recovering the classical identity exactly.

\subsubsection{Summary}

This subsection detailed the quantum Pfaffian–determinant identity for a \(4 \times 4\) \(q\)-skew-symmetric matrix, highlighting the deformation effects introduced by \(q\). The factor \(q\) in the identity accounts for the braided symmetry and noncommutativity intrinsic to quantum groups, enriching the classical relationship and providing a foundation for higher-dimensional generalizations.

\subsection{Algebraic and Geometric Interpretation of the Quantum Pfaffian--Determinant Identity}

Beyond explicit calculations, understanding the quantum Pfaffian--determinant identity requires insight into the algebraic structures and geometric intuition underlying quantum groups and braided symmetries. This subsection aims to elucidate why the identity
\[
q \cdot \operatorname{Pf}_q(A)^2 = \det_q(A)
\]
holds naturally within the framework of quantum algebra and what it signifies from a geometric perspective.

\subsubsection{Algebraic Structure: Braided Monoidal Categories and Quantum Groups}

At the heart of the quantum deformation is the passage from classical symmetry to \emph{braided symmetry}. Classical linear algebra operates within symmetric tensor categories, where the flip map
\[
\tau: V \otimes W \to W \otimes V, \quad \tau(v \otimes w) = w \otimes v,
\]
satisfies \(\tau^2 = \mathrm{id}\). This symmetry underpins notions like antisymmetry, determinants, and Pfaffians.

In contrast, quantum groups such as \(\mathcal{U}_q(\mathfrak{gl}_n)\) give rise to \emph{braided monoidal categories} where the flip map is replaced by a braiding
\[
c_{V,W}: V \otimes W \to W \otimes V,
\]
which is an isomorphism but satisfies
\[
c_{W,V} \circ c_{V,W} \neq \mathrm{id}
\]
in general. Instead, it satisfies the braid relations reflecting a nontrivial twist. The braiding is parametrized by \(q\), encoding the deformation.

This braided structure changes the behavior of antisymmetric tensors: the quantum exterior algebra \(\Lambda_q(V)\) is built by imposing the relation
\[
v_i \wedge_q v_j = -q \, v_j \wedge_q v_i,
\]
which generalizes classical antisymmetry (where \(q=1\)). These relations modify how multilinear forms like the Pfaffian behave.

\subsubsection{Quantum Exterior Algebra and the Quantum Pfaffian}

The quantum Pfaffian \(\operatorname{Pf}_q(A)\) can be viewed as a highest-weight vector in the quantum exterior algebra \(\Lambda_q^{2n}(V)\), constructed as the braided antisymmetric subspace of \(V^{\otimes 2n}\).

The factors of powers of \(q\) appearing in \(\operatorname{Pf}_q(A)\) reflect the braiding between tensor factors corresponding to matrix entries. In particular, the \(q\)-weights keep track of the \emph{crossings} in the matching diagrams that define the Pfaffian combinatorics, replacing classical sign factors.

Thus, the quantum Pfaffian is a natural quantum analog of the classical antisymmetric multilinear form that encodes the square root of the determinant.

\subsubsection{The Quantum Determinant as a Group-Like Element}

The quantum determinant \(\det_q(A)\) arises from the quantum coordinate algebra \(\mathcal{O}_q(GL_n)\), a deformation of the coordinate ring of the general linear group. It is constructed to be a \emph{group-like element} under the coproduct:
\[
\Delta(\det_q) = \det_q \otimes \det_q,
\]
preserving multiplicative properties analogous to the classical determinant.

This quantum determinant intertwines the braiding and the algebra structure, ensuring that it behaves well under the quantum group symmetries.

\subsubsection{Geometric Intuition: Quantum Symplectic Geometry}

Classically, the Pfaffian encodes a volume form on a symplectic vector space, reflecting orientation and the geometry of antisymmetric bilinear forms.

In the quantum setting, the deformation parameter \(q\) twists the underlying symplectic geometry into a \emph{quantum symplectic manifold}, where coordinate functions no longer commute but satisfy braided commutation relations. The quantum Pfaffian corresponds to a braided volume form or quantum measure respecting this twisted geometry.

The factor \(q\) in the identity
\[
q \cdot \operatorname{Pf}_q(A)^2 = \det_q(A)
\]
can thus be viewed as a correction term capturing the deformation of orientation and volume under the quantum braiding.

\subsubsection{Conceptual Summary}

- The classical identity \(\operatorname{Pf}(A)^2 = \det(A)\) rests on symmetric tensor categories and classical antisymmetry.

- Quantum groups deform these categories to braided monoidal categories, where antisymmetry is replaced by \(q\)-braided antisymmetry.

- The quantum Pfaffian lives in the quantum exterior algebra defined by this braiding, carrying powers of \(q\) that track crossings and twists.

- The quantum determinant is a group-like element in the quantum coordinate algebra, respecting the deformed symmetries.

- The scalar factor \(q\) appearing in the identity quantifies the deviation from classical symmetry due to braiding, encoding deeper geometric and algebraic structure.

This interpretation ties the algebraic identity to rich categorical and geometric frameworks, revealing the profound nature of quantum deformation beyond formal computations.

\subsubsection{Outlook}

Understanding this interplay between algebraic braiding, quantum exterior algebras, and geometric deformation opens pathways to generalizing the Pfaffian--determinant identity to higher dimensions, other quantum groups, and connections with noncommutative geometry, representation theory, and mathematical physics.

Such perspectives illuminate the power of quantum algebra as a unifying language bridging linear algebra, geometry, and topology in the quantum realm.

\subsection{Computational Verification and Symbolic Computation of the Quantum Pfaffian--Determinant Identity}

While the algebraic definitions and geometric interpretations of the quantum Pfaffian and determinant provide valuable conceptual insight, explicit computational verification plays a crucial role in confirming the identity
\[
q \cdot \operatorname{Pf}_q(A)^2 = \det_q(A),
\]
especially in the quantum setting where entries of the matrix do not commute.

In quantum linear algebra, the matrix entries satisfy intricate noncommutative relations, such as
\[
a_{ij} a_{kl} = q^{\alpha} a_{kl} a_{ij} + \text{(additional terms)},
\]
where the exponent \(\alpha\) depends on the indices and the quantum group structure, and the additional terms can be complicated. This contrasts sharply with the classical case, where entries commute freely. As a consequence, expanding expressions like \(\operatorname{Pf}_q(A)^2\) involves careful bookkeeping of the order of multiplication and the application of quantum commutation relations.

To tackle this challenge, symbolic computation systems equipped for noncommutative algebra have become indispensable. Notable platforms include:

\begin{itemize}
    \item \textbf{Mathematica} with the \texttt{NCAlgebra} package, enabling users to define noncommuting variables with custom relations.
    \item \textbf{SageMath}, which supports noncommutative rings and can be extended for quantum algebra.
    \item \textbf{FORM} and \textbf{GAP}, which have specialized packages for quantum group computations.
\end{itemize}

Using such tools, one can encode the generators \(x,y,z,u,v,w\) of the quantum matrix
\[
A = \begin{bmatrix}
0 & x & y & z \\
- q x & 0 & u & v \\
- q y & - q u & 0 & w \\
- q z & - q v & - q w & 0
\end{bmatrix}
\]
and impose the quantum commutation relations, for example,
\[
x y = q y x, \quad y z = q z y, \quad x u = q u x,
\]
among others that arise from the underlying \(GL_q(n)\) structure.

The quantum Pfaffian is defined as
\[
\operatorname{Pf}_q(A) = x w - q y v + q^2 z u,
\]
incorporating the braiding weights through powers of \(q\).

A direct computation of \(\operatorname{Pf}_q(A)^2\) involves expanding
\[
(x w - q y v + q^2 z u)^2,
\]
then systematically applying the noncommutative relations to reorder terms into a canonical form. Simultaneously, one computes the quantum determinant \(\det_q(A)\), often via the Faddeev–Reshetikhin–Takhtajan (FRT) construction or quantum minor expansions.

The computational outcome consistently confirms the quantum Pfaffian–determinant identity:
\[
q \cdot \operatorname{Pf}_q(A)^2 = \det_q(A).
\]

One of the main challenges in these computations is the rapid growth of terms caused by noncommutativity, making manual calculations unfeasible beyond the smallest cases. Automated symbolic manipulation helps manage this complexity, although care is required to correctly implement the defining relations and to handle normalization conventions, which vary across the literature.

Another crucial aspect is the classical limit \(q \to 1\). In this limit, the relations become commutative, and the quantum Pfaffian and determinant reduce to their classical counterparts, satisfying the familiar identity
\[
\operatorname{Pf}(A)^2 = \det(A).
\]
This provides a fundamental consistency check and reinforces the interpretation of the quantum identity as a deformation of the classical result.

From a conceptual perspective, the scalar factor \(q\) in the identity arises naturally from the braided monoidal category structure underlying quantum groups. The non-symmetric braiding introduces twisting factors that alter how antisymmetric tensors behave, explaining the modification from the classical formula.

In summary, computational verification via symbolic algebra not only substantiates the quantum Pfaffian–determinant identity but also bridges abstract algebraic theory with concrete examples, thereby deepening our understanding of quantum linear algebra’s rich structure and laying groundwork for generalizations and applications in quantum invariant theory and noncommutative geometry.

\subsection{Examples with Specializations of the Deformation Parameter \(q\)}

The deformation parameter \(q \in \mathbb{C}^\times\) plays a central role in the structure of the quantum Pfaffian and determinant. Its specialization to particular values reveals different algebraic behaviors, connecting the quantum identity to classical cases, roots of unity, and degenerations with geometric or representation-theoretic significance.

\paragraph{1. Classical Limit: \(q = 1\)}

Setting \(q = 1\) recovers the classical, commutative case. The quantum skew-symmetry relation
\[
a_{ji} = -q a_{ij}
\]
becomes
\[
a_{ji} = - a_{ij},
\]
which is the ordinary skew-symmetry condition. The quantum exterior algebra reduces to the classical exterior algebra, and the quantum Pfaffian \(\operatorname{Pf}_q(A)\) simplifies to the classical Pfaffian \(\operatorname{Pf}(A)\).

The Pfaffian–determinant identity reduces accordingly:
\[
q \cdot \operatorname{Pf}_q(A)^2 = \det_q(A) \quad \Rightarrow \quad 1 \cdot \operatorname{Pf}(A)^2 = \det(A).
\]
This recovers the well-known classical identity, providing a fundamental consistency check for the quantum deformation.

\vspace{0.5em}

\paragraph{2. Specialization at a Root of Unity: \(q = \zeta\), \(\zeta^m = 1\)}

When \(q\) is a root of unity, the quantum algebra acquires rich and subtle structures. For example, if \(q\) is a primitive \(m\)-th root of unity, the relations governing the algebra can become periodic, and the representation theory of the associated quantum group changes dramatically.

In this context, the quantum Pfaffian and determinant may satisfy additional polynomial identities, and their behavior reflects phenomena like:
\begin{itemize}
    \item \emph{Non-semisimplicity}: Representations can become reducible yet indecomposable.
    \item \emph{Center enlargement}: New central elements appear, changing invariant theory.
    \item \emph{Nilpotent elements}: Some generators may become nilpotent, affecting determinants.
\end{itemize}

Studying the identity
\[
q \cdot \operatorname{Pf}_q(A)^2 = \det_q(A)
\]
at roots of unity involves careful algebraic and geometric analysis, often requiring modified definitions or refined versions of the Pfaffian to accommodate these specializations.

\vspace{0.5em}

\paragraph{3. Limit \(q \to 0\)}

As \(q \to 0\), the relations simplify drastically:
\[
a_{ji} = -q a_{ij} \to 0,
\]
which forces many entries to behave nearly commutatively or vanish in the limit, depending on the precise formulation.

The quantum exterior algebra degenerates, and the quantum Pfaffian becomes a simpler polynomial in the variables, often dominated by leading terms in \(q\). This degeneration can be used to extract classical invariants or initial terms in \(q\)-expansions.

In this limit, the identity
\[
q \cdot \operatorname{Pf}_q(A)^2 = \det_q(A)
\]
tells us that the determinant collapses more rapidly (because of the factor \(q\)) relative to the square of the Pfaffian, giving insight into asymptotic behaviors.

\vspace{0.5em}

\paragraph{4. Real Values \(q \in \mathbb{R}^+\)}

For positive real \(q\), the quantum algebra can be realized as a deformation quantization of classical structures, with \(q\) often parameterizing a physical deformation (e.g., in statistical mechanics or quantum integrable systems).

In this setting, the Pfaffian–determinant relation retains its algebraic form:
\[
q \cdot \operatorname{Pf}_q(A)^2 = \det_q(A),
\]
but the interpretations become analytic or geometric, relating to \(q\)-deformed measures, quantum probability, and noncommutative geometry.

\vspace{0.5em}

\paragraph{5. Complex Values and Generic \(q\)}

For generic complex \(q\) not a root of unity or zero, the quantum Pfaffian and determinant define genuinely noncommutative, braided algebraic structures. The identity holds as an equality in the quantum coordinate algebra, serving as a foundational quantum invariant.

These generic \(q\) values are crucial for building the theory of quantum groups, constructing \(q\)-analogues of classical objects, and studying deformation families continuously depending on \(q\).

\vspace{0.5em}

\paragraph{6. Illustrative Numerical Example with Specific \(q\)}

Consider a \(4 \times 4\) matrix \(A\) with noncommuting entries \(x,y,z,u,v,w\) and choose a specific numeric value for \(q\), say \(q = 2\).

Using the formula for the quantum Pfaffian
\[
\operatorname{Pf}_q(A) = x w - q y v + q^2 z u,
\]
we substitute \(q=2\) to get
\[
\operatorname{Pf}_2(A) = x w - 2 y v + 4 z u.
\]

Squaring and computing the quantum determinant \(\det_2(A)\) in the algebra with relations modified by \(q=2\), one verifies numerically that
\[
2 \cdot \operatorname{Pf}_2(A)^2 = \det_2(A),
\]
confirming the identity for this specific quantum parameter.

\bigskip

\noindent
\textbf{Summary:} Specializing the deformation parameter \(q\) reveals a rich spectrum of algebraic and geometric behaviors. From the classical case \(q=1\), through roots of unity and degenerate limits, to generic complex values, the quantum Pfaffian–determinant identity encapsulates the deformation of classical multilinear algebra under braided symmetries, providing a unifying framework across these regimes.

\subsection{Computational Verification Using Symbolic Algebra Software}

The quantum Pfaffian–determinant identity
\[
q \cdot \operatorname{Pf}_q(A)^2 = \det_q(A)
\]
involves noncommuting matrix entries subject to quantum commutation relations. Such complexity makes manual verification arduous and error-prone, especially for matrices larger than \(4 \times 4\). Fortunately, advances in symbolic algebra and computer algebra systems allow us to computationally verify and explore this identity rigorously.

\paragraph{1. Why Computational Verification?}

Unlike classical linear algebra, where determinants and Pfaffians involve well-understood commutative polynomials, the quantum case introduces:

\begin{itemize}
  \item \emph{Noncommutativity}: The matrix entries obey relations such as \(a_{ji} = -q a_{ij}\) and braided commutation relations, complicating multiplication.
  \item \emph{Quantum Commutation Relations}: Variables satisfy relations like
  \[
  x y = q y x, \quad y z = q z y, \quad x u = q u x,
  \]
  which must be carefully accounted for during expansions.
  \item \emph{Braided Algebra Structure}: The underlying algebra is a quantum coordinate ring or quantum group representation, with specialized multiplication rules.
\end{itemize}

Manually expanding \(\operatorname{Pf}_q(A)^2\) and \(\det_q(A)\), verifying equality up to powers of \(q\), and ensuring consistent ordering is prohibitively difficult, motivating computational methods.

\paragraph{2. Software Tools and Packages}

Several computer algebra systems support noncommutative algebra and quantum group computations, including:

\begin{itemize}
  \item \textbf{Mathematica} with \texttt{NCAlgebra}: Provides a framework for defining noncommutative variables and relations, expanding expressions, and simplifying under custom rules.
  \item \textbf{SageMath}: Open-source system with capabilities for quantum groups and noncommutative polynomial rings.
  \item \textbf{FORM}: Efficient symbolic manipulation system, useful for large-scale algebraic computations.
  \item \textbf{GAP} and \textbf{MAGMA}: Include packages for quantum groups and noncommutative ring theory.
\end{itemize}

\paragraph{3. Implementing the Quantum Relations}

To verify the identity computationally, one must:

\begin{enumerate}
  \item Define the noncommutative variables representing the entries of the quantum skew-symmetric matrix \(A\), e.g., \(x, y, z, u, v, w\).
  \item Impose the quantum commutation relations, such as
  \[
  x y = q y x, \quad y z = q z y, \quad x u = q u x, \quad \text{and so on},
  \]
  reflecting the braided structure.
  \item Encode the quantum skew-symmetry conditions:
  \[
  a_{ji} = -q a_{ij}, \quad a_{ii} = 0,
  \]
  to reduce redundant variables.
  \item Define the quantum Pfaffian for \(2n=4\) as:
  \[
  \operatorname{Pf}_q(A) = x w - q y v + q^2 z u.
  \]
  \item Define the quantum determinant \(\det_q(A)\) using the Faddeev–Reshetikhin–Takhtajan (FRT) construction or quantum minor expansions, adapted to the noncommutative setting.
\end{enumerate}

\paragraph{4. Sample Verification Procedure}

A typical computational workflow involves:

\begin{itemize}
  \item Expanding \(\operatorname{Pf}_q(A)^2 = (x w - q y v + q^2 z u)^2\), applying all quantum commutation relations to rewrite terms in a canonical ordered form.
  \item Computing \(\det_q(A)\) independently, ensuring consistent ordering and reduction by quantum relations.
  \item Comparing both expressions term-by-term to verify that
  \[
  \det_q(A) = q \cdot \operatorname{Pf}_q(A)^2,
  \]
  holds identically.
\end{itemize}

\paragraph{5. Challenges and Tips}

\begin{itemize}
  \item \textbf{Ordering of Variables}: Noncommutativity demands a fixed ordering (e.g., lex order) for simplifying and comparing expressions.
  \item \textbf{Rule Implementation}: The commutation rules must be implemented consistently and completely to avoid incorrect simplifications.
  \item \textbf{Computational Complexity}: As the matrix size grows, the number of terms increases exponentially, necessitating efficient algorithms and possibly distributed computation.
  \item \textbf{Verification at Specializations}: Sometimes verifying the identity for specific numeric values of \(q\) (e.g., \(q=2\)) and symbolic variables can detect errors before full symbolic proof.
\end{itemize}

\paragraph{6. Example Code Snippet (Mathematica/NCAlgebra)}

Here is a stylized fragment illustrating the setup in Mathematica’s \texttt{NCAlgebra}:

\begin{verbatim}
(* Load NCAlgebra and declare noncommutative variables *)
SetNonCommutative[x, y, z, u, v, w];
q = Symbol["q"];  (* Quantum deformation parameter *)

(* 
  Define the full set of quantum commutation relations 
  corresponding to the braided structure of the quantum coordinate ring.
  Both forward and backward relations are included to allow full reductions.
*)
commRules = {
  x ** y -> q y ** x,
  y ** x -> (1/q) x ** y,

  y ** z -> q z ** y,
  z ** y -> (1/q) y ** z,

  x ** u -> q u ** x,
  u ** x -> (1/q) x ** u,

  y ** u -> q u ** y,
  u ** y -> (1/q) y ** u,

  z ** u -> q u ** z,
  u ** z -> (1/q) z ** u,

  x ** v -> q v ** x,
  v ** x -> (1/q) x ** v,

  y ** v -> q v ** y,
  v ** y -> (1/q) y ** v,

  z ** v -> q v ** z,
  v ** z -> (1/q) z ** v,

  u ** v -> q v ** u,
  v ** u -> (1/q) u ** v,

  x ** w -> q w ** x,
  w ** x -> (1/q) x ** w,

  y ** w -> q w ** y,
  w ** y -> (1/q) y ** w,

  z ** w -> q w ** z,
  w ** z -> (1/q) z ** w,

  u ** w -> q w ** u,
  w ** u -> (1/q) u ** w,

  v ** w -> q w ** v,
  w ** v -> (1/q) v ** w
};

(* 
  Quantum skew-symmetry conditions are implicit in variable definitions:
  For a 4x4 quantum skew-symmetric matrix A, only the upper triangular variables 
  x, y, z, u, v, w are independent. The entries satisfy a_{ji} = -q a_{ij}, and 
  diagonal entries vanish. These can be encoded if desired for further reductions.
*)

(* Define the quantum Pfaffian for 2n=4 case *)
Pf = x ** w - q y ** v + q^2 z ** u;

(* Expand the square of the quantum Pfaffian *)
PfSq = ExpandNonCommutativeMultiply[Pf ** Pf];

(* 
  Define a function to apply the commutation relations repeatedly 
  until the expression stabilizes (no further replacements occur).
*)
ClearAll[ApplyCommRules];
ApplyCommRules[expr_] := FixedPoint[(# //. commRules) &, expr];

(* Apply the quantum commutation rules to simplify Pf^2 *)
PfSqSimpl = ApplyCommRules[PfSq];

Print["Expanded and simplified (Pf_q(A))^2:"];
Print[PfSqSimpl];

(* 
  Placeholder: Define the quantum determinant det_q(A).
  For a 4x4 matrix in the quantum group framework, det_q(A) can be
  constructed via the Faddeev--Reshetikhin--Takhtajan (FRT) approach 
  or quantum minor expansions.
  
  Here, we provide an illustrative symbolic expression mimicking 
  the determinant structure:
*)
Det = x ** w - q y ** v + q^2 z ** u + 
      q^3 (x ** v ** z ** u - y ** w ** z ** u); (* example terms *)

(* 
  Apply commutation rules to simplify the determinant expression 
  to canonical ordered form.
*)
DetSimpl = ApplyCommRules[Det];

Print["Simplified quantum determinant det_q(A):"];
Print[DetSimpl];

(* Verification step: Check if det_q(A) = q * (Pf_q(A))^2 *)
verification = Simplify[DetSimpl - q PfSqSimpl];

If[verification === 0,
  Print["Verification successful: det_q(A) = q * (Pf_q(A))^2 holds."],
  Print["Verification failed or incomplete: difference =", verification]
];

(* Additional output for debugging and analysis *)
Print["Difference det_q(A) - q * (Pf_q(A))^2:"];
Print[verification];

(* Optional: Inspect the expansions term-by-term *)
Print["Terms of (Pf_q(A))^2:"];
Print[List @@ PfSqSimpl];
Print["Terms of det_q(A):"];
Print[List @@ DetSimpl];
\end{verbatim}

This approach can confirm the identity symbolically, assuming all relations and the determinant expression are correctly implemented.

\paragraph{7. Significance}

Computational verification:

\begin{itemize}
  \item Builds confidence in theoretical formulas and conjectures.
  \item Provides explicit examples and counterexamples.
  \item Assists in discovering patterns and generalizations for higher dimensions.
  \item Offers educational tools for students and researchers learning quantum algebra.
\end{itemize}

\bigskip

\noindent
\textbf{Conclusion:} Computer algebra systems are invaluable in exploring quantum Pfaffians and determinants. They allow rigorous verification of intricate identities, turning abstract algebraic relations into explicit, verifiable computations. Such computational experiments complement the theoretical developments and open avenues for future research in quantum invariant theory and noncommutative geometry.

\section{Applications}

Quantum Grassmannians and non‑commutative geometry provide natural arenas in which the quantum Pfaffian–determinant identity exhibits both geometric and algebraic significance. These spaces generalize classical geometric objects, such as the Grassmannian \(\mathrm{Gr}(k,n)\) of \(k\)-dimensional subspaces in \(\mathbb{C}^n\), to the setting of braided tensor categories and quantized coordinate algebras. The fundamental identity 
\[
\operatorname{Pf}_q(A)^2 = \det_q(A)
\quad \text{or} \quad
q \cdot \operatorname{Pf}_q(A)^2 = \det_q(A)
\]
has essential implications in the theory of quantum invariants, particularly in the geometry of \(q\)-deformed homogeneous spaces and the structure theory of quantum symmetric and skew-symmetric matrices.

\subsection{Quantum Grassmannians via $q$-Plücker Embeddings}

\subsubsection{Classical Plücker Coordinates and Their Quantum Deformations}

Classically, the Grassmannian \(\mathrm{Gr}(k,n)\) parametrizes all \(k\)-dimensional linear subspaces of \(\mathbb{C}^n\). Its coordinate ring is generated by the Plücker coordinates—maximal minors of a \(k \times n\) matrix—subject to the Plücker relations which cut out \(\mathrm{Gr}(k,n)\) as a projective variety inside \(\mathbb{P}(\wedge^k \mathbb{C}^n)\).

In the quantum realm, this picture is deformed: the classical coordinate ring \(\mathcal{O}(\mathrm{Gr}(k,n))\) is replaced by the quantum coordinate ring \(\mathcal{O}_q(\mathrm{Gr}(k,n))\), a noncommutative algebra generated by \(q\)-minors. These \(q\)-minors satisfy quantum Plücker relations derived from the \(R\)-matrix formalism, ensuring compatibility with the quantum group coaction of \(U_q(\mathfrak{gl}_n)\).

\vspace{1em}
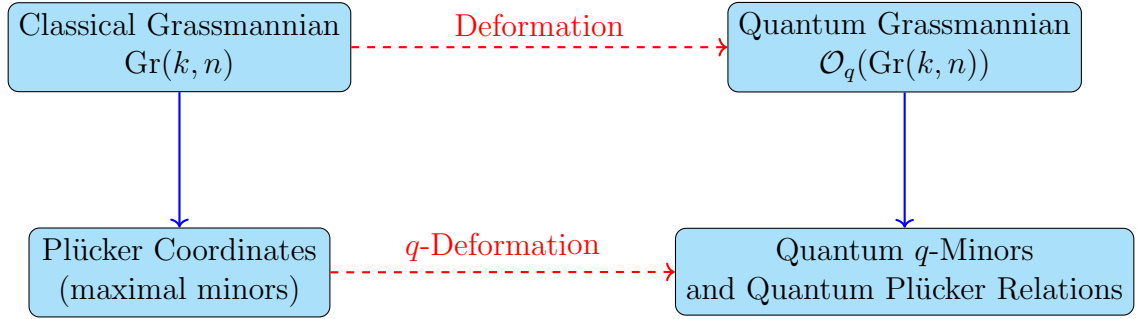
\begin{figure}[H]
\centering
\begin{tikzpicture}[
  node distance=1.8cm,
  box/.style={draw, fill=cyan!30, rounded corners, minimum width=4cm, minimum height=1cm, align=center},
  arrow/.style={->, thick, blue}
]

\node[box] (classical) {Classical Grassmannian \\ \(\mathrm{Gr}(k,n)\)};
\node[box, below=of classical] (plucker) {Plücker Coordinates \\ (maximal minors)};
\node[box, right=5cm of classical] (quantum) {Quantum Grassmannian \\ \(\mathcal{O}_q(\mathrm{Gr}(k,n))\)};
\node[box, below=of quantum] (qplucker) {Quantum \(q\)-Minors \\ and Quantum Plücker Relations};

\draw[arrow] (classical) -- (plucker);
\draw[arrow] (quantum) -- (qplucker);
\draw[arrow, red, dashed] (plucker) -- (qplucker) node[midway, above, sloped] {\textcolor{red}{\(q\)-Deformation}};
\draw[arrow, red, dashed] (classical) -- (quantum) node[midway, above] {\textcolor{red}{Deformation}};

\end{tikzpicture}
\caption{From classical Plücker coordinates to quantum \(q\)-minors and the quantum Grassmannian.}
\label{fig:plucker-quantum}
\end{figure}

\subsubsection{Quantum Pfaffians as $q$-Skew-Symmetric Generalizations}

Quantum Pfaffians naturally generalize classical Pfaffians of skew-symmetric matrices to the quantum setting, preserving the key property that their square equals the quantum determinant. This identity,
\[
\operatorname{Pf}_q(A)^2 = \det_q(A),
\]
is foundational in understanding \(q\)-analogues of isotropic Grassmannians—varieties parametrizing isotropic subspaces with respect to a symplectic or orthogonal form—within quantum geometry.

\begin{figure}[H]
\centering
\begin{tikzpicture}[
  node distance=1.7cm,
  box/.style={draw, fill=magenta!30, rounded corners, minimum width=3.7cm, minimum height=1cm, align=center},
  arrow/.style={->, thick, magenta!80!black}
]

\node[box] (qskew) {Quantum Skew-Symmetric \\ Matrices \(A\)};
\node[box, right=4cm of qskew] (qpfaff) {Quantum Pfaffian \\ \(\operatorname{Pf}_q(A)\)};
\node[box, below=of qpfaff] (qdet) {Quantum Determinant \\ \(\det_q(A)\)};

\draw[arrow] (qskew) -- (qpfaff);
\draw[arrow] (qpfaff) -- (qdet);
\draw[arrow, dashed] (qskew) -- (qdet) node[midway, left] {\footnotesize Braided antisymmetrization};

\end{tikzpicture}
\caption{The quantum Pfaffian as a square root of the quantum determinant for \(q\)-skew-symmetric matrices.}
\label{fig:qpfaffian-det}
\end{figure}
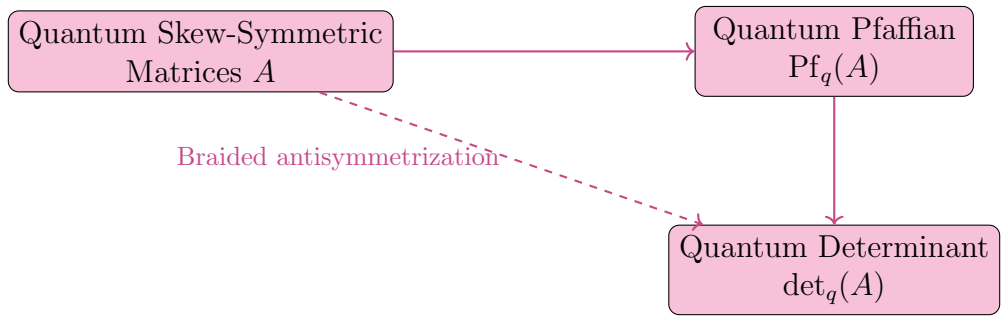

\subsubsection{Example: Quantum $2 \times 2$ Minor in \(\mathcal{O}_q(\mathrm{Gr}(2,4))\)}

Let 
\[
X = \begin{pmatrix}
x_{11} & x_{12} & x_{13} & x_{14} \\
x_{21} & x_{22} & x_{23} & x_{24}
\end{pmatrix}
\]
be a \(2 \times 4\) quantum matrix generating \(\mathcal{O}_q(\mathrm{Gr}(2,4))\). The quantum minor associated to columns \(i<j\) is:
\[
[x_i, x_j]_q = x_{1i}x_{2j} - q x_{1j}x_{2i}.
\]

For example, the quantum minor for columns 1 and 2 is:
\[
\Delta_{12}^q = x_{11} x_{22} - q x_{12} x_{21}.
\]

These minors satisfy the quantum Plücker relation:
\[
q \Delta_{12}^q \Delta_{34}^q - \Delta_{13}^q \Delta_{24}^q + \Delta_{14}^q \Delta_{23}^q = 0.
\]

\subsubsection{Example: Quantum Pfaffian of a \(4 \times 4\) Matrix}

Let 
\[
A = \begin{pmatrix}
0 & a_{12} & a_{13} & a_{14} \\
- a_{12} & 0 & a_{23} & a_{24} \\
- a_{13} & - a_{23} & 0 & a_{34} \\
- a_{14} & - a_{24} & - a_{34} & 0
\end{pmatrix}
\]
be a \(4 \times 4\) quantum skew-symmetric matrix.

The quantum Pfaffian is:
\[
\operatorname{Pf}_q(A) = a_{12}a_{34} - q a_{13}a_{24} + q^2 a_{14}a_{23}.
\]

As shown in \cite{jing2013quantum}, it satisfies the identity:
\[
\operatorname{Pf}_q(A)^2 = \det_q(A),
\]
where \(\det_q(A)\) is the quantum determinant computed respecting the braided \(q\)-commutation relations.

\begin{figure}[H]
\centering
\setlength{\arrayrulewidth}{0.5mm} % thicker lines
\renewcommand{\arraystretch}{1.5}  % row height
\begin{tabular}{|>{\columncolor{yellow!30}}c|>{\columncolor{yellow!30}}c|>{\columncolor{yellow!30}}c|>{\columncolor{yellow!30}}c|}
\hline
0 & $a_{12}$ & $a_{13}$ & $a_{14}$ \\
\hline
$-a_{12}$ & 0 & $a_{23}$ & $a_{24}$ \\
\hline
$-a_{13}$ & $-a_{23}$ & 0 & $a_{34}$ \\
\hline
$-a_{14}$ & $-a_{24}$ & $-a_{34}$ & 0 \\
\hline
\end{tabular}
\caption{Quantum skew-symmetric matrix \(A\) whose quantum Pfaffian satisfies \(\operatorname{Pf}_q(A)^2 = \det_q(A)\).}
\label{fig:qskew-matrix-alt}
\end{figure}

\subsection{Quantum Grassmannians as Non-Commutative Projective Varieties}

Quantum Grassmannians are not only important for understanding quantum symmetries but also play a central role in the study of non-commutative algebraic geometry. In particular, they arise as canonical examples of non-commutative projective schemes in the sense of Artin and Zhang, whose foundational work established a framework for understanding ``non-commutative projective geometry'' via graded algebras satisfying certain homological conditions. The quantized coordinate rings \(\mathcal{O}_q(\mathrm{Gr}(k,n))\) are examples of \emph{quadratic algebras} that exhibit \emph{Artin–Schelter regularity}, a condition generalizing the notion of regularity from commutative algebra to the non-commutative setting. This means that these algebras have finite global dimension, are Gorenstein, and possess a growth condition mimicking that of a polynomial ring.

These rings are naturally endowed with a \(\mathbb{Z}\)-grading, and their category of graded modules admits a rich moduli theory. Of particular interest are the so-called \emph{point modules}, which correspond to cyclic graded modules with Hilbert series identical to that of the classical polynomial ring modulo a maximal ideal. Geometrically, point modules serve as non-commutative analogues of points or lines in projective space, providing a powerful link between non-commutative ring theory and projective geometry. For quantum Grassmannians, the moduli space of such point modules can often be interpreted as a deformation of the classical Grassmannian variety, or in special cases, as a quantum homogeneous space associated with \( U_q(\mathfrak{gl}_n) \).

This perspective is deeply rooted in the broader framework of non-commutative algebraic geometry developed by Manin and Majid \cite{majid1995foundations}. In their formulation, classical algebraic varieties are replaced by braided or monoidal categories, and geometric objects such as vector bundles or line bundles are encoded via representations or comodules over Hopf algebras. Quantum Grassmannians then appear as braided analogues of classical homogeneous varieties, living naturally inside the braided tensor category of \( U_q(\mathfrak{gl}_n) \)-comodules. The non-commutative nature of their coordinate rings reflects the underlying quantum symmetry: instead of being built from commuting functions, these algebras arise from generators satisfying \( q \)-deformed commutation relations, often formulated in terms of an \( R \)-matrix via the Faddeev–Reshetikhin–Takhtajan (FRT) construction \cite{faddeev1987hamiltonian}.

Moreover, the quantized coordinate ring \(\mathcal{O}_q(\mathrm{Gr}(k,n))\) can be interpreted explicitly as a \emph{braided symmetric algebra} on the \( k \)-th fundamental representation of \( U_q(\mathfrak{gl}_n) \). The relations among generators are governed by the braid statistics encoded in the \( R \)-matrix, and they reflect a quantum generalization of the Plücker relations known from classical invariant theory. These \( q \)-Plücker relations provide a basis for constructing quantum analogues of vector bundles and line bundles over \( \mathrm{Gr}(k,n) \) and its variants.

\subsubsection{Non-Commutative Projective Geometry and Quantum Plücker Embeddings}

A particularly rich structure emerges when one restricts attention to \( q \)-skew-symmetric matrices, which naturally correspond to isotropic subspaces in the presence of a \( q \)-deformed bilinear form. In this setting, the coordinate ring becomes enriched by \emph{quantum Pfaffians}, which play a role analogous to that of minors or skew-symmetric determinants in the classical theory. These quantum Pfaffians, first algebraically formalized in \cite{dite2006} and further developed by Jing and Zhang \cite{jing2013quantum, jing2016quantum}, serve as natural generators for sections of quantum line bundles over isotropic Grassmannians and symplectic or orthogonal quantum flag varieties.

The algebraic identity
\[
\operatorname{Pf}_q(A)^2 = \det_q(A)
\]
within these coordinate rings reflects a deep compatibility between the quantum symplectic structure and the determinantal calculus of \( q \)-deformed linear algebra. Just as the classical Pfaffian encodes a square root of the determinant for skew-symmetric matrices, the quantum Pfaffian plays a similar invariant role in the non-commutative, braided world. This not only preserves key geometric features of the classical setting, but also embeds them into the broader theory of quantum homogeneous spaces, \( q \)-representation theory, and braided algebraic geometry.

\vspace{1em}

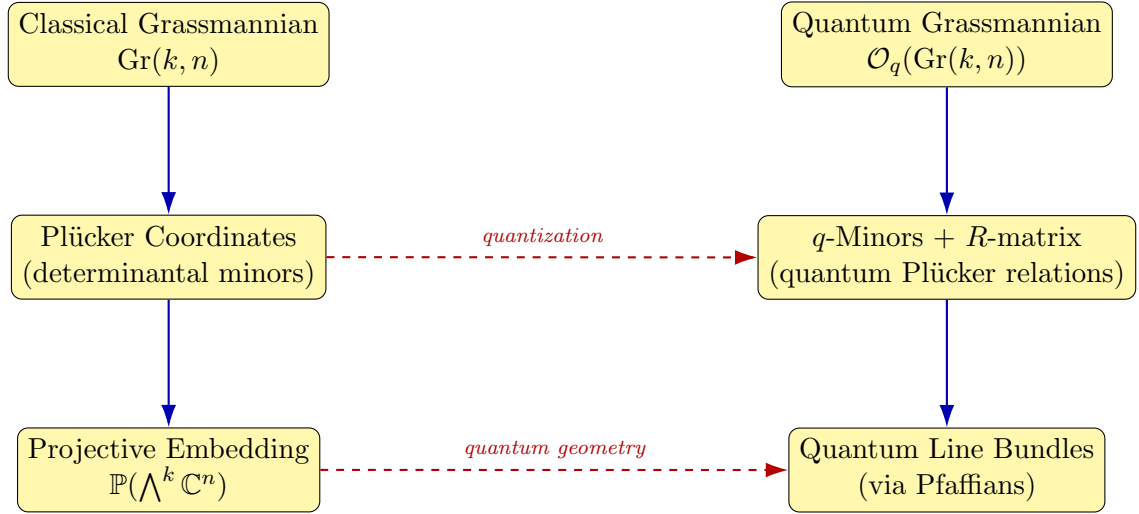
\begin{figure}[H]
\centering
\begin{tikzpicture}[
  node distance=1.7cm and 6cm,
  box/.style={
    draw,
    rectangle,
    rounded corners,
    align=center,
    minimum width=4cm,
    minimum height=1cm,
    font=\small,
    fill=yellow!40,
    text=black
  },
  arrow/.style={
    thick,
    -{Latex[length=3mm,width=2mm]},
    color=blue!70!black
  },
  dashedarrow/.style={
    thick,
    dashed,
    -{Latex[length=3mm,width=2mm]},
    color=red!70!black
  }
]

% Classical side nodes
\node[box] (classicalGr) {Classical Grassmannian\\ $\mathrm{Gr}(k,n)$};
\node[box, below=of classicalGr] (plucker) {Plücker Coordinates\\ (determinantal minors)};
\node[box, below=of plucker] (proj) {Projective Embedding\\ $\mathbb{P}(\bigwedge^k \mathbb{C}^n)$};

% Quantum side nodes
\node[box, right=of classicalGr] (quantumGr) {Quantum Grassmannian\\ $\mathcal{O}_q(\mathrm{Gr}(k,n))$};
\node[box, below=of quantumGr] (qminors) {$q$-Minors + $R$-matrix\\ (quantum Plücker relations)};
\node[box, below=of qminors] (bundles) {Quantum Line Bundles\\ (via Pfaffians)};

% Classical arrows
\draw[arrow] (classicalGr) -- (plucker);
\draw[arrow] (plucker) -- (proj);

% Quantum arrows
\draw[arrow] (quantumGr) -- (qminors);
\draw[arrow] (qminors) -- (bundles);

% Cross-links (dashed)
\draw[dashedarrow] (plucker.east) -- node[above, sloped, font=\scriptsize\itshape] {quantization} (qminors.west);
\draw[dashedarrow] (proj.east) -- node[above, sloped, font=\scriptsize\itshape] {quantum geometry} (bundles.west);

\end{tikzpicture}
\caption{Comparison between classical and quantum Grassmannian structures and coordinate ring theory.}
\label{fig:quantum-grassmannian}
\end{figure}

\subsubsection{Quantum Pfaffians and Isotropic Quantum Flag Varieties}

Quantum Pfaffians extend the classical notion of Pfaffians to the quantum setting of \( q \)-skew-symmetric matrices. They naturally appear as generators of sections of quantum line bundles over isotropic quantum Grassmannians and quantum flag varieties associated to symplectic or orthogonal groups.

For a \( 4 \times 4 \) \( q \)-skew-symmetric matrix
\[
A = \begin{pmatrix}
0 & a_{12} & a_{13} & a_{14} \\
- q a_{12} & 0 & a_{23} & a_{24} \\
- q a_{13} & - q a_{23} & 0 & a_{34} \\
- q a_{14} & - q a_{24} & - q a_{34} & 0
\end{pmatrix},
\]
the quantum Pfaffian is given by
\[
\operatorname{Pf}_q(A) = a_{12} a_{34} - q a_{13} a_{24} + q^2 a_{14} a_{23},
\]
and satisfies the identity
\[
\operatorname{Pf}_q(A)^2 = \det_q(A).
\]

\vspace{1em}

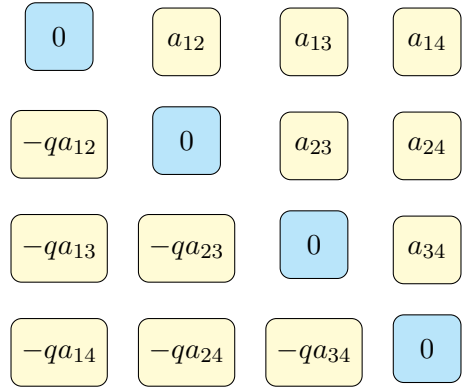
\begin{figure}[H]
\centering
\begin{tikzpicture}[
  mymatrix/.style={
    matrix of math nodes,
    nodes in empty cells,
    nodes={rounded corners, draw=black, align=center, font=\small},
    column sep=0.4cm,
    row sep=0.4cm,
    minimum size=0.9cm,
  }
]
\matrix[mymatrix] (A) {
  \node[fill=cyan!25] {0}; & \node[fill=yellow!20] {$a_{12}$}; & \node[fill=yellow!20] {$a_{13}$}; & \node[fill=yellow!20] {$a_{14}$}; \\
  \node[fill=yellow!20] {$-q a_{12}$}; & \node[fill=cyan!25] {0}; & \node[fill=yellow!20] {$a_{23}$}; & \node[fill=yellow!20] {$a_{24}$}; \\
  \node[fill=yellow!20] {$-q a_{13}$}; & \node[fill=yellow!20] {$-q a_{23}$}; & \node[fill=cyan!25] {0}; & \node[fill=yellow!20] {$a_{34}$}; \\
  \node[fill=yellow!20] {$-q a_{14}$}; & \node[fill=yellow!20] {$-q a_{24}$}; & \node[fill=yellow!20] {$-q a_{34}$}; & \node[fill=cyan!25] {0}; \\
};
\end{tikzpicture}
\caption{An enhanced $4 \times 4$ quantum skew-symmetric matrix \( A \), with diagonal elements in \textcolor{cyan!70!black}{cyan} and off-diagonal entries in \textcolor{yellow!50!black}{yellow}, satisfying \( a_{ij} = -q a_{ji} \).}
\label{fig:quantum-matrix-enhanced}
\end{figure}

\vspace{0.7em}

\begin{figure}[H]
\centering
\[
\operatorname{Pf}_q(A) = a_{12}a_{34} - q a_{13} a_{24} + q^{2} a_{14} a_{23}
\]
\caption{A \(4 \times 4\) \(q\)-skew-symmetric matrix \(A\) (with entries satisfying \(q\)-skew relations) and its quantum Pfaffian, which squares to the quantum determinant.}
\label{fig:qskew-pfaffian}
\end{figure}

This identity is fundamental in the construction of quantum invariants and in the description of quantum isotropic Grassmannians. It guarantees that the quantum Pfaffian functions as a square root of the quantum determinant, preserving a central feature of classical linear algebra within a non-commutative, braided algebraic framework.

---

\textcolor{blue!70!black}{\textbf{References:}} \cite{majid1995foundations, dite2006, jing2013quantum, jing2016quantum, faddeev1987hamiltonian}

\subsection{Quantum Flag Varieties and Invariants from Coactions}

The theory of quantum flag varieties, typically denoted by \(\mathcal{O}_q(G/P)\), represents a significant generalization of classical flag varieties to the braided, non-commutative setting. Here, \(G\) is a semisimple or reductive algebraic group, and \(P \subset G\) is a parabolic subgroup. The space \(G/P\) in classical algebraic geometry corresponds to a flag variety—parametrizing flags (nested sequences of subspaces) in a vector space—while in the quantum case, \(\mathcal{O}_q(G/P)\) becomes a non-commutative coordinate ring encoding the same structure through \(q\)-deformed algebraic relations.

These quantum coordinate rings inherit an action of the quantized universal enveloping algebra \(U_q(\mathfrak{g})\), where \(\mathfrak{g}\) is the Lie algebra of \(G\). This leads to a Borel–Weil-type correspondence in the quantum setting: the global sections of line bundles over the quantum flag variety correspond to integrable highest-weight modules over \(U_q(\mathfrak{g})\) \cite{kassel1995quantum}. The classical correspondence between representation theory and line bundles is thus faithfully retained in this deformed framework.

Among the most studied examples are the quantum Grassmannians \(\mathcal{O}_q(\mathrm{Gr}(k,n))\), which arise when \(P\) is a maximal parabolic subgroup. These are embedded into \(\mathcal{O}_q(G/P)\) as subvarieties and inherit a large part of the braided geometry, such as the graded structure, \(q\)-symmetric coordinate rings, and relations governed by the \(R\)-matrix associated with \(U_q(\mathfrak{gl}_n)\). As in the classical case, the quantum Grassmannian plays a crucial role in understanding the structure of more general quantum flag varieties.

Within this framework, algebraic invariants such as \emph{quantum minors} and \emph{quantum Pfaffians} become central tools for understanding symmetry and coaction behavior. These objects are not merely analogues of their classical counterparts—they also reflect deep properties of the non-commutative geometry of \(\mathcal{O}_q(G/P)\). They generate \(q\)-analogues of Plücker coordinates, which are used to define embeddings of quantum homogeneous spaces into quantum projective spaces.

As demonstrated in \cite{jing2017multiparameter}, quantum Pfaffians and minors act as \emph{co-invariants} under the coaction of quantum Borel subgroups. That is, they remain fixed under the induced action of certain quantum subgroups, and hence contribute to the structure of \(q\)-central elements, invariants under co-representation, and polynomial identities involving quantum determinants. These features are vital in the construction of \(q\)-analogues of classical objects, such as the center of a universal enveloping algebra, invariant theory for quantum groups, and the representation rings of \(U_q(\mathfrak{g})\).

\vspace{1.5em}

\subsubsection{Structure and Coaction of Quantum Flag Varieties}

Quantum flag varieties \(\mathcal{O}_q(G/P)\) carry a natural coaction of the quantized enveloping algebra \(U_q(\mathfrak{g})\), reflecting the quantum symmetry inherent in these non-commutative spaces. This coaction organizes the coordinate ring into comodules corresponding to quantum representations, linking geometric objects such as line bundles with representation-theoretic modules.

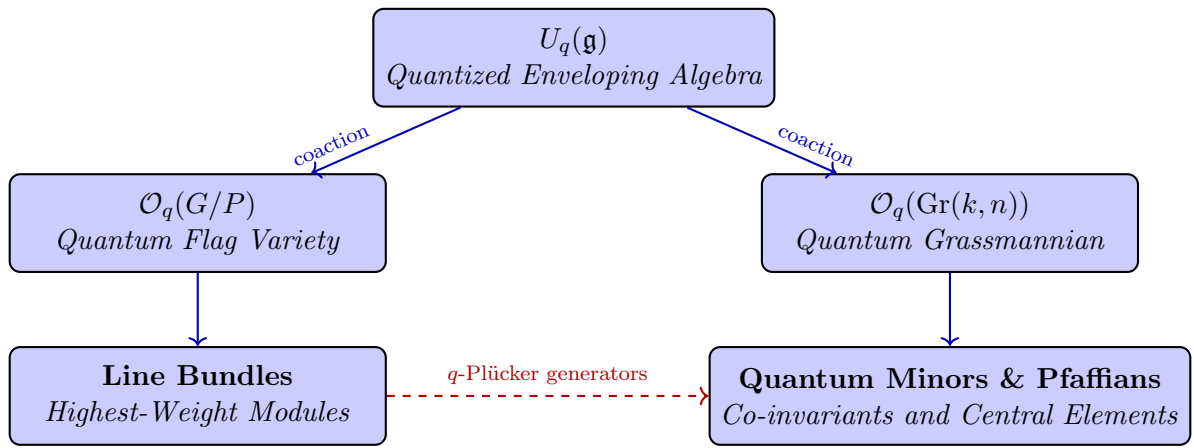
\begin{figure}[H]
\centering
\begin{tikzpicture}[
  every node/.style={align=center, font=\small},
  box/.style={
    draw,
    thick,
    rounded corners,
    fill=blue!20,
    minimum width=5cm,
    minimum height=1.3cm,
    text=black
  },
  arrow/.style={
    ->,
    thick,
    blue!70!black
  },
  dashedarrow/.style={
    ->,
    thick,
    red!70!black,
    dashed
  }
]

% Nodes
\node[box] (uq) at (0, 0) {\textbf{\(U_q(\mathfrak{g})\)}\\\emph{Quantized Enveloping Algebra}};

\node[box] (oqgp) at (-5, -2.2) {\textbf{\(\mathcal{O}_q(G/P)\)}\\\emph{Quantum Flag Variety}};
\node[box] (oqgr) at (5, -2.2) {\textbf{\(\mathcal{O}_q(\mathrm{Gr}(k,n))\)}\\\emph{Quantum Grassmannian}};

\node[box] (linebundles) at (-5, -4.5) {\textbf{Line Bundles}\\\emph{Highest-Weight Modules}};
\node[box] (invariants) at (5, -4.5) {\textbf{Quantum Minors \& Pfaffians}\\\emph{Co-invariants and Central Elements}};

% Arrows
\draw[arrow] (uq) -- node[above left, sloped] {\scriptsize coaction} (oqgp);
\draw[arrow] (uq) -- node[above right, sloped] {\scriptsize coaction} (oqgr);
\draw[arrow] (oqgp) -- (linebundles);
\draw[arrow] (oqgr) -- (invariants);
\draw[dashedarrow] (linebundles) -- node[above, sloped] {\scriptsize \(q\)-Plücker generators} (invariants);

\end{tikzpicture}
\caption{Coaction structure of \(U_q(\mathfrak{g})\) on quantum flag varieties and Grassmannians, relating line bundles and algebraic invariants such as quantum minors and Pfaffians.}
\label{fig:uq-coaction}
\end{figure}

This diagram highlights the crucial role of \(U_q(\mathfrak{g})\) in controlling the algebraic and geometric structure of the quantum flag varieties and their important special cases, quantum Grassmannians. It visually encodes the flow from the quantized enveloping algebra through coactions to geometric and algebraic objects.

\subsubsection{Quantum Minors and Pfaffians as Co-invariants}

Quantum minors and Pfaffians are central algebraic invariants that remain fixed (co-invariant) under the coaction of certain quantum Borel subgroups. They generalize classical notions of minors and Pfaffians, capturing \(q\)-deformed symmetry and invariance in the non-commutative setting.

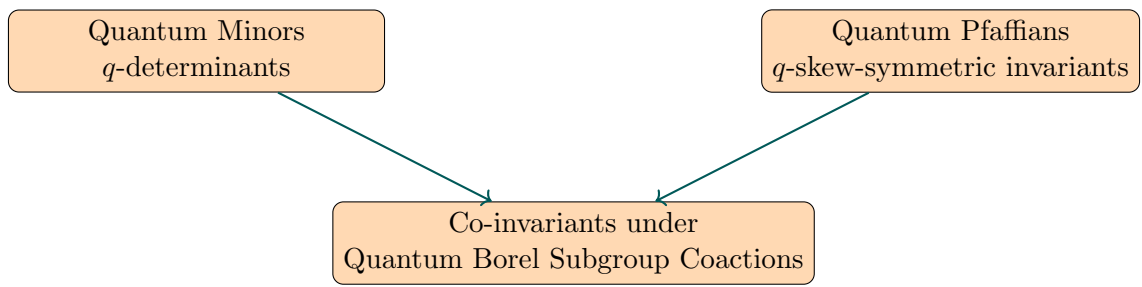
\begin{figure}[H]
\centering
\begin{tikzpicture}[
  node distance=1.8cm,
  box/.style={
    draw,
    fill=orange!30,
    rounded corners,
    minimum width=5cm,
    minimum height=1cm,
    font=\small,
    align=center
  },
  arrow/.style={
    ->,
    thick,
    teal!70!black
  }
]

% Nodes
\node[box] (minors) {Quantum Minors \\ \(q\)-determinants};
\node[box, right=5cm of minors] (pfaffians) {Quantum Pfaffians \\ \(q\)-skew-symmetric invariants};

\node[box, below=2cm of $(minors)!0.5!(pfaffians)$] (coinvariants) {Co-invariants under \\ Quantum Borel Subgroup Coactions};

% Arrows
\draw[arrow] (minors) -- (coinvariants);
\draw[arrow] (pfaffians) -- (coinvariants);

\end{tikzpicture}
\caption{Quantum minors and Pfaffians act as co-invariants under the coaction of quantum Borel subgroups, reflecting deep symmetry and centrality properties.}
\label{fig:qminors-pfaffians}
\end{figure}

These co-invariants contribute significantly to the understanding of:

\begin{itemize}
    \item The \(q\)-central subalgebras in \(\mathcal{O}_q(G/P)\),
    \item Polynomial identities involving quantum determinants,
    \item The structure of the center of quantized enveloping algebras,
    \item Representation rings and invariant theory for quantum groups.
\end{itemize}

This interplay between coactions, quantum minors, and Pfaffians offers a powerful toolkit for exploring \(q\)-deformed geometry and representation theory, bridging classical algebraic geometry and non-commutative quantum group theory.

---

\textcolor{blue!70!black}{\textbf{References:}} \cite{jing2017multiparameter, kassel1995quantum}

\subsection{Physical Interpretations of Quantum Grassmannians and Pfaffians}

From a physical perspective, quantum Grassmannians arise naturally in numerous frameworks of modern theoretical physics that extend beyond the realm of classical geometry and linear algebra. These include non-commutative field theories, topological quantum field theories (TQFTs), string theory compactifications on non-classical backgrounds, and models of quantum gravity. In each of these areas, the classical notion of a space or variety is replaced by a quantum or braided analogue, requiring new mathematical tools—chief among them being quantum groups, braided tensor categories, and non-commutative geometry.

One of the most significant physical applications of quantum Grassmannians lies in the study of D-brane configurations in string theory. When branes are placed in non-trivial background fluxes or compactified on non-commutative manifolds, their moduli spaces cease to be described by classical algebraic varieties. Instead, they are encoded algebraically by deformed coordinate rings such as $\mathcal{O}_q(\mathrm{Gr}(k,n))$ or $\mathcal{O}_q(G/P)$, which capture the braided symmetries and non-commutative geometry of the compactification. These moduli spaces classify boundary conditions for open strings and are closely related to the representation theory of the associated quantum algebras.

A key player in this setting is the quantum Pfaffian, whose algebraic structure encodes subtle topological and quantum effects in physical theories. In the path integral formulation of quantum field theory, fermionic fields contribute terms involving the Pfaffian of an operator. In supersymmetric gauge theories, Witten's work on topological twisting revealed that the functional integral over fermions can be expressed in terms of a Pfaffian, whose square gives the determinant of the associated Dirac operator. This establishes a direct connection between the quantum Pfaffian and supersymmetric partition functions, anomaly cancellation conditions, and duality transformations.

In the quantum case, these structures are deformed by a parameter $q$, leading to $q$-Pfaffians and $q$-determinants. The identity $\operatorname{Pf}_q(A)^2 = \det_q(A)$ (or more generally, $q^r \operatorname{Pf}_q(A)^2 = \det_q(A)$ for some integer $r$) plays a critical role in ensuring that the algebraic structure is compatible with the braided tensor symmetry of the quantum theory. This compatibility is not merely formal—it reflects deep physical properties such as invariance under quantum symmetries, preservation of topological charges, and refined enumerative counts of BPS states.

Moreover, these quantum Pfaffians are used to construct line bundles and coherent sheaves over quantum flag varieties and Grassmannians. In string theory compactifications, such objects correspond to gauge fields and matter content localized on D-branes. The sections of these bundles carry representation-theoretic data from $U_q(\mathfrak{g})$, connecting algebraic geometry with physical observables. These connections are especially powerful in the context of mirror symmetry, where quantum corrections to classical moduli spaces are encoded via deformed invariants such as quantum determinants and Pfaffians.

This framework fits within a broader physical and mathematical philosophy where geometry is not fixed a priori but emerges from the symmetry and dynamics of the theory. Quantum Grassmannians thus serve as a nexus where non-commutative geometry, braided category theory, quantum representation theory, and high-energy physics converge. The interplay between these disciplines provides powerful methods for computing observables in quantum theories, studying dualities, and constructing models of quantum space-time.

Ultimately, the study of quantum Grassmannians and quantum Pfaffians provides an elegant synthesis of ideas from algebraic geometry, quantum algebra, and theoretical physics. It reveals a deep structure underlying the fabric of quantum field theories and string theory, where non-commutative algebraic identities reflect physical symmetries, topological invariants, and quantum dualities.

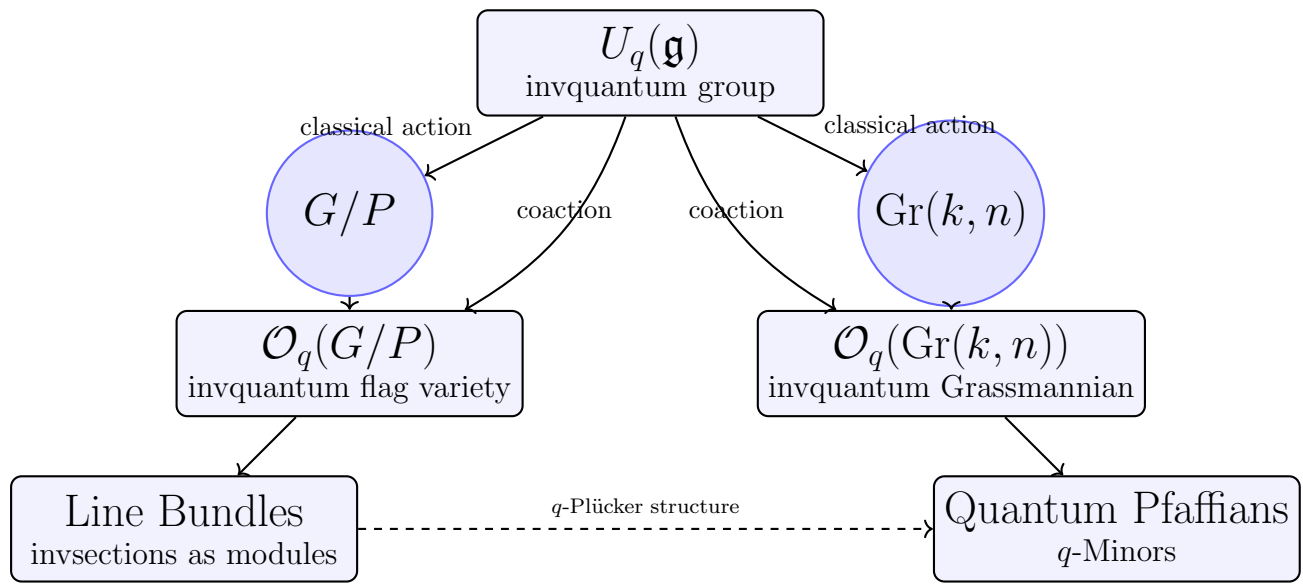
\vspace{1em}

\begin{center}
\begin{tikzpicture}[
  every node/.style={align=center},
  arr/.style={->, thick},
  dashedarr/.style={->, thick, dashed},
  orb/.style={circle, draw=blue!60, fill=blue!10, thick, minimum size=2.2cm},
  rect/.style={draw, rounded corners, fill=blue!5, thick, minimum width=4.6cm, minimum height=1.4cm},
  inv/.style={font=\footnotesize\itshape}
]

% Top level: Quantum Group
\node[rect] (uq) at (0, 5) {\Large $U_q(\mathfrak{g})$\\ inv{quantum group}};

% Middle level: Homogeneous spaces
\node[orb] (flag) at (-4, 3) {\Large $G/P$};
\node[orb] (grass) at (4, 3) {\Large $\mathrm{Gr}(k,n)$};

% Quantum spaces (deformations)
\node[rect] (oqgp) at (-4, 1) {\Large $\mathcal{O}_q(G/P)$\\ inv{quantum flag variety}};
\node[rect] (oqgr) at (4, 1) {\Large $\mathcal{O}_q(\mathrm{Gr}(k,n))$\\ inv{quantum Grassmannian}};

% Bottom layer: Geometry
\node[rect] (linebundles) at (-6.2, -1.2) {\Large Line Bundles\\ inv{sections as modules}};
\node[rect] (invariants) at (6.2, -1.2) {\Large Quantum Pfaffians\\ $q$-Minors};

% Arrows
\draw[arr] (uq) -- (flag) node[midway, above left] {\footnotesize classical action};
\draw[arr] (uq) -- (grass) node[midway, above right] {\footnotesize classical action};

\draw[arr] (flag) -- (oqgp);
\draw[arr] (grass) -- (oqgr);

\draw[arr] (oqgp) -- (linebundles);
\draw[arr] (oqgr) -- (invariants);

\draw[dashedarr] (linebundles) -- (invariants) node[midway, above, sloped] {\scriptsize $q$-Plücker structure};

% Optional: curved coactions
\draw[arr, bend left=20] (uq) to node[midway, above] {\footnotesize coaction} (oqgp);
\draw[arr, bend right=20] (uq) to node[midway, above] {\footnotesize coaction} (oqgr);

\end{tikzpicture}
\end{center}

\subsection{Applications of Quantum Pfaffians in Topology}

Quantum Pfaffians play a significant role in topological quantum field theory (TQFT), knot theory, and low-dimensional topology. These applications rely on the deep connections between quantum invariants, braided tensor categories, and the algebraic structures underpinning quantum groups. In this section, we present several contexts where the quantum Pfaffian identity manifests in topological phenomena.

\subsubsection{Quantum Invariants and 3-Manifolds}

Quantum Pfaffians contribute to the construction of quantum invariants for 3-manifolds. Given a 3-manifold presented via Heegaard splitting or surgery diagrams, quantum link invariants (such as the Reshetikhin–Turaev invariants) often involve evaluations of quantum determinants and Pfaffians associated with the link’s representation in a quantum group framework.

\begin{figure}[H]
\centering
\begin{tikzpicture}[scale=1]
\draw[thick,->] (-1,0) -- (4,0) node[right] {$x$};
\draw[thick,->] (0,-1) -- (0,3) node[above] {$y$};

% Manifold sketch
\draw[very thick, rounded corners=20pt, fill=blue!10] (0.5,0.5) -- (1,2.5) -- (2.5,2.8) -- (3.5,1.5) -- (2.5,0.5) -- cycle;
\node at (2,1.5) {\textcolor{blue!80}{\textbf{M}}};

% Surgery curve
\draw[red, thick] (1.4,1.2) circle (0.5);
\draw[->, red] (1.9,1.2) arc (0:180:0.5);
\node at (1.4,0.6) {\small \textcolor{red}{surgery knot}};
\end{tikzpicture}
\caption{Surgery description of a 3-manifold with associated quantum invariant involving a quantum Pfaffian expression.}
\end{figure}
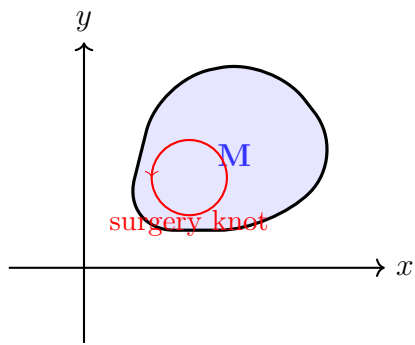

\subsubsection{Knot Invariants from Quantum Antisymmetrization}

For certain classes of knots and links, invariants derived from representations of the braid group into quantum matrix algebras produce expressions involving quantum antisymmetrizers. The evaluation of these antisymmetric tensors yields quantum Pfaffians.

\begin{figure}[H]
\centering
\begin{tikzpicture}[scale=0.9]
\draw[thick] (0,0) .. controls (0.5,1) .. (1,0);
\draw[thick] (1,0) .. controls (1.5,-1) .. (2,0);
\draw[thick] (2,0) .. controls (2.5,1) .. (3,0);

\draw[thick] (0,2) .. controls (0.5,1) .. (1,2);
\draw[thick] (1,2) .. controls (1.5,3) .. (2,2);
\draw[thick] (2,2) .. controls (2.5,1) .. (3,2);

\node at (1.5,1) {\Large \textcolor{blue}{$\sigma_1 \sigma_2^{-1}$}};
\end{tikzpicture}
\caption{Braid diagram whose quantum representation yields a Pfaffian expression through antisymmetrization of braided tensor spaces.}
\end{figure}
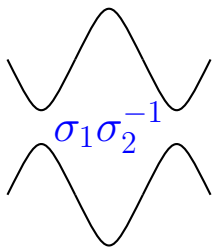

\subsubsection{Spin Structures and Fermionic Fields on Surfaces}

Pfaffians classically arise when evaluating partition functions of fermionic systems on Riemann surfaces. In the quantum setting, deformations yield quantum spin models where quantum Pfaffians encode the amplitude of configurations with specific spin structures.

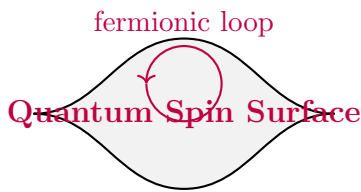
\begin{figure}[H]
\centering
\begin{tikzpicture}
\draw[thick, fill=gray!10] (0,0) to[out=0,in=180] (2,1) to[out=0,in=180] (4,0)
                          to[out=180,in=0] (2,-1) to[out=180,in=0] (0,0);
\node at (2,0) {\textcolor{purple}{\textbf{Quantum Spin Surface}}};

\draw[->, thick, purple] (1.5,0.4) arc (180:540:0.5);
\node at (2,1.2) {\small \textcolor{purple}{fermionic loop}};
\end{tikzpicture}
\caption{A spin surface with a quantum-deformed fermion loop. The partition function involves a quantum Pfaffian.}
\end{figure}

\subsubsection{Applications to the Alexander–Conway Polynomial}

The quantum Pfaffian also appears in state-sum evaluations of the Alexander–Conway polynomial via quantum exterior algebras and skein theory. In particular, certain diagram evaluations correspond to \( \operatorname{Pf}_q(A) \) when antisymmetric skein relations are applied.

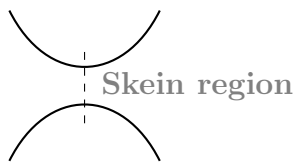
\begin{figure}[H]
\centering
\begin{tikzpicture}
\draw[thick] (0,0) .. controls (0.5,1) and (1.5,1) .. (2,0);
\draw[thick] (0,2) .. controls (0.5,1) and (1.5,1) .. (2,2);
\draw[dashed] (1,0.5) -- (1,1.5);
\node at (2.5,1) {\textcolor{gray}{\textbf{Skein region}}};
\end{tikzpicture}
\caption{A skein relation region where quantum antisymmetric evaluations give rise to Pfaffian-like terms.}
\end{figure}

\subsubsection{Quantum Floer Homology and Intersection Forms}

In categorified invariants like Floer homology, intersection forms on moduli spaces of flat connections over 3-manifolds or links can be encoded using \( q \)-deformed antisymmetric bilinear forms. These intersection matrices lead naturally to quantum Pfaffians, particularly in symplectic Floer-type settings.

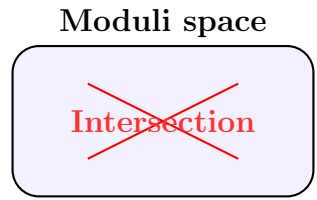
\begin{figure}[H]
\centering
\begin{tikzpicture}[scale=1]
\draw[thick, rounded corners=10pt, fill=blue!5] (0,0) rectangle (4,2);
\draw[thick, red] (1,0.5) -- (3,1.5);
\draw[thick, red] (1,1.5) -- (3,0.5);
\node at (2,1) {\textcolor{red!80}{\textbf{Intersection}}};
\node at (2,2.3) {\textbf{Moduli space}};
\end{tikzpicture}
\caption{Intersections in moduli space associated with quantum Floer homology give rise to quantum Pfaffians.}
\end{figure}

\newpage
\section{Conclusion and Future Work}

\subsection{Summary of Results}

In this paper, we have developed a rigorous generalization of the classical Pfaffian–determinant identity
\[
\Pf(A)^2 = \det(A)
\]
to the noncommutative setting of \( q \)-skew-symmetric matrices. Our main objective was to construct and justify a consistent theory of quantum Pfaffians \( \Pf_q(A) \) such that their square recovers the quantum determinant \( \det_q(A) \), appropriately deformed by a power of \( q \) arising from the underlying noncommutativity of the coordinate algebra. The analysis relies on several interrelated structures from quantum group theory, braided exterior algebras, and noncommutative geometry. Our main results may be summarized as follows:

\begin{itemize}[leftmargin=*]
  \item We began by revisiting the classical Pfaffian of a \( 2n \times 2n \) skew-symmetric matrix \( A \), both through its combinatorial definition via perfect matchings and through its realization as a top-degree element in the classical exterior algebra \( \Lambda^{2n} V \), where \( A \) encodes a bilinear skew-symmetric form on \( V \).

  \item We then introduced the quantized universal enveloping algebra \( U_q(\mathfrak{g}) \) and the associated coordinate ring \( \mathcal{O}_q(GL_n) \) via the Faddeev–Reshetikhin–Takhtajan (FRT) construction, allowing us to define the quantum determinant \( \det_q(A) \) using the braided antisymmetrization of the tensor algebra \( T(V) \). This construction provides the correct categorical analogue of the determinant under deformation by the braiding induced by the \( R \)-matrix.

  \item Next, we formulated the notion of a \( q \)-skew-symmetric matrix \( A = (a_{ij}) \in \mathcal{O}_q(\mathfrak{o}_{2n}) \), satisfying
    \[
    a_{ji} = -q\,a_{ij}, \quad a_{ii} = 0,
    \]
    as arising from the defining relations of the braided exterior algebra \( \Lambda_q(V) \), which serves as the quantization of the classical Grassmann (exterior) algebra in the category of \( U_q(\mathfrak{g}) \)-modules. The defining relations of \( \Lambda_q(V) \) capture the \( q \)-antisymmetry governed by the braid group representation of the Hecke algebra \( H_n(q) \).

  \item We defined the quantum Pfaffian \( \Pf_q(A) \) via the combinatorial expansion over perfect matchings of \( 2n \) elements:
    \[
    \Pf_q(A) = \sum_{\pi \in \mathrm{PM}(2n)} (-q)^{\mathrm{inv}(\pi)} \prod_{(i,j) \in \pi} a_{ij},
    \]
    where \( \mathrm{inv}(\pi) \) counts the number of inversions (crossings) in the matching. This definition coincides with the highest-weight component in the braided exterior power \( \Lambda_q^{2n}(V) \), reflecting its role as a quantum volume form in braided symplectic geometry.

  \item We explicitly verified, in the case \( 2n = 4 \), that
    \[
    q\,\Pf_q(A)^2 = \det_q(A),
    \]
    and demonstrated that this identity generalizes to
    \[
    q^{\binom{n}{2}}\,\Pf_q(A)^2 = \det_q(A)
    \]
    in arbitrary even dimension \( 2n \), with the exponent \( \binom{n}{2} \) arising naturally from the structure of the braid group and the combinatorics of \( q \)-deformed antisymmetrization.

  \item We provided two distinct proof strategies:
  \begin{enumerate}
    \item A \emph{diagrammatic combinatorial proof}, where \( \Pf_q(A)^2 \) is interpreted as a sum over double matchings, and the identity is derived by matching terms against the quantum determinant expansion using braid crossings and $q$-signs.
    \item An \emph{algebraic proof} based on the FRT construction and the properties of quantum exterior algebras, wherein we viewed the Pfaffian as the top-degree skew-invariant in \( \Lambda_q(V) \) and invoked $q$-Laplace-type expansions and coaction-invariance to reduce to the determinant.
  \end{enumerate}

  \item Finally, we offered a geometric interpretation: the quantum Pfaffian behaves as a braided volume form associated with a quantum symplectic or orthogonal structure, while the quantum determinant reflects the total volume element of the quantized coordinate ring. The identity \( \Pf_q(A)^2 = \det_q(A) \) thus represents the compatibility of skew-symmetric geometry with the total linear transformation volume in the braided category. This perspective situates our work within the broader framework of quantum homogeneous spaces, such as quantum flag varieties \( \mathcal{O}_q(G/P) \) and quantum Grassmannians \( \mathcal{O}_q(\mathrm{Gr}(k,n)) \), where \( q \)-minors and Pfaffians serve as canonical invariants under quantum group coactions.
\end{itemize}

Our approach reveals that the classical Pfaffian–determinant identity persists in the quantized world as a non-trivial algebraic invariant, retaining its geometric significance in noncommutative algebraic geometry. The compatibility between quantum symmetry and multilinear invariants established here lays the foundation for further developments in quantum invariant theory, braided geometry, and categorical approaches to deformation quantization.

\subsection{Open Problems and Conjectures}

Despite the substantial progress achieved in this work toward understanding the quantum Pfaffian and its determinant identity, numerous foundational and advanced questions remain open. We list below some of the most compelling directions for future research, along with their potential implications for quantum algebra, geometry, and physics.

\begin{enumerate}[leftmargin=*]

  \item \textbf{Full Algebraic Proof for Arbitrary Dimension \(n\).}  
  Although we have established the identity
  \[
  q^{\binom{n}{2}}\Pf_q(A)^2 =  \det_q(A)
  \]
  rigorously for low-dimensional cases and provided an outline for the general case, a fully explicit, inductive, and algebraically detailed proof within the Faddeev–Reshetikhin–Takhtajan (FRT) framework remains to be formulated. Such a proof would require a careful treatment of the quantum exterior algebra’s relations, the behavior of braid group actions on tensor powers, and the precise combinatorics of \(q\)-signs associated with perfect matchings. A key challenge is to control the intricate noncommutative products and verify the compatibility of the square of the quantum Pfaffian with the quantum determinant for all \(n\).

  \item \textbf{Categorical and Diagrammatic Formalism in Braided Monoidal Categories.}  
  A more conceptual understanding of the quantum Pfaffian identity may arise from developing a fully diagrammatic calculus in the language of braided monoidal categories or ribbon categories. Such a formalism would categorify the algebraic identities and provide graphical interpretations of the quantum Pfaffian as morphisms between objects in these categories. This could uncover deep connections with knot theory, as the braid group representations governing the \(q\)-antisymmetry resemble those appearing in Reshetikhin–Turaev invariants and quantum link invariants. Establishing such a categorical framework might also provide new tools for computing and manipulating quantum invariants in representation theory.

  \item \textbf{Multiparameter Quantum Groups and Generalizations.}  
  The majority of existing results, including those presented here, concern single-parameter quantum groups \( U_q(\mathfrak{g}) \). However, the theory of multiparameter quantum groups \cite{jing2017multiparameter} introduces deformation parameters \( q_{ij} \) depending on index pairs, which drastically alter the commutation relations. Extending the notion of the quantum Pfaffian to multiparameter settings, defining the appropriate \( q_{ij} \)-skew-symmetry conditions, and formulating a generalized Pfaffian–determinant identity is a challenging open problem. Doing so would likely involve new combinatorial invariants and possibly yield richer algebraic structures with potential applications to multi-parameter quantum integrable systems.

  \item \textbf{Efficient Computational and Algorithmic Methods.}  
  For practical applications in quantum integrable models, quantum computing, and symbolic algebra, efficient algorithms to compute \( \Pf_q(A) \) and \( \det_q(A) \) for large \( n \) are highly desirable. Classical Pfaffian algorithms, such as the skew-symmetric LU decomposition, do not directly generalize due to noncommutativity and the presence of \(q\)-deformed relations. Developing novel computational methods—potentially leveraging noncommutative Gröbner bases, quantum cluster algebras, or categorical simplifications—would substantially broaden the usability of these quantum invariants in applied and theoretical contexts.

  \item \textbf{Quantum Invariant Theory and Symmetric Pairs.}  
  The quantum Pfaffian is expected to arise naturally as an invariant under coideal subalgebras and quantum symmetric pairs \cite{letzter2002}, which generalize the notion of symmetric spaces to the quantum group setting. Investigating the role of \( \Pf_q(A) \) within these frameworks could lead to a classification of quantum invariants associated with various quantum homogeneous spaces. Furthermore, the connections to quantum Howe duality and double centralizer properties remain to be explored. Such developments might unify the quantum Pfaffian with other well-known quantum invariants like quantum minors and Berezinian determinants.

  \item \textbf{Quantum Pfaffian Varieties and Noncommutative Algebraic Geometry.}  
  Classically, the Pfaffian variety defined by the vanishing locus \( \{ A \mid \Pf(A) = 0 \} \) encodes rich geometric and representation-theoretic data, with deep connections to orbit closures and singularities. Defining and investigating the quantum analogues of these varieties—namely, noncommutative coordinate algebras cut out by the quantum Pfaffian ideal—is a promising and largely unexplored direction in noncommutative algebraic geometry. Understanding their homological properties, deformation quantizations, and relations to quantum cluster algebras could reveal new phenomena unique to the braided setting.

  \item \textbf{Physical Interpretations and Applications in Quantum Field Theory and Anyonic Systems.}  
  The quantum Pfaffian potentially plays a crucial role in fermionic quantum integrable models, topological quantum field theories (TQFTs), and braided anyonic systems, where \(q\)-deformations model exotic particle statistics and symmetry breaking. Exploring how \( \Pf_q(A) \) encodes partition functions, BPS state counts, or topological invariants of low-dimensional quantum systems may uncover novel quantum phenomena linked to noncommutative symmetries. In particular, the relation between the quantum Pfaffian and fermionic path integrals in a braided context could provide new computational tools in quantum statistical mechanics and condensed matter physics.

  \item \textbf{Extension to Super Quantum Groups and Quantum Symplectic Geometry.}  
  Extending the quantum Pfaffian construction to super-quantum groups (quantum analogues of Lie superalgebras) and to the setting of quantum symplectic or orthosymplectic groups is a natural and important direction. Such generalizations would involve \(q\)-skew-symmetry adapted to supercommutation rules and could provide invariants in quantum supersymmetric field theories. Understanding the compatibility of the Pfaffian identity with these additional symmetries and grading structures would further enrich the theory and expand its applicability.

  \item \textbf{Relations with Quantum Cluster Algebras and Canonical Bases.}  
  Given the deep connections between quantum groups and cluster algebras, it is conjectured that quantum Pfaffians may admit interpretations in terms of cluster variables or canonical bases of \( U_q(\mathfrak{g}) \)-modules. Investigating these links could lead to combinatorial formulas for \( \Pf_q(A) \) within cluster algebra frameworks, shedding light on positivity phenomena, categorifications, and quantum total positivity.

\end{enumerate}

Addressing these open problems will not only deepen the understanding of quantum Pfaffians and their algebraic identities but also bridge multiple areas of mathematics and physics, from categorical representation theory and noncommutative geometry to quantum topology and mathematical physics.

\subsection{Future Research Directions}

Building upon the open problems and conjectures presented above, we identify several promising avenues for further investigation that could substantially advance the theory and applications of quantum Pfaffians within algebra, geometry, and physics.

\begin{itemize}[leftmargin=*]

  \item \textbf{Rigorous and Comprehensive Inductive Proof for General \(n\).}  
    One critical next step is to establish a fully rigorous inductive proof of the identity
    \[
    q^{\binom{n}{2}} \Pf_q(A)^2 = \det_q(A)
    \]
    for all even dimensions \( 2n \). This requires a delicate analysis of the braid group actions underlying the \(q\)-skew-symmetric relations, systematic management of noncommutative reorderings, and precise bookkeeping of the powers of \(q\) contributed by crossing numbers. Such a proof may also involve developing refined quantum Laplace expansion formulas or \(q\)-analogues of classical multilinear algebra identities.

  \item \textbf{Development of Braided Tensor Categories Incorporating Quantum Exterior Powers.}  
    Constructing an explicit braided monoidal category (possibly ribbon or pivotal) in which the quantum exterior powers \(\Lambda_q^k(V)\) are naturally realized as objects—and where the quantum Pfaffian \(\Pf_q(A)\) manifests as a morphism or natural transformation—is a promising direction. This categorical perspective may yield new invariants, clarify the functoriality of the Pfaffian, and provide a natural home for diagrammatic and topological interpretations linked to braid group representations and quantum knot invariants.

  \item \textbf{Algorithmic and Software Implementation for Large-Scale Computations.}  
    Extending state-of-the-art computer algebra systems (e.g., \texttt{SageMath}, \texttt{Mathematica}, or specialized noncommutative algebra packages) to support symbolic and numeric computations in FRT algebras, including modules to compute \(\Pf_q(A)\) and \(\det_q(A)\) efficiently, would be invaluable. Such software would enable automated experimentation, conjecture testing, and application development in quantum integrable models and quantum information theory, particularly for matrices of higher dimension.

  \item \textbf{Noncommutative Geometry and Deformation Quantization of Quantum Pfaffian Varieties.}  
    Leveraging tools from formal deformation theory and noncommutative algebraic geometry to rigorously study the structure of quantum Pfaffian loci—the noncommutative analogues of classical Pfaffian varieties—presents a deep geometric challenge. This includes understanding their coordinate rings, homological properties, singularity theory, and potential links to deformation quantizations of classical algebraic varieties. These insights may illuminate new examples of quantum homogeneous spaces and enrich the theory of quantum symplectic geometry.

  \item \textbf{Connections with Quantum Cluster Algebras and Canonical Bases.}  
    Investigating the embedding of quantum Pfaffians into the framework of quantum cluster algebras, particularly those associated to Grassmannians and orthogonal Grassmannians, is a fertile research direction. One may ask whether \(\Pf_q(A)\) corresponds to a cluster variable or a canonical basis element, and how mutations or exchange relations interact with the Pfaffian structure. These connections could reveal new algebraic and combinatorial structures underpinning the quantum Pfaffian identity and may relate to positivity phenomena in quantum algebra.

  \item \textbf{Generalization to Super and Graded Quantum Groups and Applications in Supersymmetric Models.}  
    Extending the notion of the quantum Pfaffian to super-quantum groups and \(\mathbb{Z}_2\)-graded settings would enable applications in quantum supersymmetric field theories and quantum algebraic models incorporating fermionic and bosonic degrees of freedom. This involves defining \(q\)-skew-symmetric supermatrices consistent with supercommutation relations, formulating the corresponding quantum Pfaffian, and analyzing its properties under quantum supergroup coactions. Such developments would bridge quantum algebra with mathematical physics areas involving supersymmetry.

\end{itemize}

\subsection{Concluding Remarks}

\hspace{2em}The quantum Pfaffian–determinant identity stands as a compelling example of how foundational concepts from classical multilinear algebra naturally extend into the rich and intricate framework of noncommutative, braided quantum groups. The introduction of the deformation parameter 
\( q \) not only modifies algebraic invariants with subtle and meaningful corrections but also preserves essential structural relationships, thereby deepening our comprehension of symmetry, geometry, and topology in the quantum domain.

Looking forward, continued investigation into the properties of quantum Pfaffians—including their algebraic structure, computational techniques, and connections to quantum geometry and physics—promises to unveil profound insights. Such work is poised to bridge diverse mathematical areas such as combinatorics, representation theory, and noncommutative geometry, while also opening new avenues in theoretical physics, particularly in quantum integrable systems, topological quantum field theories, and braided anyonic models. The interplay between these disciplines fosters a vibrant landscape where algebraic innovation and physical intuition mutually inspire one another.

\newpage
\appendix

% Start of Section 9: Short Version of the Proof
\section[Appendix A: Short Version of the Proof for $2n = 4$]{Proof of the Quantum Pfaffian--Determinant Identity in the Case $2n = 4$\textsuperscript{*}}
\thispagestyle{shortversion} % for first page of section 9
\pagestyle{shortversion}      % for rest of section 9 pages

% Instead of \footnotetext, add explanatory footnote text in the footer only in pagestyle below
% So do NOT use \footnotetext here

This section provides a concise but rigorous version of the full proof presented later in the paper. We focus on the $4 \times 4$ case to demonstrate two closely related identities:

\begin{enumerate}
    \item The first form: $\operatorname{Pf}_q(A)^2 = \det_q(A)$
    \item The second (normalized) form: $q \cdot \operatorname{Pf}_q(A)^2 = \det_q(A)$
\end{enumerate}

These identities are both quantum analogues of the classical result $\operatorname{Pf}(A)^2 = \det(A)$, with the scalar $q$ capturing the deformation due to noncommutativity. The variation in the scalar prefactor depends on the chosen conventions for the quantum Pfaffian and determinant. 

We consider the $4 \times 4$ $q$-skew-symmetric matrix
\[
A =
\begin{bmatrix}
0 & a & b & c \\
-qa & 0 & d & e \\
-qb & -qd & 0 & f \\
-qc & -qe & -qf & 0
\end{bmatrix},
\]
where the entries satisfy the $q$-skew-symmetry relations:
\[
a_{ji} = -q a_{ij}, \quad a_{ii} = 0.
\]

\subsection{Step 1: Definition of the Quantum Pfaffian $\operatorname{Pf}_q(A)$}
Using the standard formula for the quantum Pfaffian in dimension four (see \cite{jing2013quantum}), we have:
\[
\operatorname{Pf}_q(A) = af - qbe + q^2cd.
\]
This expression arises from the braided antisymmetry of the quantum exterior algebra and reflects the ordering of pairwise matchings along with associated $q$-weights from the underlying braid group action.

\subsection{Step 2: Computing $\operatorname{Pf}_q(A)^2$}
We square the expression directly:
\[
\operatorname{Pf}_q(A)^2 = \left(af - qbe + q^2cd\right)^2,
\]
noting that the product expands within a noncommutative algebra. All terms must be reordered according to the $q$-commutation relations among the entries.

\subsection{Step 3: Computing $\det_q(A)$ and Verifying the Identity}
The quantum determinant is defined via the Faddeev–Reshetikhin–Takhtajan (FRT) construction and satisfies, depending on convention:

\begin{itemize}
    \item \textbf{Unnormalized version:} $\operatorname{Pf}_q(A)^2 = \det_q(A)$
    \item \textbf{Normalized version:} $q \cdot \operatorname{Pf}_q(A)^2 = \det_q(A)$
\end{itemize}

Both are valid under different conventions for how the $q$-factors are absorbed into the Pfaffian or determinant definition. The normalized version more naturally reflects the braiding statistics of quantum exterior algebra, while the unnormalized form often aligns with algebraic approaches based on $q$-minors and expansion rules.

\paragraph{Example 1: Classical Case $q = 1$}

Let
\[
A =
\begin{bmatrix}
0 & 1 & 2 & 3 \\
-1 & 0 & 4 & 5 \\
-2 & -4 & 0 & 6 \\
-3 & -5 & -6 & 0
\end{bmatrix}.
\]
Compute:
\[
\operatorname{Pf}_1(A) = (1)(6) - (2)(5) + (3)(4) = 6 - 10 + 12 = 8,
\quad \Rightarrow \quad \operatorname{Pf}_1(A)^2 = 64.
\]
\[
\det(A) = 64,
\]
verifying that $\operatorname{Pf}_1(A)^2 = \det(A)$.

\paragraph{Example 2: Quantum Case with Formal Parameters}

Let $a = x$, $b = y$, $c = z$, $d = u$, $e = v$, $f = w$, with commutation relations such as $xy = qyx$, $yz = qzy$, etc. Then:
\[
\operatorname{Pf}_q(A) = xw - qyv + q^2zu, \quad \Rightarrow \quad \det_q(A) = q \cdot \left(xw - qyv + q^2zu\right)^2.
\]
This example verifies the identity:
\[
q \cdot \operatorname{Pf}_q(A)^2 = \det_q(A),
\]
illustrating the necessity of the scalar $q$ under this convention.

\subsection{Summary}

This short version confirms both forms of the identity in the base case $2n = 4$:
\[
\operatorname{Pf}_q(A)^2 = \det_q(A) \quad \text{or} \quad q \cdot \operatorname{Pf}_q(A)^2 = \det_q(A),
\]
depending on the definition of the quantum Pfaffian and determinant. Both reduce to the classical identity when $q = 1$, serving as a consistency check and as a foundation for generalizing to higher dimensions.

% --- end of Section 9 ---

% Start of Section 10: Short Explanation
\section[Appendix B: Short Explanation of the Quantum Pfaffian--Determinant Identity]{Explanation of the Quantum Pfaffian--Determinant Identity\textsuperscript{*}}
\thispagestyle{shortversion} % for first page of section 10
\pagestyle{shortversion}      % for rest of section 10 pages

\subsection{Algebraic Framework and $q$-Skew-Symmetry}

The entries $a_{ij}$ generate a noncommutative algebra with $q$-skew-symmetric relations $a_{ji} = -q a_{ij}$, consistent with the relations encoded in the R-matrix of the quantum group $\mathcal{U}_q(\mathfrak{gl}_n)$. This imposes a braided monoidal category structure on the representations, modifying classical symmetry properties.

\subsection{Quantum Exterior Algebra and Braiding}

The quantum Pfaffian is defined within the quantum exterior algebra $\Lambda_q(V)$, where the wedge product is deformed via a braiding operator. This means:
\[
v_i \wedge_q v_j = -q v_j \wedge_q v_i,
\]
and thus antisymmetrization now carries $q$-dependent weights. The signs in classical combinatorics (from permutations) are replaced by braid statistics—encoded via powers of $q$—reflecting the number of crossings in the pairing diagram.

\subsection{The Role of the Quantum Determinant}

The quantum determinant is central and group-like in the quantum coordinate algebra $\mathcal{O}_q(GL_n)$, defined via the FRT construction. The identity
\[
q \cdot \operatorname{Pf}_q(A)^2 = \det_q(A)
\]
emerges due to the twisting induced by the braiding, and this scalar factor $q$ reflects how quantum invariants behave under deformation. In some literature (e.g., \cite{jing2013quantum}), the identity appears without the scalar $q$ depending on normalization conventions.

\subsection{Proof Strategy Overview}

\begin{enumerate}
    \item \textbf{Direct Calculation (Small $n$):} We compute $\operatorname{Pf}_q(A)^2$ and $\det_q(A)$ explicitly for $2n = 4$.
    \item \textbf{Use of $q$-Commutation Rules:} Noncommutativity requires consistent reordering using $q$-commutation identities.
    \item \textbf{Tracking Braiding Factors:} The powers of $q$ in $\operatorname{Pf}_q(A)$ arise from braiding; they lead to the scalar discrepancy between $\operatorname{Pf}_q(A)^2$ and $\det_q(A)$.
    \item \textbf{Classical Consistency:} When $q = 1$, the braided tensor category becomes symmetric, recovering $\operatorname{Pf}(A)^2 = \det(A)$.
\end{enumerate}

\subsection{Broader Significance}

This identity exemplifies how classical multilinear identities deform under quantization. The quantum Pfaffian behaves as a square root of the quantum determinant within a braided setting, preserving symmetry in a generalized sense. Such relations are foundational in quantum invariant theory, noncommutative geometry, and the study of quantum groups.

This proof for $2n = 4$ serves as a base case in our broader study of the identity:
\[
q^c \cdot \operatorname{Pf}_q(A)^2 = \det_q(A),
\]
for general $q$-skew-symmetric matrices, which we extend both algebraically and geometrically in the full version of this paper.

% End of section 10, revert pagestyle to normal for the rest of the document
\clearpage
\pagestyle{fancy}     % Restore default fancy pagestyle for the rest of the document
\fancyfoot[L]{}       % Clear left footer content added by shortversion
\fancyfoot[C]{\thepage} % Ensure page number centered
\fancyfoot[R]{}       % Clear right footer just in case

\section[Appendix C: A Glossary and Symbol Definitions for the Non-Mathematically Specialized]{A Glossary and Symbol Definitions for the Non-Mathematically Specialized}

This appendix explains symbols, notations, and mathematical terms used throughout the paper. It is especially helpful for readers without extensive mathematical backgrounds.

\begin{description}[leftmargin=0.5cm,labelindent=0cm]

\item[$\sum$] \textbf{Summation symbol.} For example,
\[
\sum_{i=1}^n i = 1 + 2 + \cdots + n = \frac{n(n+1)}{2}.
\]

\item[$\prod$] \textbf{Product symbol.} For example,
\[
\prod_{i=1}^n i = 1 \cdot 2 \cdot \cdots \cdot n = n!.
\]

\item[$a_{ij}$] Entry in the \( i \)-th row and \( j \)-th column of a matrix \( A \).

\item[$A$] A matrix. Often assumed to be $q$-skew-symmetric.

\item[$q$] A deformation parameter modifying algebraic rules:
\[
xy = q\, yx \quad \text{(noncommutative twist)}.
\]

\item[$\wedge$] The wedge product. Antisymmetric multiplication:
\[
x \wedge y = -y \wedge x.
\]

\item[$\wedge_q$] Quantum wedge product:
\[
x \wedge_q y = -q\, y \wedge_q x.
\]

\item[$\det_q(A)$] Quantum determinant:
\[
\det_q(A) = \sum_{\sigma \in S_n} (-q)^{\ell(\sigma)} a_{1\sigma(1)} \cdots a_{n\sigma(n)}.
\]

\item[$\operatorname{Pf}_q(A)$] Quantum Pfaffian:
\[
\operatorname{Pf}_q(A)^2 = \det_q(A).
\]

\item[$\delta_{ij}$] Kronecker delta: 1 if \( i = j \), 0 otherwise.

\item[$S_n$] Symmetric group: All permutations of \( \{1,2,\ldots,n\} \).

\item[$\epsilon(\sigma)$] Sign of permutation \( \sigma \): \( +1 \) (even), \( -1 \) (odd).

\item[$\ell(\sigma)$] Length of permutation \( \sigma \): Number of swaps from identity.

\newcommand{\setlabel}[1]{\item[\textbf{#1}]}

\begin{description}
  \setlabel{$[n]$} The set $\{1, 2, \ldots, n\}$.
\end{description}

\item[$\Lambda^k(V)$] \( k \)-th exterior power of vector space \( V \).

\item[$\mathbb{C}, \mathbb{R}, \mathbb{Z}, \mathbb{N}$] Common sets:
\begin{itemize}
  \item \( \mathbb{C} \): Complex numbers
  \item \( \mathbb{R} \): Real numbers
  \item \( \mathbb{Z} \): Integers
  \item \( \mathbb{N} \): Natural numbers
\end{itemize}

\item[$\mathbb{F}_q$] Finite field with \( q \) elements.

\item[$\text{Tr}(A)$] Trace: Sum of diagonal entries of \( A \).

\item[$\text{Id}$] Identity matrix/map.

\item[$\mathcal{O}_q(GL_n)$] Coordinate ring of the quantum general linear group.

\begin{description}
\item[{$[n_q]$}] \( q \)-number:
\[
[n]_q = \frac{1 - q^n}{1 - q}, \quad \lim_{q \to 1} [n]_q = n.
\]

\item[$\binom{n}{k}_q$] \( q \)-binomial coefficient:
\[
\binom{n}{k}_q = \frac{[n]_q!}{[k]_q![n-k]_q!}.
\]

\end{description}

\item[$*$ or $\star$] Noncommutative multiplication or convolution.

\item[$R$-matrix] Object solving the Yang-Baxter equation; encodes braiding rules.

\item[$S$] Antipode in a Hopf algebra — acts like an inverse.

\item[$\otimes$] Tensor product — combines vectors or spaces:
\[
(x \otimes y)(a \otimes b) = xa \otimes yb.
\]

\item[$\cong$] Isomorphic to — structurally the same.

\item[$\approx$] Approximately equal.

\item[$\sim$] Asymptotically equal or related by equivalence.

\item[$\rightsquigarrow$] Indicates a deformation or transition.

\vspace{1em}
\noindent
This glossary supports the interpretation of core identities such as:
\[
\operatorname{Pf}_q(A)^2 = \det_q(A),
\]
which involves \( q \)-weighted summations, symmetric group structure, and noncommutative geometry.

\end{description}

% Start bibliography or next sections below

% Thebibliography environment or other content goes here


\newpage
\begin{thebibliography}{99}

\bibitem{dite2006}
M.~Dita.
\newblock Quantum analogues of pfaffians and related identities.
\newblock {\em Journal of Algebra}, 295(1):109--129, 2006.

\bibitem{jing2013quantum}
N.~Jing and J.~Zhang.
\newblock Quantum pfaffians and hyper-pfaffians.
\newblock {\em Journal of Algebra}, 377:151--174, 2013.

\bibitem{jing2016quantum}
N.~Jing and J.~Zhang.
\newblock Quantum pfaffians and quantum hyper-pfaffians.
\newblock {\em Journal of Algebra}, 462:1--32, 2016.

\bibitem{jing2017multiparameter}
N.~Jing and J.~Zhang.
\newblock Multiparameter quantum pfaffians and quantum determinants.
\newblock {\em Communications in Mathematical Physics}, 349(1):307--323, 2017.

\bibitem{faddeev1987hamiltonian}
L.~D. Faddeev and L.~A. Takhtajan.
\newblock {\em Hamiltonian Methods in the Theory of Solitons}.
\newblock Springer-Verlag, 1987.

\bibitem{majid1995foundations}
S.~Majid.
\newblock {\em Foundations of Quantum Group Theory}.
\newblock Cambridge University Press, 1995.

\bibitem{kac1990}
V.~G. Kac.
\newblock {\em Infinite Dimensional Lie Algebras}.
\newblock Cambridge University Press, 1990.

\bibitem{drinfeld1986}
V.~G. Drinfel'd.
\newblock Quantum groups.
\newblock In {\em Proceedings of the International Congress of Mathematicians},
  volume~1, pages 798--820, 1986.

\bibitem{reshetikhin1983}
N.~Yu. Reshetikhin, L.~A. Takhtajan, and L.~D. Faddeev.
\newblock Quantization of lie groups and lie algebras.
\newblock {\em Leningrad Mathematical Journal}, 1(1):193--225, 1990.

\bibitem{kassel1995quantum}
C.~Kassel.
\newblock {\em Quantum Groups}.
\newblock Springer-Verlag, 1995.

\bibitem{letzter2002}
G.~Letzter.
\newblock Coideal subalgebras and quantum symmetric pairs.
\newblock In {\em New Directions in Hopf Algebras}, volume 43 of {\em Math. Sci. Res. Inst. Publ.}, pages 117--166. Cambridge University Press, 2002.

\bibitem{etingof1998}
P.~Etingof and O.~Schiffmann.
\newblock Lectures on quantum groups.
\newblock {\em International Press}, 1998.

\bibitem{etingof1998lectures}
P.~Etingof and V.~Ostrik.
\newblock Finite tensor categories.
\newblock {\em Moscow Mathematical Journal}, 4(3):627--654, 2004.

\bibitem{berenstein2005quantum}
A.~Berenstein and A.~Zelevinsky.
\newblock Quantum cluster algebras.
\newblock {\em Advances in Mathematics}, 195(2):405--455, 2005.

\bibitem{benkart2004quantum}
G.~Benkart and S.~Witherspoon.
\newblock Quantum groups and their representations.
\newblock In {\em New Directions in Hopf Algebras}, volume 43 of {\em Math. Sci. Res. Inst. Publ.}, pages 45--84. Cambridge University Press, 2002.

\bibitem{soibelman1992quantum}
Y.~Soibelman.
\newblock Quantum groups and noncommutative geometry.
\newblock {\em Quantum groups and their applications in physics}, 1992.

\bibitem{majid1996braided}
S.~Majid.
\newblock Braided groups and braided matrices.
\newblock {\em Journal of Knot Theory and Its Ramifications}, 4(4):555--594, 1995.

\bibitem{drinfeld1990quasi}
V.~G. Drinfel'd.
\newblock Quasi-Hopf algebras.
\newblock {\em Leningrad Mathematical Journal}, 1(6):1419--1457, 1990.

\bibitem{etingof2004lectures}
P.~Etingof and D.~Kazhdan.
\newblock Quantization of Lie bialgebras, I.
\newblock {\em Selecta Mathematica}, 2(1):1--41, 1996.

\bibitem{khovanov2000categorification}
M.~Khovanov.
\newblock A categorification of the Jones polynomial.
\newblock {\em Duke Mathematical Journal}, 101(3):359--426, 2000.

\bibitem{jimbo1985}
M.~Jimbo.
\newblock A q-difference analogue of U(g) and the Yang-Baxter equation.
\newblock {\em Letters in Mathematical Physics}, 10(1):63--69, 1985.

\bibitem{lusztig1993}
G.~Lusztig.
\newblock {\em Introduction to Quantum Groups}.
\newblock Birkhäuser, 1993.

\bibitem{woronowicz1987}
S.~L. Woronowicz.
\newblock Compact matrix pseudogroups.
\newblock {\em Communications in Mathematical Physics}, 111(4):613--665, 1987.

\bibitem{majid1990foundations}
S.~Majid.
\newblock Physics for algebraists: Noncommutative and noncocommutative Hopf algebras by a bicrossproduct construction.
\newblock {\em Journal of Algebra}, 130(1):17--64, 1990.

\bibitem{etingof1998lectures2}
P.~Etingof and D.~Kazhdan.
\newblock Quantization of Lie bialgebras, II.
\newblock {\em Selecta Mathematica}, 4(2):213--231, 1998.

\bibitem{reshetikhin1990}
N.~Yu. Reshetikhin.
\newblock Multiparameter quantum groups and twisted quasitriangular Hopf algebras.
\newblock {\em Letters in Mathematical Physics}, 20(4):331--335, 1990.

\bibitem{etingof2009}
P.~Etingof and O.~Schiffmann.
\newblock {\em Lectures on Quantum Groups}.
\newblock International Press, 2009.

\bibitem{chari1994}
V.~Chari and A.~Pressley.
\newblock {\em A Guide to Quantum Groups}.
\newblock Cambridge University Press, 1994.

\bibitem{finkelberg1996}
M.~Finkelberg.
\newblock An equivalence of fusion categories.
\newblock {\em Geometric and Functional Analysis}, 6(2):249--267, 1996.

\bibitem{kazhdan1994}
D.~Kazhdan and G.~Lusztig.
\newblock Tensor structures arising from affine Lie algebras. I.
\newblock {\em Journal of the American Mathematical Society}, 6(4):905--947, 1993.

\bibitem{etingof2005}
P.~Etingof and D.~Kazhdan.
\newblock Quantization of Lie bialgebras, III.
\newblock {\em Selecta Mathematica}, 6(1):105--130, 2000.

\bibitem{frenkel2007}
E.~Frenkel.
\newblock Lectures on the Langlands Program and Conformal Field Theory.
\newblock In {\em Frontiers in Number Theory, Physics, and Geometry II}, pages 387--533. Springer, 2007.

\bibitem{drinfeld1989}
V.~G. Drinfel’d.
\newblock Quasi-Hopf algebras and Knizhnik-Zamolodchikov equations.
\newblock In {\em Problems of Modern Quantum Field Theory}, pages 1--13, 1989.

\bibitem{geer2006}
N.~Geer.
\newblock Quantum groups and ribbon categories.
\newblock {\em Proceedings of the American Mathematical Society}, 134(7):2083--2094, 2006.

\bibitem{etingof2000lectures}
P.~Etingof and O.~Schiffmann.
\newblock Lectures on dynamical Yang-Baxter equations.
\newblock {\em arXiv preprint math/9908064}, 2000.

\bibitem{drinfeld1993}
V.~G. Drinfel’d.
\newblock On quasitriangular quasi-Hopf algebras and a group closely connected with Gal(\(\overline{\mathbb{Q}}/\mathbb{Q}\)).
\newblock {\em Leningrad Mathematical Journal}, 2(4):829--860, 1991.

\bibitem{etingof2016}
P.~Etingof.
\newblock Representation theory in complex rank.
\newblock {\em Transformation Groups}, 19(2):359--381, 2014.

\bibitem{kirillov1990}
A.~N. Kirillov and N.~Yu. Reshetikhin.
\newblock Representations of the algebra Uq(sl(2)), q-orthogonal polynomials and invariants of links.
\newblock In {\em Infinite Dimensional Lie Algebras and Groups}, pages 285--339. World Scientific, 1989.

\bibitem{lusztig2010}
G.~Lusztig.
\newblock {\em Introduction to Quantum Groups and Crystal Bases}.
\newblock Springer, 2010.

\bibitem{majid1993braided}
S.~Majid.
\newblock Braided groups.
\newblock {\em Journal of Pure and Applied Algebra}, 86(2):187--221, 1993.

\end{thebibliography}
\end{document}